\documentclass[11pt]{article}
\usepackage{amsmath,amssymb}
\usepackage{graphicx}
\usepackage{wrapfig}
\usepackage{hyperref}

\newcommand{\qed}{{\unskip\nobreak\hfil\penalty50\hskip2em\vadjust{}
       \nobreak\hfil$\Box$\parfillskip=0pt\finalhyphendemerits=0\par}}

\newtheorem{theorem}{Theorem}[section]
\newtheorem{lemma}{Lemma}[section]
\newtheorem{definition}{Definition}[section]
\newtheorem{cor}{Corollary}
\newtheorem{prop}{Proposition}

\newcommand{\eps}{\epsilon}

\newcommand{{\Z}}{\mathbb Z}
\newcommand{\R}{\mathbb R}
\newcommand{\N}{\mathbb N}

\newcommand{\dist}{\vert \vert}

\renewcommand{\P} {{\mathcal P}}

\newcommand{\cir}{\Gamma_0}
\newcommand{\reg}{{\rm RG} \big( \Gamma_0 \big)}
\newcommand{\reggam}{{\rm RG} \big( \Gamma \big)}
\newcommand{\conv}{{\rm conv}\big( \Gamma_0 \big)}
\newcommand{\delconv}{\partial {\rm conv}\big( \Gamma_0 \big)}
\newcommand{\mfl}{{\rm MFL}\big( \Gamma_0 \big)}
\newcommand{\mlr}{{\rm MLR}\big( \Gamma_0 \big)}
\newcommand{\mlrs}{{\rm MLRF}\big( \Gamma_0 \big)}

\newcommand{\vcir}{V\big( \Gamma_0 \big)}

\newcommand{\acon}{{\rm AREA}_{\bo{0},n^2}}

\newcommand{\intg}{{\rm INT}\big( \Gamma_0 \big)}

\newcommand{\globdis}{{\rm GD}\big( \Gamma_0 \big)}
\newcommand{\mc}{\mathcal}
\newcommand{\bo}{\mathbf}
\newcommand{\argu}{{\rm arg}}

\newcommand{\axy}{A_{\bo{x},\bo{y}}}
\newcommand{\axye}{E\big(A_{\bo{x},\bo{y}}\big)}

\newcommand{\xo}{\bo{x_0}}
\newcommand{\yo}{\bo{y_0}}

\newcommand{\gacc}{{\rm GAC}}
\newcommand{\zoz}{\{0,1 \}^{E(\Z^2)}}

\newcommand{\ang}{\angle}
\newcommand{\qzero}{q_0}
\newcommand{\clu}[1]{C_{\pi/2 - \qzero}^F \big( #1 \big)}
\newcommand{\clum}[1]{C_{\pi/2 - \qzero}^B \big( #1 \big)}

\newcommand{\mar}{\theta_{{\rm RG}_{\qzero,c_0}}^{\rm MAX}\big(\cir\big)}
\newcommand{\margam}{\theta_{{\rm RG}}^{\rm MAX}\big(\Gamma\big)}

\newcommand{\ccone}{c_1}
\newcommand{\cctwo}{C_1}

\newcommand{\cfour}{s_3(n)}
\newcommand{\cfive}{s_2(n)}
\newcommand{\csix}{s_1(n)}

\newcommand{\cstar}{C^*}

\newcommand{\cona}{c'}
\newcommand{\conb}{c''}
\newcommand{\cgac}{2}

\newcommand{\clemgac}{C_2}

\newcommand{\crwm}{C}

\newcommand{\ctheta}{\chi}
\newcommand{\qone}{s_1(n)}
\newcommand{\qtwo}{s_2(n)}
\newcommand{\qthr}{s_3(n)}
\newcommand{\cnew}{c_3}
\newcommand{\epsone}{\epsilon_1}

\newcommand{\res}{{\rm RES}}
\newcommand{\lr}[1]{{\rm LR} \big( \bo{#1} , \gamsw \big)}
\newcommand{\gamsw}{\Gamma_{\rm(sw)}}
\newcommand{\rggamsw}{{\rm RG}\big(\Gamma_{\rm (sw)}\big)}
\newcommand{\rggamswcl}{{\rm RG}\big(\overline\Gamma_{\rm (sw)}\big)}
\newcommand{\rggamswclex}{{\rm RG}_{\qzero/2,c_0/2}\big(\overline\Gamma_{\rm (sw)}\big)}
\newcommand{\csw}{\bo{c_{sw}}}
\newcommand{\sw}[1]{{\rm SW}\big( #1 \big)}
\newcommand{\areatrap}{{\rm AREA}_{n^2}^{\rm sw}}
\newcommand{\smtu}{{\rm MBT}}
\newcommand{\smtup}{\overline\smtu_0}
\newcommand{\unfav}{{\rm UNFAV}}
\newcommand{\sid}{{\rm SID}}
\newcommand{\circl}{\overline\Gamma_0}
\newcommand{\regcl}{{\rm RG}\big( \circl \big)}
\newcommand{\gamswcl}{\overline\Gamma_{\rm (sw)}}

\newcommand{\centre}{{\rm cen}}

\newcommand{\areacon}{\big\vert \intg \big\vert \geq n^2}

\newcommand{\wulff}{\mathcal{W}_\beta}
\newcommand{\outfluc}{{\rm OUTFLUC}}
\newcommand{\marcl}{\theta_{{\rm RG}_{\qzero/2,c_0/2}}^{\rm MAX}\big( \overline\Gamma_{\rm (sw)} \big)}
\newcommand{\marclo}{\theta_{{\rm RG}_{\qzero,c_0}}^{\rm MAX}\big( \overline\Gamma_0 \big)}
\newcommand{\secajn}{A_j^n}
\newcommand{\thetan}{\theta_n}
\newcommand{\mn}{m_n}
\newcommand{\functht}{u(\chi)}
\newcommand{\csid}{\epsilon_0}
\newcommand{\epso}{\epsilon_0}
\newcommand{\cluh}[1]{C_{\pi/2 - \qzero/2}^F \big( #1 \big)}
\newcommand{\clumh}[1]{C_{\pi/2 - \qzero/2}^B \big( #1 \big)} 
\newcommand{\sentier}{\tau}

\def\build#1_#2^#3{\mathrel{ \mathop{\kern 0pt#1}\limits_{#2}^{#3}}}

\def\fff#1{&{{\pageref{#1}}}\cr}
\def\hfff#1{\label{#1}}

\begin{document}
\title{Phase separation in random cluster models II: \\ the droplet at equilibrium, and local deviation lower bounds}
\author{Alan Hammond\thanks{Department of Statistics, University of Oxford. This work was undertaken during a visit to the Theory Group at Microscoft Research in Redmond, WA, and at {\'E}cole Normale Sup{\'e}rieure in Paris.}} 
 \maketitle

\begin{abstract}
We study the droplet that results from conditioning the planar subcritical Fortuin-Kasteleyn random
cluster model on the presence of an open circuit $\cir$ encircling the origin and enclosing an
area of at least (or exactly) $n^2$. 
We consider local deviation of the droplet boundary, measured in a radial sense by the maximum local roughness, $\mlr$, this being the maximum distance from a point in the circuit $\cir$ to the boundary $\delconv$ of the circuit's convex hull; and in a longitudinal sense by what we term {\it maximum facet length}, namely, the length of the longest line segment of which the polygon $\delconv$ is formed. 
We prove that there exists a constant $c > 0$
such that the conditional probability that the normalised quantity 
$n^{-1/3}\big( \log n \big)^{-2/3} \mlr$ 
exceeds $c$ 
tends to $1$ in the high $n$-limit; and that the same statement holds for  
$n^{-2/3}\big( \log n \big)^{-1/3}  \mfl$. 
To obtain these bounds, we exhibit the random cluster measure conditional on the presence of an open circuit trapping high area as the invariant measure of a Markov chain that resamples sections of the circuit boundary. We analyse the chain at equilibrium to prove the local roughness lower bounds. 
Alongside complementary upper bounds provided in \cite{hammondone}, the fluctuations $\mlr$ and $\mfl$ are determined up to a constant factor.
\end{abstract}

\setlength{\baselineskip}{16pt}

\newpage
\tableofcontents
\newpage
\begin{subsection}{Glossary of notation}
A fair amount of notation will be required in the proofs appearing in this article. For the reader's convenience, we begin by listing much of the notation, and provide a summarizing phrase for each item, as well as the page number at which the concept is introduced.

\bigskip
\def\qq{&}

\begin{center}
\halign{
#\quad\hfill&#\quad\hfill&\quad\hfill#\cr
$\Gamma$ \qq a generic circuit \fff{gencir}  
${\rm INT}(\Gamma)$ \qq the region enclosed by $\Gamma$ \fff{intcir}
$\cir$ \qq  the outermost open circuit enclosing $\bo{0}$ \fff{outcir}
$\conv$ \qq the convex hull of $\cir$ \fff{conv}
$\mlr$ \qq  maximum local roughness \fff{mlr}    
$\mfl$ \qq  maximum facet length \fff{mfl}    
$T_{\bo{0},\bo{x},\bo{y}}$ \qq the triangle with vertices $\bo{0}$, $\bo{x}$ and $\bo{y}$ \fff{txy}
$\axy$ \qq  the sector with apex $\bo{0}$ bounded by $\bo{x}$ and $\bo{y}$ \fff{axy}   
$W_{\bo{v},c}$ \qq  the cone about $\bo{v}$ with apex $\bo{0}$ and aperture $2c$ \fff{wvc}
$\wulff$ \qq  the unit-area Wulff shape \fff{wulff}   
$\globdis$ \qq  global distortion (from the Wulff shape) \fff{globdis}
$\centre(\Gamma)$ \qq  the centre of a circuit \fff{centre}
$\acon$ \qq the event of $n^2$-area capture by a centred circuit \fff{acon}
$C^F_{\pi/2 - q} \big( \bo{v} \big)$ 
\qq the $q$-forward cone with apex $\bo{v}$ \fff{forback}
$C^B_{\pi/2 - q} \big( \bo{v} \big)$ 
\qq the $q$-backward cone  with apex $\bo{v}$ \fff{forback}
 $\reg$ \qq  the set of circuit regeneration sites \fff{reg}   
 $\regcl$ \qq  the set of cluster regeneration sites \fff{regcl}   
 $\marclo$ \qq  the maximal angular separation of cluster regeneration sites \fff{marclo}   
$\sw{\Gamma}$ \qq  the southwest corner of a circuit \fff{sw}
$\areatrap$ \qq the event of $n^2$-area capture by a southwest-centred circuit \fff{areatrap}  
$\gamsw$ \qq  the outermost open circuit as defined by southwest-centering \fff{gamsw}  
$\gamswcl$ \qq  the outermost open cluster as defined by southwest-centering \fff{gamswcl}  
$\gamma_{\bo{x},\bo{y}}$ \qq  the outermost open path in $A_{\bo{x},\bo{y}}$ from $\bo{x}$ to $\bo{y}$ \fff{gamoop}   
 $\overline\gamma_{\bo{x},\bo{y}}$ \qq  the open cluster of  $\bo{x}$ and $\bo{y}$ in $A_{\bo{x},\bo{y}}$ \fff{ogamoop}
$I_{\bo{x},\bo{y}} \big( \gamma_{\bo{x},\bo{y}} \big)$ \qq the bounded component of $\axy \setminus \gamma_{\bo{x},\bo{y}}$ \fff{ixy} 
 $\res_j$ \qq  the resampling procedure in the $j$-th sector $A_j^n$ \fff{resj}  
$\smtu$ \qq   the index-set of sectors with moderate boundary turning \fff{mbt}  
$\functht$ \qq a function of $\ctheta$ for which $\functht = o\big(\ctheta^{1/2}\big)$ \fff{functht}
$\unfav$ \qq  the index-set of unfavourable sectors (having low local roughness) \fff{unfav}   
 $\gacc \big( \bo{x},\bo{y} ,\epsilon\big)$ \qq  configurations in $\axy$ realizing $\epsilon$-good area capture \fff{gac}  
 $\sid \big( \bo{x},\bo{y} ,\epsilon\big)$ \qq $\epsilon$-significant inward deviation (of a connected set in $\axy$ from $\bo{x}$ to $\bo{y}$) \fff{sid}
 $\res$ \qq  the complete resampling procedure \fff{res}  
}\end{center}

\end{subsection}

\begin{section}{Introduction}

The theory of phase separation is concerned with the geometry of the random boundary between 
populations of distinct spins in a statistical mechanical model such as percolation, the Potts model or the random cluster model. For example, if the two-dimensional Ising model at supercritical inverse temperature $\beta > \beta_c$
in a large box with negative boundary conditions is conditioned by the presence of a significant excess of plus signs, then those excess signs typically gather together in a single droplet having the opposite magnetisation to its exterior. The object of study of phase separation is then the droplet boundary. As explained in \cite{alexcube} and \cite{hammondone}, a close relative of this problem is that of the behaviour of the circuit that arises by conditioning a subcritical random cluster model on the presence of a circuit encircling the origin and trapping a high area.

Wulff \cite{wulff} proposed that the profile of such constrained circuits would macroscopically resemble a dilation of an isoperimetrically optimal curve that now bears his name. For the Ising problem,
this claim was first verified rigorously in \cite{dks} at low temperature, 
the derivation being extended up to the critical temperature by \cite{ioffeschonmann}. Fluctuations from this profile have been classified into global or long-wave effects, measured by the deviation (in the Hausdorff metric) of the convex hull of the circuit from an optimally placed dilate of the Wulff shape. Local fluctuations have been measured by the inward deviation of the circuit from the boundary of its convex hull. 

In this article, we prove lower bounds on the local fluctuations of  the outermost open circuit in a subcritical random cluster model, when this circuit is conditioned to entrap a large area.  

We recall the definition of the random cluster model.
\begin{definition}
For $\Lambda \subseteq \Z^2$, let $E(\Lambda)$ denote the set of nearest-neighbour edges whose endpoints lie in $\Lambda$
and write $\partial_{\rm int} \big( \Lambda \big)$ for the interior vertex boundary of $\Lambda$, namely, the subset of $\Lambda$ each of whose elements is an endpoint of some element of 
 $E(\Lambda)^c$. Fix a choice of $\Lambda \subseteq \Z^2$ that is finite. 
The free random cluster model on $\Lambda$ with parameters $p \in [0,1]$ and $q > 0$ on $\Lambda$ 
is the probability space over $\eta: E(\Lambda) \to \{0,1\}$ with measure
$$
\phi_{p,q}^f(\eta) = \frac{1}{Z_{p,q}}   p^{\sum_e \eta(e)} 
\big( 1 - p \big)^{\sum_e (1 - \eta(e))} q^{k(\eta)},
$$
where $k(\eta)$ denotes the number of connected components in the subgraph of 
$\big(\Lambda,E(\Lambda)\big)$ containing all vertices and all edges $e$ such that $\eta(e) = 1$. (The constant $Z_{p,q}$ is a normalization.) The wired random cluster model 
$\phi_{p,q}^w$ is defined similarly, with $k(\eta)$ now denoting the number of such connected components none of whose edges touch $\partial_{\rm int} \big( \Lambda \big)$.

For parameter choices $p \in [0,1]$ and $q \geq 1$, 
either type of random cluster measure $\P$ satisfies the FKG inequality: suppose that 
$f,g:  \{0,1\}^{E(\Lambda)} \to \R$
are increasing functions with respect to the natural partial order on $\{0,1\}^{E(\Lambda)}$. 
Then $\mathbb{E}_\P \big( fg \big) \geq \mathbb{E}_\P \big( f \big) \mathbb{E}_\P \big( g \big)$,
where $\mathbb{E}_\P$ denotes expectation with respect to $\P$.

Consequently, we define the infinite-volume free and wired random cluster measures $\P^f$ and $\P^w$ as limits of the finite-volume counterparts taken along any increasing sequence of finite sets $\Lambda \uparrow \Z^2$. The measures $\P^f$ and $\P^w$ are defined on the space of functions $\eta: E(\Z^2) \to \{0,1\}$ with the product $\sigma$-algebra. In a realization $\eta$, the edges $e \in E(\Z^2)$
 such that $\eta(e) = 1$ are called open; the remainder are called closed.  A subset of $E(\Z^2)$ will be called open (or closed) if all of its elements are open (or closed).
We will record a realization in the form $\omega \in \zoz$, where 
the set of coordinates that are equal to $1$ under $\omega$ is the set of open edges under $\eta$.
Any $\omega \in \zoz$ will be called a configuration. For $\bo{x},\bo{y} \in \Z^d$, we write $\bo{x} \leftrightarrow \bo{y}$ to indicate that $\bo{x}$ and $\bo{y}$ lie in a common connected component of open edges.  
 
Set $\beta \in (0,\infty)$ according to $p = 1 - \exp\{ - 2 \beta \}$. 
In this way, the infinite volume measures are parameterized by $\P_{\beta,q}^w$ and $\P_{\beta,q}^f$ with $\beta > 0$ and $q \geq 1$.
For any $q \geq 1$, $\P^w_{\beta,q} = \P^f_{\beta,q}$ for all but at most countably many values of $\beta$ \cite{grimmett}. 
We may thus define
$$
\beta_c^1 = \inf \big\{  \beta > 0: \P^*_{\beta,q} \big(  0 \leftrightarrow \infty \big) > 0 \big\}
$$
obtaining the same value whether we choose $* = w$ or $* = f$.
\end{definition}
There is a unique random cluster model for each subcritical $\beta < \beta_c^1$ \cite{grimmett}, that we will denote by $\P_{\beta,q}$.
\begin{definition}
Let $\hat{\beta}_c$ denote the supremum over $\beta > 0$ such that the following holds: letting $\Lambda = \big\{ -N,\ldots,N \big\}^2$, there exist constants $C > c > 0$ such that, for any $N$,
$$
\P^w_{\beta,q} \Big( \bo{0} \leftrightarrow \Z^2 \setminus  \Lambda_N \Big) \leq C \exp \big\{ - c N  \big\}.
$$
\end{definition}
In the two-dimensional case that is the subject of this article, 
it has been established 
that $\beta_c^1 = \hat{\beta}_c$ for $q=1$, $q=2$  
and for $q$ sufficiently high, by \cite{alexmix}, and respectively \cite{ab}, \cite{abf} and \cite{lmmrs}.
A recent advance \cite{beffaraduminilcopin} showed that, on the square lattice, in fact, $\beta_c^1 = \hat{\beta}_c$ holds for all $q \geq 1$. The common value, which is $2^{-1} \log \big( 1 + \sqrt{q} \big)$, we will denote by $\beta_c$.

The droplet boundary is now defined:
\begin{definition}\label{defcirsur}
A circuit \hfff{gencir} $\Gamma$ is a nearest-neighbour path in $\Z^2$ whose endpoint coincides with its start point, but for which no other vertex is visited twice.
We set $E(\Gamma)$ equal to the set of nearest-neighbour edges between successive elements of $\Gamma$. 
For notational convenience, when we write $\Gamma$, we refer to the closed subset of $\R^2$ given by the union of the topologically closed intervals corresponding to the elements of $E(\Gamma)$. 
We set $V(\Gamma) = \Gamma \cap \Z^2$.

Let $\omega \in \zoz$. A circuit $\Gamma$ is called open if $E(\Gamma)$ is open. 
For any circuit $\Gamma$, we write \hfff{intcir} ${\rm INT} \big( \Gamma \big) \subseteq \R^2$ for the bounded component of $\R^2 \setminus \Gamma$, that is, for the set of points in $\R^2$ enclosed by $\Gamma$. The notation $\vert \cdot \vert$, when applied to subsets of $\R^2$, will denote Lebesgue measure. 

An open circuit $\Gamma$ is called outermost if any open circuit $\Gamma'$ satisfying 
${\rm INT} \big( \Gamma \big) \subseteq {\rm INT} \big( \Gamma' \big)$ is equal to $\Gamma$. Note that, if
a point $\bo{z} \in \R^2$ is enclosed by a positive but finite number of open circuits in a configuration 
$\omega \in \zoz$ , it
is enclosed by a unique outermost open circuit.

We write \hfff{outcir} $\cir$ for the outermost open circuit $\Gamma$ for which $\bo{0} \in {\rm INT}(\Gamma)$, taking $\cir = \emptyset$ if no such circuit exists.
\end{definition}
\noindent{\bf Remark.} 
Under any subcritical random cluster measure $P = \P_{\beta,q}$, with $\beta < \beta_c$  and $q \geq 1$, there is an exponential decay in distance for the probability that two points lie in the same open cluster. (See Theorem $A$ of \cite{civ}.) As such, $P$-a.s., no point in $\R^2$ is surrounded by infinitely many open circuits, so that $\cir$ exists (and is non-empty) whenever $\bo{0}$ is surrounded by an open circuit. 

Our object of study is the subcritical random cluster model given the event $\areacon$, with $n \in \N$ high. The resulting circuit $\cir$ has a local geometry given by the Gaussian fluctuations of a subcritical path conditioned between two points. Globally, however, this geometry is constrained by curvature. 
To focus on the interplay between these local fluctuation and global curvature effects, it is natural to consider the curve near a given point, rescaled by $n^{2/3}$ longitudinally and $n^{1/3}$ radially. 
We prove lower bounds on the greatest deviations present in the circuit on these scales. 
The two notions of maximum fluctuation are now given.
\begin{definition}
We write \hfff{conv} $\conv$ for the convex hull of $V\big( \cir \big)$, 
$$
\conv = \Big\{ \bo{z} \in \R^2: \ \exists \bo{x},\bo{y} \in V\big( \cir \big) \, \, \textrm{for which $\bo{z} \in \big[ \bo{x}, \bo{y} \big]$} \Big\}.
$$ 
The maximum local roughness \hfff{mlr} $\mlr$ is defined to be the maximal distance of an element of $V(\cir)$
 to the boundary of the convex hull: that is
$$
 \mlr = \sup \Big\{ d \big( \bo{x} , \delconv \big): \bo{x} \in \vcir \Big\},
$$ 
where $d:\R^2 \to [0,\infty)$ denotes the Euclidean distance.
We will denote by \hfff{mfl} {\it maximum facet length}
  $\mfl$ the length of the longest line segment of which the polygon $\delconv$ is comprised.
  \end{definition}
It is proved in  \cite{uzunalex} that $P_p \big( \mlr \geq n^{1/3} (\log n)^{-2/3} \big\vert \areacon \big) \to 1$ as $n \to \infty$, with $P_p$, $p \in (0,1/2)$, denoting subcritical bond percolation on $\Z^2$. The power-law term here was expected to be sharp. 
For a broad class of subcritical models, local roughness was proved in \cite{alexcube} to be bounded above by $O \big(n^{1/3} (\log n)^{2/3} \big)$ 
in an $L^1$-sense, validating the sharpness of the power-law term for an averaged form of local roughness.  

Our central conclusion for radial local deviation is the following lower bound on maximum local roughness.
\begin{theorem}\label{thmmlrubd}
Let $P = \P_{\beta,q}$, with $\beta < \beta_c$ and $q \geq 1$. 
For any $\epsilon > 0$, there exists $\delta > 0$ such that 
$$
P \Big( \mlr < \delta n^{1/3} \big( \log n \big)^{2/3}  \Big\vert \areacon \Big)  \leq
  \exp \Big\{ - n^{\frac{1}{13} - \epsilon} \Big\}.
$$
\end{theorem}
Maximum facet length was not an object explicitly considered by K. Alexander and H. Uzun, but it plays a central role in our approach. The conclusion  is:
\begin{theorem}\label{thmmflubd}
Let $P = \P_{\beta,q}$, with $\beta < \beta_c$  and $q \geq 1$.
There exists $\epsilon > 0$ and a function $\phi: (0,\epsilon) \to (0,\infty)$ satisfying 
 $\lim_{c \downarrow 0}{\phi(c)} = \infty$ 
such that, for all $c \in (0,\epsilon)$, and for $n \in \N$ sufficiently high,
$$
P \Big( \mfl < c n^{2/3} \big( \log n \big)^{1/3}  \Big\vert \areacon \Big) \leq n^{-\phi(c)}.
$$
\end{theorem}
The article \cite{hammondone} proves upper bounds that complement Theorems \ref{thmmlrubd} and \ref{thmmflubd}.
Taken together, our conclusion is:
\begin{cor}\label{cormlrsum}
Let $P = \P_{\beta,q}$, with $\beta < \beta_c$  and $q \geq 1$. 
Then there exist constants $0 < c < C < \infty$ such that
$$
P \bigg(  c \leq \frac{\mlr}{n^{1/3} \big( \log n \big)^{2/3}} \leq C  \bigg\vert \areacon \bigg)
  \to 1,
\qquad \textrm{as $n \to \infty$,} 
$$
and 
$$
P \bigg(  c \leq \frac{\mfl}{n^{2/3} \big( \log n \big)^{1/3}} \leq C  \bigg\vert \areacon \bigg)
  \to 1, 
\qquad \textrm{as $n \to \infty$.} 
$$
\end{cor}
That is, the techniques of this paper and its counterpart \cite{hammondone}  are sufficient to derive the conjectured exponents for the power-laws in radial and longitudinal local deviation, and to identify and prove exponents for the logarithmic correction for these quantities. 

An important ingredient in obtaining these results is an understanding that the conditioned circuit is highly regular, with little backtracking from its overall progress. 
The third article \cite{hammondthr}
in the series presents such a result, on the renewal structure of the conditioned circuit.
We will state its main conclusion as Theorem \ref{thmmaxrg}. 
Equipped with this tool, it is straightforward to derive the results under a conditioning on a fixed area:
\begin{theorem}\label{thmfixedarea}
Theorems \ref{thmmlrubd}, \ref{thmmflubd} and Corollary  \ref{cormlrsum} are valid (with verbatim statements) under the conditional measure $P \big( \cdot \big\vert \vert \intg \vert = n^2 \big)$.  
\end{theorem}
The proof of Theorem \ref{thmfixedarea} appears in \cite{hammondone}.
\begin{subsection}{Techniques and relations}
To prove Theorem \ref{thmmlrubd}, we will introduce a random resampling procedure that updates the circuit under the conditional measure along a segment of length $n^{1/3}( \log n)^{2/3}$. We design this procedure so that it leaves invariant the conditional distribution. As such, the conditional measure may be viewed as the equilibrium of a process of updates, in which the circuit is updated along successive stretches of length 
 $n^{1/3}( \log n)^{2/3}$. We prove Theorem \ref{thmmlrubd} by analysing the behaviour of this updating process at equilibrium: the updating occasionally produces configurations realizing the lower bound on $\mlr$ that subsequent resamplings are not likely to undo. The technique of proof has interest from a probabilistic point of view, because it analyses a measure by regarding it as the equilibrium of a random process that is introduced for this purpose. 
Arguments in which a convenient Markov chain is defined and analysed have previously been used to prove such results as correlation inequalities in percolation \cite{vdbhk}: see 
\cite{haggstrom}
for a review.
Of course, the Metropolis algorithm is very commonly used to sample approximately a measure (as reviewed in Chapter 3 of \cite{lpw}). In our case, the process of resampling might also be adapted to sample the conditional measure.

Of the numerous definitions of boundary fluctuation, we believe that local deviation is of particular interest. We briefly explain why. As we have mentioned, the radial $n^{1/3}$ and longitudinal $n^{2/3}$ scalings for this deviation are also characterized as being that scale in which the competition between locally Gaussian fluctuation in the circuit and its global need to curve takes place on a roughly equal footing. Rescaling the interface by these radial and longitudinal factors, we obtain a random function similar to a Brownian excursion conditioned on the event of capturing an area of at least a certain random quantity, where the law of this quantity has been chosen so that the event of trapping this amount in area has typically a bounded probability of being satisfied. The distribution of the area under a unit-time Brownian excursion, known as the Airy distribution, has been computed in \cite{darling} and \cite{louchard}. This law is derived non-rigorously in \cite{majcomt} as the limiting law of the maximum height of a kinetic interface in the KPZ universality class \cite{kpz}. In seeking to understand the relationship between the random geometry of the droplet boundary in models of phase separation, and the functionals universal in models of interfaces subject to smoothing and random roughening, it thus seems very natural to make local deviation, and such functionals as its maximum, the object of attention. See \cite{hammondone} for further discussion of this theme. 
\end{subsection}
\begin{subsection}{The structure of the paper}
In Section \ref{secnott}, notation is established and some tools from \cite{hammondone} and \cite{hammondthr} are recalled. Among these is the regeneration structure theorem of \cite{hammondthr}. This result requires some definitions related to centering of the circuit under the conditional measure. Thus, after introducing general notation in Section \ref{secnot}, we explain in Section 
\ref{seccircen} the notation for centering used for the regeneration structure theorem, which is stated in Section \ref{secrrg}. This preliminary section also concludes with a large deviations' statement for the macroscopic profile of the circuit in Section \ref{secldgd} and the statement of upper bounds on maximum local deviation from \cite{hammondone}, since we will need these assertions in the proofs. 

In Section \ref{secmain}, we present the principal proof in this paper, that of Theorem \ref{thmmlrubd}. 
Given Theorem \ref{thmmlrubd}, and by means of surgical techniques presented in \cite{hammondone}, the proof of Theorem \ref{thmmflubd} is fairly straightforward. We explain why and present the argument at the end of Section \ref{secmain}. Section \ref{sectools} presents some technical proofs deferred from the proof of Theorem \ref{thmmlrubd}. 
The concluding Section \ref{secconcrem} discusses some questions posed by our approach, including how the decay rates in Theorems \ref{thmmlrubd} and \ref{thmmflubd} might be improved.
\\
\noindent{\bf Acknowledgments.} I would like to thank Kenneth Alexander, Dmitry Ioffe, Yuval Peres and Senya Shlosman for helpful discussions, and Horatio Boedihardjo and a referee for their detailed comments on the paper.
\end{subsection}
\end{section}
\begin{section}{Notation and tools}\label{secnott}
\begin{subsection}{Notation}\label{secnot}
\begin{definition}\label{defpathedge}
Elements of $\R^2$ will be denoted by boldface symbols. 
By a discrete path, we mean a list of elements of $\Z^2$, each being a nearest-neighbour of the preceding one, and without repetitions. 
In referring to a path, we mean a subset of $\R^2$ given by the union of the topologically closed edges formed from the set of consecutive pairs of vertices of some discrete path. 
(As such, a path is defined to be self-avoiding, including at its vertices.)
In a similar vein, any subset of $\R^2$ that is introduced as a connected set is understood to be a union of closed intervals $\big[ \bo{u},\bo{v} \big]$ corresponding to nearest-neighbour edges $(\bo{u},\bo{v})$.
For such a set $A$, we write $V(A) = A \cap \Z^2$ and $E(A)$ for the set of edges of which $A$
is comprised. 

For a general subset $A \subseteq \R^2$, we write $E(A)$ for the set of nearest-neighbour edges 
$(\bo{u},\bo{v}) \in E(\Z^2)$ such that $\big[ \bo{u},\bo{v} \big] \subseteq A$. (This is of course consistent with the preceding definition.) 
\end{definition}
\begin{definition}\label{defto}
For  $\bo{x},\bo{y} \in \Z^2$, $\bo{y} \not= \bo{x}$,
we write $\ell_{\bo{x},\bo{y}}$ for the planar line containing $\bo{x}$ and $\bo{y}$, and $\ell^+_{\bo{x},\bo{y}}$ for the semi-infinite line segment that contains $\bo{y}$ and has endpoint $\bo{x}$. We write $\big[ \bo{x},\bo{y} \big]$ for the line segment whose endpoints are $\bo{x}$ and $\bo{y}$. 
We write \hfff{txy} $T_{\bo{0},\bo{x},\bo{y}}$ for the
closed triangle with vertices $\bo{0}$, $\bo{x}$ and $\bo{y}$.
For $\bo{x},\bo{y} \in \R^2$, 
we write $\ang\big(\bo{x},\bo{y} \big) \in [0,\pi]$ for the angle between these two vectors.
Borrowing complex notation, we write $\argu\big(\bo{x}\big)$ for the argument of $\bo{x}$.
In many derivations, the cones, line segments and points in question all lie in a cone, rooted at the origin, whose aperture has angle strictly less than $2\pi$. As such, it is understood that $\argu$
denotes a continuous branch of the argument that is defined throughout the region under consideration.
\end{definition}
Sometimes we wish to specify a cone by a pair of boundary points, and sometimes by the argument-values of its boundary lines:
\begin{definition}
For $\bo{x},\bo{y} \in \Z^2$, $\argu(\bo{x}) < \argu(\bo{y})$, write \hfff{axy}
$$
 A_{\bo{x},\bo{y}} = \Big\{ \bo{z} \in \R^2: \argu\big( \bo{x} \big) \leq  \argu\big( \bo{z} \big) \leq  \argu\big( \bo{y} \big)  \Big\} \cup \big\{ \bo{0} \big\}.
$$ 
To specify a cone by the argument-values of its boundary lines, take $\bo{v} \in \Z^2$ and $c \in [0,\pi)$, and let \hfff{wvc}
\begin{equation}\label{wnot}
W_{\bo{v},c} = \Big\{ \bo{z} \in \R^2: \argu(\bo{v}) - c \leq  \argu(\bo{z}) \leq  \argu(\bo{v}) + c \Big\}
\cup \big\{ \bo{0} \big\}
\end{equation}
denote the cone of points whose angular displacement from $\bo{v}$ is at most $c$. 
Extending
this notation, for any $\bo{x} \in \Z^2$ and $c \in [0,\pi)$, we write 
$W_{\bo{v},c}\big( \bo{x} \big) = \bo{x} + W_{\bo{v},c}$.
We also write, for $\bo{x} \in \R^2$ and $c \in (0,2\pi)$,
$$
W_{\bo{x},c}^+  = \Big\{ \bo{z} \in \R^2: \argu(\bo{x})  \leq  \argu(\bo{z}) \leq  \argu(\bo{x}) + c \Big\}
\cup \big\{ \bo{0} \big\}
$$
and 
$$
W_{\bo{x},c}^-  = \Big\{ \bo{z} \in \R^2: \argu(\bo{x}) - c \leq  \argu(\bo{z}) \leq  \argu(\bo{x})  
\Big\}
\cup \big\{ \bo{0} \big\}.
$$
\end{definition}
\begin{definition}
Write $S^1 \subseteq \R^2$ for the boundary of the unit ball in the Euclidean metric. 
For $\bo{v} \in \R^2$, let $\bo{v}^{\perp} \in S^1$ denote the vector 
obtained by a counterclockwise turn of $\pi/2$ from the direction of $\bo{v}$. 
\end{definition}
\begin{definition}\label{defmarg}
For $P$ a probability measure on $\zoz$ and for $\omega' \in \{ 0,1 \}^A$ for some $A \subseteq E(\Z^2)$, we write $P_{\omega'}$ for the conditional law of $P$ given $\omega\big\vert_A = \omega'$. We will also write $P \big( \cdot \big\vert \omega' \big)$ for $P_{\omega'}$.
\end{definition}
\begin{definition}
Given a subset $A \subseteq \R^2$, two elements $\bo{x},\bo{y} \in \Z^2 \cap A$, 
  we write 
$\big\{ \bo{x} \build\leftrightarrow_{}^A \bo{y} \big\}$ for the subset of $\omega \in \{ 0,1 \}^{E(A)}$
 for which there exists an $\omega$-open path from $\bo{x}$ to $\bo{y}$ all of
whose edges lie in $E(A)$. By the ($\omega$-)open cluster of $\bo{x}$ in $A$, 
we mean the connected subset of
$A$ 
whose members lie in an edge 
belonging to some ($\omega$-)open path in $E(A)$ 
that begins at $\bo{x}$.
\end{definition}
Throughout, the notation $\dist \cdot \dist$ and $d \big( \cdot , \cdot \big)$ refers to the Euclidean metric on $\R^2$.  For $\gamma \subseteq \R^2$, we set ${\rm diam}(\gamma) = \sup \big\{ d\big(\bo{x},\bo{y}\big): \bo{x},\bo{y} \in \gamma \big\}$. For $K > 0$, we set $B_K = \big\{ \bo{x} \in \R^2: \dist \bo{x} \dist \leq K \big\}$.
\end{subsection}
\begin{subsection}{The Wulff shape and circuit centering}\label{seccircen}
We require an assertion on the profusion of regeneration sites in the conditioned circuit, for which some preliminaries are needed.
The macroscopic profile of the conditioned circuit is given by the boundary of the Wulff shape.
\begin{definition}
 We define the {\it inverse correlation length}: for $\bo{x} \in \R^2$,
$$
\xi(\bo{x}) = - \lim_{k \to \infty} k^{-1} \log \P_{\beta,q} \big( \bo{0} \leftrightarrow \lfloor k\bo{x} \rfloor  \big),
$$
where $\lfloor \bo{y} \rfloor \in \Z^2$ is the component-wise integer part of $\bo{y} \in \R^2$. 
\end{definition}
\begin{definition}
The unit-area Wulff shape \hfff{wulff} $\wulff$ is the compact set given by
$$
\wulff 
= \lambda \bigcap_{\bo{u} \in S^1} \Big\{ \bo{t} \in \R^2: \big( \bo{t},\bo{u} \big) \leq \xi\big(\bo{u}\big)  \Big\},
$$
where $(\cdot,\cdot)$ denotes the scalar product on $\R^2$, and where
the dilation factor $\lambda > 0$ is chosen to ensure that $\big\vert \wulff \big\vert = 1$.
\end{definition}
The following appears in Theorem B of \cite{civ}:
\begin{lemma}\label{lemozstr}
Let $P = \P_{\beta,q}$ with $\beta < \beta_c$  and $q \geq 1$.
Then $\wulff$ has a locally analytic, strictly convex boundary.
\end{lemma}
Global deviations of the conditioned circuit from the Wulff shape may be measured in the following way.
\begin{definition}
Let $\Gamma \subseteq \R^2$ denote a circuit.  Define its
global distortion ${\rm GD} \big( \Gamma \big)$ (from a factor $n$ dilate of the Wulff shape boundary) by means of \hfff{globdis}
\begin{equation}\label{eqdefgd}
{\rm GD} \big( \Gamma \big) = \inf_{\bo{z} \in \Z^2} d_H \Big(  n \partial \wulff + \bo{z}, \Gamma  \Big),
\end{equation}
where $d_H$ denotes the Hausdorff distance on sets in $\R^2$. 
\end{definition}
(In a general context, this would be a peculiar definition. However, we will work with this quantity only in the case of circuits that are conditioned to trap an area of at least, or exactly, $n^2$.)

To formulate a statement regarding the presence of radial regeneration sites, it is natural to work with circuits that are centred at the origin in a sense that we now specify.
\begin{definition}\label{defncentr}
Let $\Gamma \subseteq \R^2$ denote a circuit. 
The lattice point $\bo{z}$ attaining the minimum in (\ref{eqdefgd}) will be called the centre \hfff{centre} $\centre(\Gamma)$ of $\Gamma$. In the case that the minimum is not uniquely attained, we take $\bo{z}$ to be the lexicographically minimal among those points in $\Z^2$ that attain the minimum. 
\end{definition}
In fact, we may adopt any deterministic rule for breaking ties in the definition of  $\centre(\Gamma)$. The stated rule has been given simply for definiteness.
\begin{definition}
We write \hfff{acon} $\acon$ for the event $\big\{ \areacon \big\} \cap \big\{ \centre(\cir) = \bo{0} \big\}$.
\end{definition}
\end{subsection}
\begin{subsection}{Radial regeneration structure}\label{secrrg}
We are now ready to state the assertion on circuit regularity which we need.
\begin{definition}\label{defrg}
Let  $q \in (0,\pi/2)$ and $c \in (0,\pi)$.  

\hfff{forback} 
The $q$-forward cone $C^F_{\pi/2 - q} \big( \bo{v} \big)$ denotes the set of vectors $\bo{w} \in \R^2$ for which $\ang\big( \bo{w} - \bo{v}, \bo{v}^{\perp} \big) \leq \pi/2 - q$. The $q$-backward cone $C^B_{\pi/2 - q} \big( \bo{v} \big)$
 denotes the set of vectors $\bo{w} \in \R^2$ for which $\ang\big( \bo{w} - \bo{v}, - \bo{v}^{\perp} \big) \leq \pi/2 - q$.

Let $\Gamma$ denote a circuit for which $\bo{0} \in {\rm INT}\big(\Gamma \big)$.
A site $\bo{v} \in \Gamma$ is called a $(q,c)-\Gamma$-regeneration site (of $\Gamma$) if
\begin{equation}\label{eqdefreg}
\Gamma \cap 
 W_{\bo{v},c} \subseteq C^F_{\pi/2 - q} \big(  \bo{v} \big) \cup C^B_{\pi/2 - q} \big( \bo{v} \big).
\end{equation}
\end{definition} 
See Figure 1. Clearly, any $(q,c)-\Gamma$-regeneration site $\bo{v}$ has the property that the semi-infinite line from the origin through $\bo{v}$ cuts $\Gamma$ only at $\bo{v}$. As we will see, this absence of backtracking will permit us to perform surgery on circuit segments bounded by such regeneration sites. The form of the definition of regeneration sites also entails that the two parts of the circuit near such a site are well-separated, in a sense that will be more apparent  when we come to discuss the ratio-weak-mixing property in Section \ref{sechyp}.
\begin{figure}\label{figregdef}
\begin{center}
\includegraphics[width=0.3\textwidth]{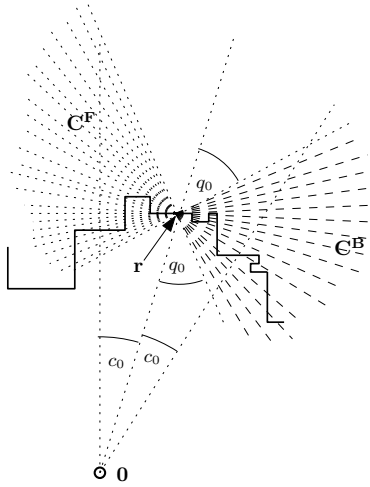} \\
\end{center}
\caption{A $(\qzero,c_0)-\Gamma$-regeneration site $\bo{r}$ and the nearby circuit. We write $C^F = C^F_{\pi/2 - \qzero}(\bo{r})$, and $C^B = C^B_{\pi/2 - \qzero}(\bo{r})$, and highlight these cones in dotted and dashed lines. Inside the cone $W_{\bo{r},c_0}$, the circuit remains within them.}
\end{figure}
\begin{definition}\label{defmar}
From now on, we fix $q_0 > 0$ and $c_0 \in  (0,\qzero/2)$ to be two small constants. (The precise requirements that they must satisfy are given in Definition 2.11 of \cite{hammondthr}. In addition to these requirements, we impose that $\qzero < 3\ccone/(4\cctwo)$, where $\ccone$ and $\cctwo$ are specified in the upcoming Lemma \ref{lemmacat}.)
Let $\Gamma$ denote a circuit for which $\bo{0} \in {\rm INT}\big(\Gamma \big)$.
We write \hfff{reg} ${\rm RG}(\Gamma)$ for the set of  $(\qzero,c_0)-\Gamma$-regeneration sites.
We write $\margam \in [0,2\pi]$ for the angle of the largest angular sector rooted at the origin that contains no element of ${\rm RG}(\Gamma)$. That is,
\begin{equation}\label{eqdefmar}
\margam = \sup \Big\{ r \in [0,2\pi): \exists \bo{a} \in S^1, 
 W_{\bo{a},r/2} \big( \bo{0} \big) \cap \reggam = \emptyset  \Big\}.
\end{equation}
\end{definition}
Theorem 2.1 of \cite{hammondthr} states that:
\begin{theorem}\label{thmmaxrg}
Let $P = \P_{\beta,q}$ with $\beta < \beta_c$  and $q \geq 1$. 
There exist $c > 0$ and $C > 0$ such that
$$
P \Big(  \mar > u/n   \Big\vert \acon \Big)
 \leq  \exp \Big\{ - c u \Big\} 
$$
for $n \in \N$ and $C \log n   \leq u \leq c n$. 
\end{theorem}
An alternative and stronger definition of regeneration site arises by instead considering the open cluster to which $\cir$ belongs. We denote this set by $\circl$. We define the set \hfff{regcl} $\regcl$ of {\it cluster regeneration sites} 
according to (\ref{eqdefreg}) (with $c = c_0$ and $q = \qzero$ understood), in which $\Gamma$ is replaced by $\circl$, and \hfff{marclo} $\marclo$ by (\ref{eqdefmar}) with the same change. We have Theorem 1.2 of \cite{hammondthr}:
\begin{theorem}\label{thmmaxrgcl}
Let $P = \P_{\beta,q}$ with $\beta < \beta_c$  and $q \geq 1$. 
There exist  $c > 0$ and $C > 0$ such that,
for $C \log n   \leq u \leq c n$, and for $n \in \N$,
$$
P \Big(  \marclo > u/n   \Big\vert \areacon \Big)
 \leq  \exp \Big\{ - c u \Big\}.  
$$ 
\end{theorem}
\end{subsection}
\begin{subsection}{Large deviations of global distortion}\label{secldgd}
A large deviations' estimate (Proposition 1 of \cite{hammondone}) on the macroscopic profile of the conditioned circuit will be valuable.
\begin{prop}\label{propglobdis}
Let $P = \P_{\beta,q}$ with $\beta < \beta_c$  and $q \geq 1$. 
There exists $c > 0$ and $n_0:(0,c) \to (0,\infty)$
such that, for any  $\epsilon \in \big( 0, c \big)$, and for all $n \geq n_0(\epsilon)$,
\begin{equation}\label{eqgd}
P \Big( \globdis > \epsilon n \Big\vert \big\vert {\rm INT} \big( \cir \big) \big\vert \geq n^2  \Big)  \leq \exp \big\{ - c \epsilon n \big\}. 
\end{equation}
Under this measure, $\centre\big( \cir \big) \in {\rm INT} \big( \cir \big)$ except with  probability that is exponentially decaying in $n$. Moreover, (\ref{eqgd}) holds under the conditional measure
$P \big( \cdot \big\vert \acon \big)$. 
\end{prop}

We mention that an analogue of this result in dimensions $d \geq 3$ is much more subtle. Proofs of such analogues have been undertaken by \cite{cerf} and \cite{cerfpisztora}. 

The following is an immediate consequence of Proposition \ref{propglobdis}. 
\begin{lemma}\label{lemmac}
There exist $\epsilon > 0$, $\ccone > 0$ and $\cctwo > 0$ such that, for $n \in \N$,
$$
P \Big(   \cir \subseteq B_{\cctwo n} \setminus B_{\ccone n}    \Big\vert \areacon \Big)
   \geq 1 - \exp \big\{ - \epsilon n \big\}.
$$
The same statement holds under 
$P \big( \cdot \big\vert \acon \big)$. 
\end{lemma}
\end{subsection}
\begin{subsection}{The upper bounds on maximum local deviation}
We will need to control from above maximum local roughness and maximum facet length by means of the main conclusions Theorem 1.1 and 1.2 of \cite{hammondone}:
\begin{theorem}\label{thmmlrbd}
Let $P = \P_{\beta,q}$, with $\beta < \beta_c$  and $q \geq 1$. 
Then there exist $C> c >0$ and $t_0 \geq 1$ such that, for $t \geq t_0$, 
$t = O \big( n^{5/36} (\log n)^{-C}  \big)$,
$$
P \Big( \mlr \geq n^{1/3} \big( \log n \big)^{2/3} t \Big\vert \areacon \Big)
 \leq \exp \Big\{ - c t^{6/5} \log n  \Big\}. 
$$
\end{theorem}
\begin{theorem}\label{thmmflbd}
Let $P = \P_{\beta,q}$, with $\beta < \beta_c$  and $q \geq 1$. 
There exist $0 < c < C < \infty$
such that, for $t \geq C$, 
$t = o \big( n^{1/3} (\log n)^{-5/6}  \big)$,
$$
P \Big( \mfl \geq n^{2/3} \big( \log n \big)^{1/3} t \Big\vert \areacon \Big)
 \leq \exp \Big\{ - c t^{3/2} \log n  \Big\}. 
$$
\end{theorem}
\end{subsection}
\begin{subsection}{Hypotheses on the configuration measure \texorpdfstring{$P$}{P}}\label{sechyp}
Most arguments in this paper and in \cite{hammondone} and \cite{hammondthr} depend on weaker hypotheses than being a subcritical random cluster measure $P = \P_{\beta,q}$, with $\beta < \beta_c$  and $q \geq 1$.
See Section 2.6 of \cite{hammondone} for a discussion. In this paper, we will sometimes refer to three properties that are satisfied by any $P = \P_{\beta,q}$, with $\beta < \beta_c$  and $q \geq 1$. We now introduce these properties.

\noindent{\bf Exponential decay of connectivity.}
The measure $P$ satsifes exponential decay of connectivity if there exists $c > 0$ such that $P_{\omega} \big( \bo{0} \leftrightarrow \partial B_n \big) \leq \exp \big\{ -cn \big\}$ for all $n \in \N$ and $\omega \in \{0,1\}^{E(\Z^2) \setminus E(B_n)}$. In this case, the property is satisfied by $P = \P_{\beta,q}$, with $\beta < \beta_c$  and $q \geq 1$ due to Theorem 1.2 in \cite{alexmix}.\\ 
\noindent{\bf Bounded energy.}
The following property  is trivially shown to be  satisfied by  $P = \P_{\beta,q}$, with $\beta < \beta_c$  and $q \geq 1$.
\begin{definition}\label{defbden}
A probability measure $P$ on $\{ 0,1  \}^{E(\Z^2)}$ satisfies 
the bounded energy property if there exists a constant $c > 0$ such that,  for any $\omega' \in \zoz$ and an edge $e \in E(\Z^2)$, the conditional probability that $\omega(e) = 1$ given the marginal $\omega' \big\vert_{E(\Z^2) \setminus \{ e \}}$ is bounded between $c$ and $1 - c$.
\end{definition}
\noindent{\bf Ratio-weak-mixing.}
The following spatial decorrelation property is well-suited to analysing the conditioned circuit.
\begin{definition}\label{defrwm}
A probability measure $P$ on $\{ 0,1  \}^{E(\Z^2)}$ is said to satisfy the
ratio-weak-mixing property if, for some $\crwm,\lambda > 0$, and for all sets 
$\mc{D}, \mc{F} \subseteq E \big( \Z^2 \big)$,
\begin{equation}\label{eqws}
\sup \Big\{ \Big\vert \frac{P \big( D \cap F \big)}{P\big( D \big) P \big(
  F \big)} - 1 \Big\vert: D \in \sigma_{\mc{D}}, F \in \sigma_{\mc{F}}, 
  P \big( D \big) P \big( F \big) > 0   \Big\} 
\leq \crwm \sum_{x \in V(\mc{D}), y \in V(\mc{F})} e^{- \lambda \vert x - y \vert},
\end{equation}
whenever the right-hand side of this expression is less than one. Here, for $A \subseteq E(\Z^2)$, $\sigma_A$ denotes the set of configuration events measurable with respect to the variables $\big\{ \omega(e): e \in A \big\}$.
\end{definition}
The ratio-weak-mixing property is satisfied by any  $P = \P_{\beta,q}$, with $\beta < \beta_c$  and $q \geq 1$,
 by Theorem 3.4 of \cite{alexon}.
We will make use of the property on only one occasion in the present article. We mention, however, that it has an important role to play in the surgeries performed in \cite{hammondone} and \cite{hammondthr}. Indeed, by making a definition of regeneration site that was a little stronger than the mere requirement that the circuit visit the radial line through the site only once, we obtain a useful quasi-independence for the configuration on either side of a regeneration site. This is because the regions $C^F$ and $C^B$ shown in Figure 1 are well-separated, in the sense that the sum on the right-hand-side of (\ref{eqws}) is bounded when the choice $\mc{D}=E(C^F)$ and $\mc{F} = E(C^B)$ is made. This means that the cone between two regeneration sites is comparatively unsullied by the details of conditioning in the exterior of the cone, making the cone a useful region for circuit surgery.
\end{subsection}
\begin{subsection}{Convention regarding constants}
An upper case  $C$ will be used for large positive constants, and a lower case $c$ for small positive constants.
The values of these two constants may change from line to line.
Some constants are fixed in all arguments: in particular, $c_0$ and $\qzero$ in the Definition \ref{defrg} of 
$\cir$-regeneration site, as well as $\ccone$ and $\cctwo$ in Lemma \ref{lemmac} (and a closely related Lemma \ref{lemmacat} to be presented later). 
\end{subsection}


\end{section} 
\begin{section}{The main argument}\label{secmain}
In this section, we present the proof of Theorem \ref{thmmlrubd}. We will begin by discussing the ideas of the proof, the considerations that determine the form of resampling of the circuit that we will undertake, and the outcomes of resampling that will provide the lower bound on maximum local roughness. With these elements motivated and formally introduced, we will then give the actual proof of Theorem \ref{thmmlrubd}.
Theorem \ref{thmmflubd} follows from Theorem  \ref{thmmlrubd} by a straightforward application of surgical techniques developed in \cite{hammondone}.  Its proof is presented at the end of the section. 
\begin{subsection}{The approach in outline}
We begin by describing the idea of the proof. More accurately, we will describe what would seem to be a natural way to trying to obtain Theorem \ref{thmmlrubd}, before mentioning some difficulties that will motivate some variations on the outlined approach.  Let $n \in \N$. We set 
\begin{equation}\label{thtval}
\theta_n = \ctheta n^{-1/3} (\log n)^{1/3},
\end{equation}
where $\ctheta$ denotes a constant that will eventually be fixed at a sufficiently small value.
We partition $\R^2$ into a set of consecutive sectors, $\secajn : = A_{j\theta_n,(j+1)\theta_n}$, $1 \leq j \leq \mn : =  \lfloor 2\pi/\theta_n \rfloor - 1$, rooted at the origin, of angle $\theta_n$. 
(In fact, the partition may include a narrower sector bordering the $x$-axis in the fourth quadrant. However, we will not make use of this sector.)
Note that
\begin{equation}\label{mtilde}
  \mn \sim
\frac{2\pi}{\ctheta} 
n^{1/3} 
\big( \log n \big)^{-1/3} 
\qquad 
\textrm{as $n \to \infty$}.
\end{equation}
For each sector, we seek to define a random replacement procedure $\zoz \to \zoz$, which, in acting on a configuration having the measure $P \big( \cdot \big\vert \acon \big)$, keeps the configuration in the complement of the sector, and resamples the sector configuration subject to maintaining the condition $\acon$. By definition, for each sector, the associated procedure leaves the measure  $P \big( \cdot \big\vert \acon \big)$ invariant. In outline, we will begin with a sample of the conditioned measure, and apply the procedure associated to the sectors in  order, running through all of them. 
At the end, we will still have a copy of the conditioned measure, which we will analyse. We will argue that each sector resampling has a probability of producing a section of the circuit $\cir$ that realizes the desired lower bound on $\mlr$ that is decaying as $n^{-\epsilon}$, for a small $\epsilon > 0$. There being a polynomial number of sectors, such a circumstance will arise often in the overall procedure. We will further argue that, once such a favourable configuration has arisen, the subsequent resamplings, further around the circle, cannot undo it, so that, to establish Theorem \ref{thmmlrubd}, we require only one favourable resampling to occur.

To implement this plan, note the following difficulties. In insisting that the sector resampling leaves invariant the law  $P \big( \cdot \big\vert \acon \big)$, we have no {\it a priori} way to describe the law of the resampled configuration inside the sector. We want to have an explicit description. To satisfy this wish, a natural suggestion is to attempt the resampling only if the two boundary lines of the sector in question cut through regeneration points (that is, elements of $\reg$) of the circuit $\cir$ in the input. 
The updated configuration in the sector would then appear in essence to be given by conditioning on the sector containing an open path between this pair of points in such a way that the circuit of which this path forms a part traps an area of at least $n^2$. 

A second difficulty is that the requirement that the circuit be centred at the origin, 
in the sense of Definition \ref{defncentr}, makes it difficult to describe the law of the resampled configuration, because this requirement imposes a constraint
 on the new open path in the sector
 that is difficult to express in a simple explicit form.
  To resolve this problem, we alter the definition of centering for a circuit in such a way that the centre is computed from data that we never resample: we will use only information in the left-hand half-plane to define the centre, and we will change our proposed resampling so as only to use sectors in the right-hand half-plane. This will permit us to describe the law of the resampled configuration in a sector without the distraction caused by the details of centering.
\end{subsection}
\begin{subsection}{Southwest centering and its relationship with the usual centering}
We now give 
the alterative definition of centering. 
\begin{definition}\label{defcsw}
For any circuit or path $\Gamma$, let \hfff{sw} $\sw{\Gamma} \in \Z^2$ denote the element of $V(\Gamma)$ 
that is minimal in the lexicographical ordering on $\Z^2$, in other words, 
whose $y$-coordinate is minimal among those vertices in $V(\Gamma)$ having minimal $x$-coordinate. Let $\csw$ denote the element $\bo{v} \in \Z^2$ for which $P \big( \sw{\cir} = \bo{v}  \big\vert \acon \big)$ is maximal. 
(If there are several elements at which the maximum is attained, we choose for $\bo{v}$ the lexicographically minimal among these.)
Note that, by the second statement of Lemma \ref{lemmac}, 
this probability is at least $cn^{-2}$, where $c = \frac{1}{2 \pi \cctwo^2}$. 

We say that a circuit $\Gamma$ is southwest-centred at $\bo{x} \in \Z^2$ if
$\sw{\Gamma} = \bo{x}$. Write ${\rm AREA}_{n^2,\bo{x}}^{\rm sw}$ for the event that there exists 
an open circuit 
$\Gamma$ satisfying $\big\vert {\rm INT} (\Gamma) \big\vert \geq n^2$ and $\sw{\Gamma} = \bo{x}$.
We further set \hfff{areatrap} $\areatrap = {\rm AREA}_{n^2,\csw}^{\rm sw}$.
Let
$\omega \in \areatrap$ be a configuration in which no point in $\R^2$ is enclosed by infinitely many open circuits (a condition which, as remarked after Definition \ref{defcirsur}, is satisfied $P$-a.s. whenever 
$P = \P_{\beta,q}$, with $\beta < \beta_c$  and $q \geq 1$).
We set \hfff{gamsw} $\gamsw = \gamsw(\omega)$ equal to the 
outermost open circuit southwest-centred at $\csw$. (Note that $\gamsw$ is well-defined, because it is given by the circuit $\Gamma$ for which ${\rm INT} (\Gamma)$ is the union of ${\rm INT} (\phi)$ over all open circuits $\phi$ 
that are southwest-centred at $\csw$.) We write \hfff{gamswcl} $\gamswcl$ for the open cluster in which $\gamsw$ is contained. 
Note that $\gamsw$ and $\gamswcl$ implicitly depend on the configuration $\omega \in \areatrap$. We make this dependence explicit only when doing so eliminates an ambiguity.
\end{definition}
There are some matters to take care of in switching from conditioning on $\acon$ to conditioning on  $\areatrap$. We discuss these before continuing the proof. The proofs of the lemmas stated here will be given in Section \ref{sectools}. 

We will work by conditioning $P$ on the event $\areatrap$, establishing the statement of Theorem \ref{thmmlrubd} for the conditional measure $P \big( \cdot \big\vert \areatrap \big)$. This done, we must infer the actual statement of Theorem \ref{thmmlrubd}. The device that we will use for this last step is now presented.
\begin{lemma}\label{lemlas}
 A set $\mathcal{M}$  of circuits is said to be translation-invariant if $\Gamma \in \mathcal{M}$
and $\bo{v} \in \Z^2$ imply that $\Gamma + \bo{v} \in \mathcal{M}$. There exist constants $C > c > 0$
such that, for any translation-invariant subset $\mathcal{M}$ of circuits, and for all $n \in \N$,
$$
P \Big( \cir \in \mathcal{M} \Big\vert \areacon \Big) \leq C n^4 P \Big( \gamsw \in \mathcal{M} \Big\vert \areatrap \Big)   + \exp \big\{ - c n  \big\}.
$$ 
\end{lemma}
We also need to exploit information about the circuit under the conditional measure known from \cite{hammondone}. There, results are stated for the circuit $\cir$ under the measure $P \big( \cdot \big\vert \areacon \big)$, whereas the present proof will consider the circuit $\gamsw$ under the measure $P \big( \cdot \big\vert \areatrap \big)$. The following lemma will be used to translate from the first framework into the second.
\begin{lemma}\label{lemoswtrans}
There exist constants $C > c > 0$
such that, for any translation-invariant subset $\mathcal{M}$ of circuits,
$$
P \Big( \gamsw \in \mathcal{M} \Big\vert \areatrap \Big) \leq C n^4 
P \Big( \cir \in \mathcal{M} \Big\vert \areacon \Big)   + \exp \big\{ - c n  \big\}.
$$ 
\end{lemma}
We will also use the following lemma, in which  the conditional probability of a collection of circuits that is not translation-invariant is considered.
\begin{lemma}\label{lemmacat}
There exist constants $\ccone$,$\cctwo$ and $c' > 0$ such that
$$
P \Big( \gamswcl \subseteq B_{\cctwo n} \setminus B_{\ccone n}, \bo{0} \in {\rm INT} \big( \gamsw \big) \Big\vert \areatrap  \Big) \geq 1 - \exp \big\{ - c' n \big\}.
$$
\end{lemma}
In Lemma \ref{lemmacat}, we retain the use of the notation $\ccone$ and $\cctwo$ for the two constants from the analogous Lemma \ref{lemmac}, although a change of value in the constants may be required to do so. 
\end{subsection}
\begin{subsection}{The use of cluster regeneration sites in resampling}
In light of the discussion so far, a natural proposal is to define the resampling procedure associated to a sector so as to leave invariant the law 
$P \big( \cdot \big\vert \areatrap \big)$, and not to attempt to resample (that is, to choose the resampling operation to act identically) unless the sector boundary lines intersect the input circuit $\gamsw$ at elements of $\rggamsw$.
There are unwanted complications in defining the resampling even now, however.
See Figure 2.
The sector boundary lines run through $\xo$ and $\yo$, each of which lies in $\rggamsw$. The input-configuration 
cluster $\gamswcl$ intersects 
these boundary lines also at $\bo{x_1}$ and $\bo{y_1}$. 
The sector resampling will be defined to map $\areatrap$ to itself,
so that, its output, for the input depicted, must contain a circuit that traps area at least $n^2$ and is southwest-centred at $\csw$. This eventuality will be effected  by the new randomness of the resampling
forging an open path across the sector that reconnects the fragment of the existing large circuit in the complement of the sector. 
In the case depicted, this circumstance may arise not only by an open path from $\xo$ to $\yo$ but also by such a path from, for example, $\bo{x_1}$ to $\bo{y_1}$. We prefer to work with a definition that eliminates these different cases, since this makes the law of the resampled configuration easier to describe. 
Recall the paragraph following Theorem \ref{thmmaxrg}.
The solution that we adopt is to attempt the sector resampling 
only under the stronger condition 
that the sector boundaries cut the circuit $\gamsw$ 
at cluster regeneration points. 
Since $\bo{x_0},\bo{y_0} \not\in \rggamswcl$ in the depicted case,
the sector resampling will not attempt a resampling for the input configuration in question, thereby eliminating the problem.  We will need a version of the cluster regeneration structure Theorem \ref{thmmaxrgcl} valid for the cluster $\gamswcl$ under the measure $P\big ( \cdot \big\vert \areatrap \big)$. 
\begin{theorem}\label{thmmaxrgclgw}
Let $P = \P_{\beta,q}$ with $\beta < \beta_c$  and $q \geq 1$. 
There exist $c > 0$ and $C > 0$ such that
\begin{equation}\label{eqmarcl}
P \Big(  \marcl > u/n   \Big\vert \areatrap \Big)
 \leq  \exp \Big\{ - c u \Big\} 
\end{equation}
for $C \log n   \leq u \leq c n$.
\end{theorem}
\noindent{\bf Remark.} In considering cluster regeneration sites of $\gamswcl$,
we will always use the parameter values $\qzero/2$ and $c_0/2$. As such, we set
$\rggamswcl = {\rm RG}_{\qzero/2,c_0/2}\big( \gamswcl \big)$.

The set of regeneration sites of a circuit depends on the choice of origin and is not translation-invariant. This means that Theorem \ref{thmmaxrgclgw} is not implied by Theorem \ref{thmmaxrgcl} by means of Lemma \ref{lemoswtrans}. We will present a different argument in Section \ref{sectools}.
\begin{figure}\label{figclreg}
\begin{center}
\includegraphics[width=0.4\textwidth]{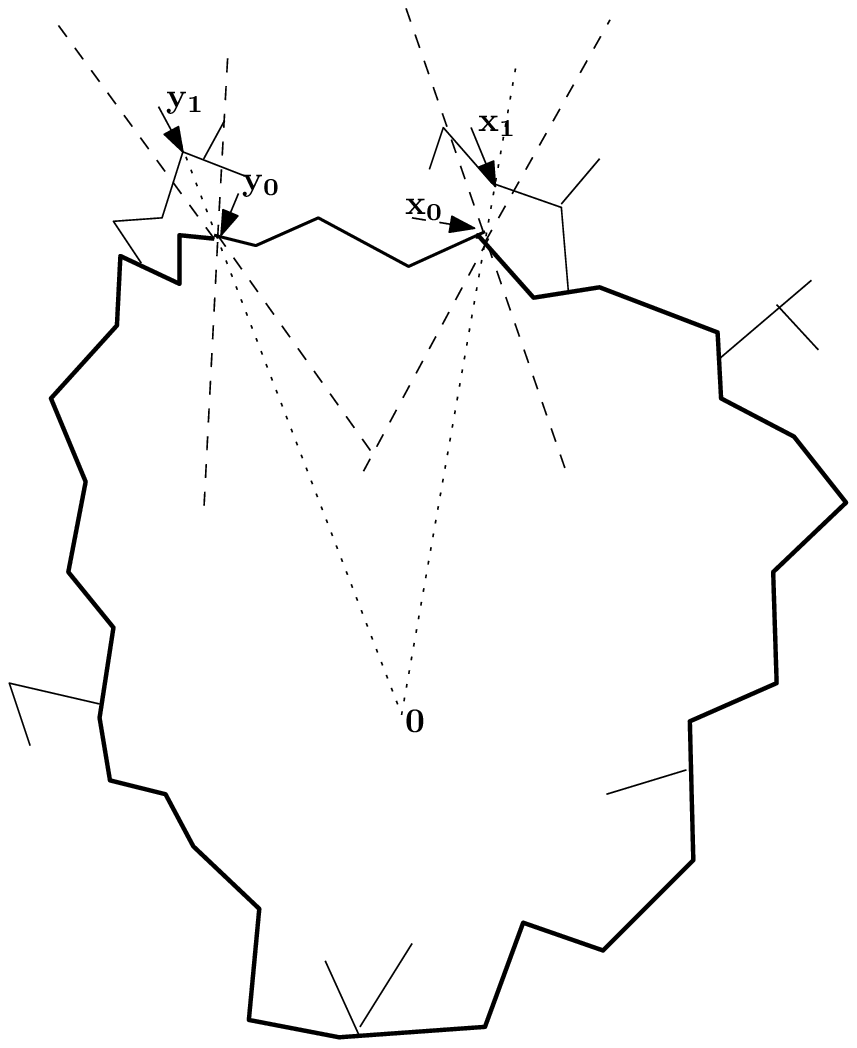} \\
\end{center}
\caption{A schematic depiction of the input for a sector-resampling. 
In solid is shown the cluster $\gamswcl$, with the circuit $\gamsw \subseteq \gamswcl$ highlighted in bold. 
The sector boundary lines (that are shown as dotted lines) intersect $\rggamswcl$ at $\xo,\bo{x_1},\yo$ and $\bo{y_1}$.
The dashed lines indicate the boundaries of 
$\cluh{\xo} \cup \clumh{\xo}$ 
and  
$\cluh{\yo} \cup \clumh{\yo}$. If we chose to perform a resampling procedure for a sector such as this, a new path might be forged between $\xo$ and $\yo$, or between $\xo$ and $\bo{y_1}$, for example. To avoid this ambiguity, we prefer not to resample such a configuration. We resample only if $\xo$ and $\yo$ are cluster regeneration sites, which, in the depicted case, they are not.}
\end{figure}
\end{subsection}
\begin{subsection}{Defining the resampling procedure in a sector}
To specify the sector-resampling procedure, we define:
\begin{definition}\label{defgagab}
Let $\bo{x},\bo{y} \in \Z^2$ satisfy $\argu(\bo{x}) < \argu(\bo{y})$.
Let $\omega \in \zoz$ realize the events that $\bo{x} \build\leftrightarrow_{}^{A_{\bo{x},\bo{y}}} \bo{y}$ and that the common cluster of $\bo{x}$ and $\bo{y}$ in $A_{\bo{x},\bo{y}}$ is finite. We define the outermost open path from $\bo{x}$ to $\bo{y}$ in $A_{\bo{x},\bo{y}}$ to be  the open path $\sentier$ 
from $\bo{x}$ to $\bo{y}$ in $A_{\bo{x},\bo{y}}$ such that 
the bounded component of 
$\axy \setminus \sentier$ is maximal. We denote this path by \hfff{gamoop} $\gamma_{\bo{x},\bo{y}} = \gamma_{\bo{x},\bo{y}}(\omega)$. (Note that there is only one bounded component, because a path is, by our definition, self-avoiding.)
We further write  \hfff{ogamoop} $\overline{\gamma}_{\bo{x},\bo{y}} =  \overline{\gamma}_{\bo{x},\bo{y}}(\omega)$ for the common $\omega$-open cluster of $\bo{x}$ and $\bo{y}$ in $\axy$. 
\end{definition}
\noindent{\bf Remark.} Note that $\gamma_{\bo{x},\bo{y}}$ is almost surely well-defined under $P$ given
  $\bo{x} \build\leftrightarrow_{}^{A_{\bo{x},\bo{y}}} \bo{y}$, with $P = \P_{\beta,q}$, $\beta < \beta_c$  and $q \geq 1$, since all open clusters are finite $P$-a.s. (as implied by the remark after Definition \ref{defcirsur}).
Note also that there is no difficulty in regard to the uniqueness of $\gamma_{\bo{x},\bo{y}}$, 
because we have decided to identify any path with the subset of $\R^2$ given by the union of the edges contained in the corresponding discrete path. 

We are now ready to specify the sector-resampling procedure. 
\begin{definition}\label{defrespsi}
Let $\bo{x},\bo{y} \in \Z^2$ lie in the right-hand half-plane and satisfy $\argu(\bo{x}) < \argu(\bo{y}) < \argu(\bo{x}) + c_0/2$. 
We define the circuit resampling operation $\psi_{\bo{x},\bo{y}}: \areatrap \to \zoz$. If the input $\omega \in \areatrap$ is such that $\bo{x} \not\in {\rm RG} \big( \gamswcl \big)$, or  $\bo{y} \not\in {\rm RG} \big( \gamswcl \big)$,
or $\bo{0} \not\in {\rm INT}\big( \gamsw \big)$,
then $\psi_{\bo{x},\bo{y}}\big( \omega \big) = \omega$. If  $\bo{x}, \bo{y} \in {\rm RG} \big( \gamswcl \big)$ and
 $\bo{0} \in {\rm INT}\big( \gamsw \big)$, then $\psi_{\bo{x},\bo{y}}\big( \omega \big)$ is taken to coincide with $\omega$ on $E(\Z^2) \setminus \axye$. On $\axye$, $\psi_{\bo{x},\bo{y}}\big( \omega \big)$ is chosen to be random, independently of 
$\omega \big\vert_{\axye}$, having the conditional distribution of the marginal of $P$ on $\axye$ given $\omega \big\vert_{E(\Z^2) \setminus \axye}$ and the following events: 
$\bo{x} \build\leftrightarrow_{}^{\axy} \bo{y}$,   
$\overline{\gamma}_{\bo{x},\bo{y}} \subseteq \cluh{\bo{x}} \cap \clumh{\bo{y}}$ and 
$\big\vert {\rm INT} \big( (\gamsw \cap \axy^c) \cup \gamma_{\bo{x},\bo{y}} \big) \big\vert \geq n^2$.
\end{definition}
\noindent{\bf Remark: the two-step formation of $\psi_{\bo{x},\bo{y}}$.} For later use, we introduce the following notation for describing the definition of $\psi_{\bo{x},\bo{y}}\big\vert_{\axye}$, in the case that it acts non-trivially. In step $A$, we condition the marginal of $P$ on $\axye$ by the data $\omega \big\vert_{E(\Z^2) \setminus \axye}$ and the event that  $\bo{x} \build\leftrightarrow_{}^{\axy} \bo{y}$. In step $B$, we further condition the marginal $P$ on $\axye$ by the events  
 $\big\vert {\rm INT} \big( (\gamsw \cap \axy^c) \cup \gamma_{\bo{x},\bo{y}} \big) \big\vert
 \geq n^2$ and $\overline\gamma_{\bo{x},\bo{y}} \subseteq \cluh{\bo{x}} \cap \clumh{\bo{y}}$.\\
\noindent{\bf Remark.} Note that, in the case that  $\psi_{\bo{x},\bo{y}}$ acts non-trivially,
 the intersection of the three events on which we condition necessarily has positive probability. This is because the input $\omega$ satisfies each of the conditions, and any so any configuration in $\axye$ that coincides with $\omega$ in a large finite neighbourhood of $\bo{0}$ will also do so.  
\begin{lemma}\label{leminv}
Let $\bo{x},\bo{y} \in \Z^2$ lie in the right-hand half-plane and satisfy $\argu(\bo{x}) < \argu(\bo{y}) < \argu(\bo{x}) + c_0/2$. 
The map  $\psi_{\bo{x},\bo{y}}$ leaves invariant the law 
$P \big( \cdot \big\vert \areatrap \big)$. If $\omega \in \areatrap$ is an input on which $\psi_{\bo{x},\bo{y}}$ acts non-identically, i.e., if
 $\bo{x}, \bo{y} \in {\rm RG} \big( \gamswcl(\omega) \big)$ and
 $\bo{0} \in {\rm INT}\big( \gamsw(\omega) \big)$, 
then, writing $\omega' = \psi_{\bo{x},\bo{y}}(\omega)$, we have that
$$
   \gamsw \big( \omega'  \big) =  \Big( \gamsw \big( \omega \big) \cap A_{\bo{x},\bo{y}}^c \Big) \cup \gamma_{\bo{x},\bo{y}}\big( \omega' \big). 
$$
\end{lemma}
The map $\psi_{\bo{x},\bo{y}}$ has been designed to ensure that this lemma holds. 
Verifying that the lemma holds requires some care in working with the definition of southwest-centering. The proof appears in Section \ref{sectoolsone}.

For each sector $\secajn$, we will select appropriate $\bo{x_j},\bo{y_j} \in \secajn$, and apply $\psi_{\bo{x_j},\bo{y_j}}$. We must have a chance of a favourable outcome that decays as a slow power. Specifically, we must have a means of choosing $\bo{x_j}$ and $\bo{y_j}$ that ensures that both belong to 
$\rggamswcl$ in the input with such a probability. For this purpose, deterministic choices of $\bo{x_j}$ and $\bo{y_j}$ appear not to suffice. Our definition of the sector-resampling procedure including the random search for its parameters is now given. The parameter $\epsilon$ that appears is fixed at a value in $(0,2/3)$
and will be specified in estimates in the proof of Theorem \ref{thmmlrubd}.
\begin{definition}\label{defres}
Let $j \in \{1,\ldots,\mn \}$ be such that $\secajn$ lies in the right-hand half-plane. We define a random resampling 
\hfff{resj} $\res_j:\areatrap \to \zoz$. 
Given $\omega \in \areatrap$, 
let $U_j^-$ and  $U_j^+$  
denote random variables whose respective laws are uniform on the intervals 
$\big[ j\thetan + \thetan/4 - n^{\epsilon - 1}, j \thetan + \thetan/4 \big]$ 
and  
$\big[ (j+1)\thetan - \thetan/4, (j+1)\thetan - \thetan/4 + n^{\epsilon - 1} \big]$.
These random variables are independent of each other and of $\omega$.
If $\bo{0} \not\in {\rm INT} \big( \gamsw \big)$, set
$\res_j(\omega) = \omega$. 
If $\bo{0} \in {\rm INT} \big( \gamsw \big)$, 
let $v_j^-$ denote the first edge in $\gamsw$ 
encountered on the line emanating from 
$\bo{0}$ 
at polar angle 
$U_j^-$, and let $\bo{x_j}$ denote the endpoint of $v_j^-$ of greater argument. Let $v_j^+$ be defined correspondingly with angle $U_j^+$ in place of $U_j^-$, and let $\bo{y_j}$ denote the endpoint of $v_j^+$ of smaller argument. 
We then set $\res_j(\omega) = \psi_{\bo{x_j},\bo{y_j}}\big( \omega \big)$. We will denote by $\big(\Omega,\mathcal{F},\mathbb{P}\big)$
the probability space in which $\res_j$ acts on an input having the distribution of $P \big(\cdot \big\vert \areatrap \big)$.
\end{definition} 
We require that:
\begin{lemma}\label{lemresjinv}
If $j \in \{1,\ldots,\mn \}$ is such that $\secajn$ lies in the right-hand half-plane,
then the law  $P \big( \cdot  \big\vert 
\areatrap \big)$ is invariant under  the map $\res_j$.
\end{lemma}
For the proof, we need the following.
\begin{lemma}\label{lemauvo}
Let $\bo{x},\bo{y}$, $\argu(\bo{x}) < \argu(\bo{y})$, lie in the right-hand half-plane.
 Under the conditional measure $P \big( \cdot \big\vert \areatrap \cap \big\{ \bo{x},\bo{y} \in {\rm RG} \big( \gamswcl \big) \big\} \cap \big\{ \bo{0} \in {\rm INT}(\gamsw) \big\}  \big)$,
the set $\gamsw \cap \axy^c$ is equal to the outermost open path from $\bo{x}$ to $\bo{y}$ in $\axy^c$
that is southwest-centred at $\csw$.
\end{lemma}
\noindent{\bf Proof.}
This statement is supplied during the proof of Lemma \ref{leminv} in Section \ref{sectoolsone}. See the paragraph following (\ref{gamtxy}). \qed
\noindent{\bf Proof of Lemma \ref{lemresjinv}.}
Let $\bo{u},\bo{v} \in \secajn$ be such that 
$\mathbb{P}\big( \bo{x_j} = \bo{u}, \bo{y_j} = \bo{v} \big) > 0$.
Set $\Theta_{\bo{u},\bo{v}} = \areatrap \cap \big\{ \bo{u},\bo{v} \in {\rm RG} \big( \gamswcl \big) \big\} \cap \big\{ \bo{0} \in {\rm INT}(\gamsw) \big\}  $.
Let $\omega_0 \in \{ 0,1 \}^{E(\Z^2) \setminus E(A_{\bo{u},\bo{v}})}$
take the form $\omega_0 = \omega \big\vert_{E(\Z^2) \setminus E(A_{\bo{u},\bo{v}})}$
for some $\omega \in \zoz$ realizing $\Theta_{\bo{u},\bo{v}}$. 

By the first assertion of Lemma \ref{leminv}, it suffices to show that, for any such $\omega_0$,
\begin{eqnarray}
 & & \textrm{the conditional distribution of $\omega\big\vert_{E(A_{\bo{u},\bo{v}})}$ is the same under} 
\label{doublepeq} \\
 & & \, \mathbb{P}\Big( \cdot \Big\vert \omega\big\vert_{E(\Z^2) \setminus E(A_{\bo{u},\bo{v}})} = \omega_0, \, \omega \in \Theta_{\bo{u},\bo{v}}, \, \bo{x_j} = \bo{u}, \, \bo{y_j} = \bo{v}  \Big) \nonumber \\
 & &   \textrm{and under} \, \,
 P\Big( \cdot \Big\vert \omega\big\vert_{E(\Z^2) \setminus E(A_{\bo{u},\bo{v}})} = \omega_0, \, \omega \in  \Theta_{\bo{u},\bo{v}}  \Big).
\nonumber
\end{eqnarray}
(Under the first conditional measure, as elsewhere when an operation is acting, $\omega \in \zoz$ denotes the input of the operation, in this case, $\res_j$.)
To see (\ref{doublepeq}), note that Lemma \ref{lemauvo} implies that, under 
$\mathbb{P}\big( \cdot \big\vert \omega\vert_{E(\Z^2) \setminus E(A_{\bo{u},\bo{v}})}  = \omega_0,
\omega \in \Theta_{\bo{u},\bo{v}}  \big)$,
the set $\gamsw(\omega) \cap A_{\bo{u},\bo{v}}^c$ is determined by the data $\omega_0$.
Let $\phi_1$ denote the nearest-neighbour edge in $\gamsw(\omega)  \cap A_{\bo{u},\bo{v}}^c$
that touches $\bo{u}$, (and is closer to the origin, if there is more than one). 
Let $\phi_2$ denote the edge in this set that touches $\bo{v}$, and is closer to the origin. 
It is easily seen that the collection $\phi_1'$ of points $\bo{w} \in \phi_1$ that are visible from the origin by an observer of $\gamsw(\omega)$, i.e, for which 
$\big[ \bo{0},\bo{w} \big] \cap \gamsw(\omega) = \big\{ \bo{w} \big\}$, is a subinterval of $\phi_1$ having $\bo{u}$ as an endpoint. Let  $\argu(\phi_1')$ 
denote the interval of argument-values of the elements of $\phi_1'$. 
The interval $\phi_2' \subseteq \phi_2$, defined correspondingly, has an endpoint at $\bo{v}$.
Then, under 
$\mathbb{P}\big( \cdot \big\vert \omega\vert_{E(\Z^2) \setminus E(A_{\bo{u},\bo{v}})}  = \omega_0, \omega \in \Theta_{\bo{u},\bo{v}}  \big)$,
the event  $\big\{ \bo{x_j} = \bo{u} \big\} \cap \big\{ \bo{y_j} = \bo{v} \big\}$ is given by
$\big\{ U_j^- \in \argu(\phi_1') \big\} \cap \big\{ U_j^+ \in \argu(\phi_2') \big\}$. 
However, this latter event is, under the same conditional measure, expressible in terms of the data $\omega_0$
and the independent randomness that generates $U_j^-$ and $U_J^+$. Hence, we obtain (\ref{doublepeq}), as required. \qed 
\begin{definition}
Let $j \in \{1,\ldots,\mn \}$ be such that $\secajn$ lies in the right-hand half-plane. 
If, in acting on an input $\omega \in \zoz$ for which $\bo{0} \in {\rm INT} \big( \gamsw \big)$, the pair $(\bo{x_j},\bo{y_j})$ located by $\res_j$ satisfies $\bo{x_j},\bo{y_j} \in {\rm RG}\big(\gamswcl\big)$, 
 we say that $\res_j$ selects $(\bo{x_j},\bo{y_j})$ successfully.
\end{definition}
\begin{lemma}\label{lemsucsel}
If $\omega \in \zoz$ satisfies $\big\{ \marcl \leq n^{\epsilon - 1}/2 \big\} \cap \big\{  \gamsw \subseteq B_{\cctwo n} \setminus B_{\ccone n} \big\}$, 
and $j \in \N$ is such that $\secajn$ lies in the right-hand half-plane, 
then the conditional probability that $\res_j$ in acting on $\omega$ selects  $(\bo{x_j},\bo{y_j})$ successfully is at least $(4 \cctwo)^{-2} \sin^2 \big( \qzero/2 \big) n^{-2\epsilon}$.
\end{lemma}
\noindent{\bf Proof.} 
Recall that $[ j \thetan + \thetan/4 - n^{\epsilon - 1}, j \thetan + \thetan/4 ]$
is the angular interval in which the direction $U_j^-$ used by $\res_j$ is chosen. 
The smaller interval  $\big[ j \thetan + \thetan/4 - n^{\epsilon - 1}/2, j \thetan + \thetan/4 \big]$
necessarily  
contains 
at least one argument value of an element of ${\rm RG}\big(\gamswcl\big)$. Select any such regeneration site and label it $\bo{r^-}$. Similarly, let $\bo{r^+}$ denote a regeneration site whose argument lies in the interval  $\big[ (j+1) \thetan - \thetan/4, (j+1) \thetan - \thetan/4 + n^{\epsilon - 1}/2 \big]$.

There are either one or two edges in $E(\gamsw)$ whose endpoint of greater argument equals $\bo{r^-}$.
Let $v^- \in E(\gamsw)$ denote the edge among these that is closer to the origin. 
We will argue that there exists a subinterval $\ell^-$ of $v^-$ having an endpoint at $\bo{r^-}$ such that an observer who stands at the origin and views $\gamsw$ sees all of $\ell^-$, with this subinterval occupying an angle of at least $\frac{\sin(\qzero/2)}{2\cctwo n}$ in the observer's field of vision. That is, writing $\ell^- = \big[ \bo{r_0^-},\bo{r^-} \big]$, then $\bo{v} \in \ell^-$, $\bo{v} \not= \bo{r_0^-}$ implies that
\begin{equation}\label{eqobsview}
\big[\bo{0},\bo{v} \big] \cap \gamsw = \{ \bo{v} \}, \quad \textrm{and} \quad 
\ang \big( \bo{r^-} , \bo{r_0^-} \big) \geq \frac{\sin(\qzero/2)}{2\cctwo n}. 
\end{equation}
We defer momentarily the derivation of (\ref{eqobsview}). From this statement, we find that, if $U_j^-$ is chosen so that its argument lies in $\big[ \argu(\bo{r_0^-}),\argu(\bo{r^-}) \big)$, then the choice of $\bo{x_j} = \bo{r^-}$ would be made by $\res_j$. Since $U_j^-$ is chosen on an angular interval of length $n^{\epsilon - 1}$ that, by the construction of $\bo{r^-}$, includes the angular interval occupied by $\ell^-$, this eventuality occurs with probability at least $\frac{1}{2\cctwo} \sin\big( \qzero/2 \big) n^{-\epsilon}$. 
A counterpart to (\ref{eqobsview}) shows that $\res_j$ independently makes a choice for $U_j^+$ that leads to $\bo{y_j} = \bo{r^+}$ with a probability that satisfies the same bound. This yields the statement of the lemma.

To derive (\ref{eqobsview}), suppose that the edge $v^-$ is vertical, so that 
$v^- = \big[ \bo{r^-} - \bo{e_2},\bo{r^-} \big]$, with $\big(\bo{e_1}, \bo{e_2}\big)$ denoting the  Cartesian unit vectors. See Figure 3.
By $\bo{r^-} \in \rggamswcl$ and $\gamsw \subseteq \gamswcl$, we know that
\begin{equation}\label{eqgamsw}
 \gamsw \cap W_{\bo{r^-},c_0/2}^- \subseteq C_{\qzero/2}^B\big(\bo{r^-} \big).
\end{equation}
Let $\ell^*$ denote that one among the pair of  semi-infinite boundary line segments of $C_{\qzero/2}^B\big(\bo{r^-} \big)$ that attains the closer approach to $\bo{0}$. Let $\bo{q}$ denote the point of intersection of $\ell^*$
and  the line segment 
$\big\{ \bo{r^-} - \bo{e_2} - t \bo{e_1}: t \geq 0 \big\}$. 
The observer at $\bo{0}$ who looks in a direction with argument in the interval 
$\big(  \argu(\bo{q}),\argu(\bo{r^-}) \big)$ sees no point of $\gamsw$ on the near-side of the line $\ell^*$ (by (\ref{eqgamsw})), while, beyond $\ell^*$, the first part of $\gamsw$ in the line of sight of the observer lies in $v^-$, by construction. That is, the observer viewing $\gamsw$ from $\bo{0}$ sees the edge $v^-$  in this angular interval.

\begin{figure}\label{figextra}
\begin{center}
\includegraphics[width=0.5\textwidth]{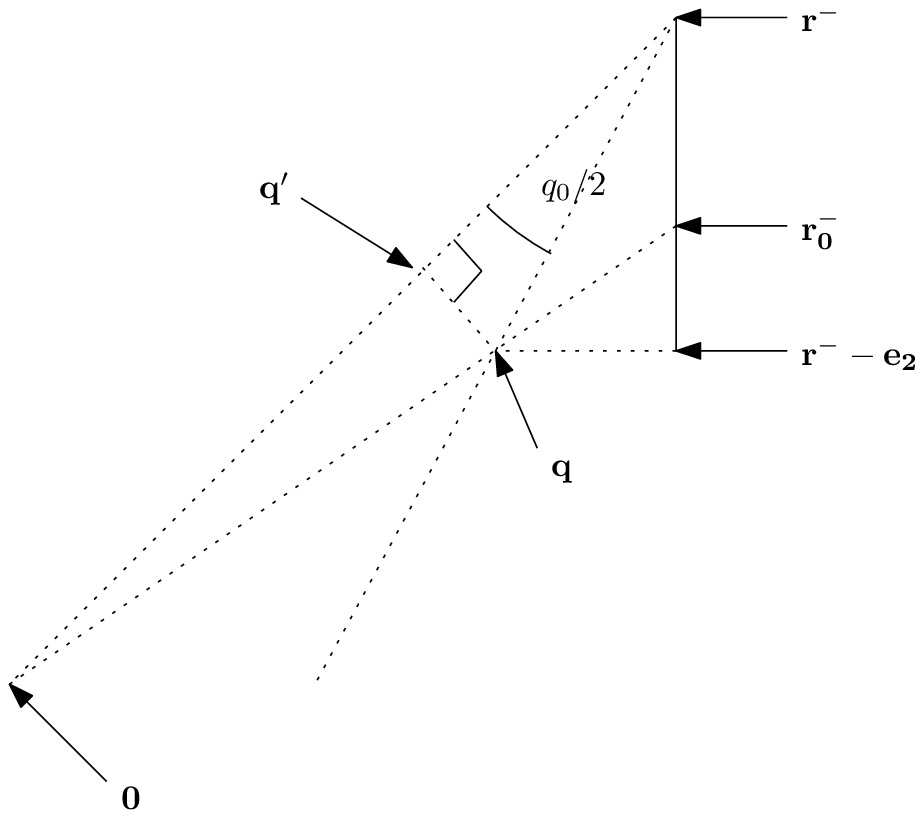} \\
\end{center}
\caption{Verifying (\ref{eqobsview}).}
\end{figure}
Setting  $\bo{r_0^-}$ to be the point in $v^-$ of argument $\argu(\bo{q})$,  
it remains to confirm the inequality in (\ref{eqobsview}).
Set $\bo{q'}$ equal to the point in $\ell_{\bo{0},\bo{r^-}}$ that is closest to $\bo{q}$.
Since $\bo{r-}$, $\bo{q'}$ and $\bo{q}$ form the vertices of a right-angled triangle, we have that
$\dist \bo{q'} - \bo{q} \dist = \dist \bo{r^-} - \bo{q} \dist \sin \big( \qzero/2 \big)$.
Note that $\dist \bo{r^-} - \bo{q} \dist \geq \dist \bo{r^-} - \big( \bo{r^-} - \bo{e_2} \big) \dist = 1$,
so that $\dist \bo{q'} - \bo{q} \dist \geq \sin \big( \qzero/2 \big)$. We have then that
$\ang \big( \bo{q},\bo{r^-} \big) \geq \dist \bo{q'} - \bo{q} \dist/\dist \bo{q} \dist
  = \dist \bo{q'} - \bo{q} \dist/\dist \bo{r^-} \dist \big( 1 + o(1) \big)$,
where the error term arises because $\bo{r^-} \in \gamsw \subseteq B_{\ccone n}^c$ is distant from the origin,
and because $\dist \bo{r^-} - \bo{q} \dist \leq \csc \big( \qzero/2 \big)$.
From $\bo{r^-} \in \gamsw \subseteq B_{\cctwo n}$, we obtain 
$\ang \big( \bo{q},\bo{r^-} \big) \geq \frac{\sin ( \qzero/2 )}{2\cctwo n}$, as required to complete the derivation of
(\ref{eqobsview}) in the case that $v^-$ is vertical. 
The case that this edge is horizontal is very similar, and is omitted. \qed
\end{subsection}
\begin{subsection}{Explaining the choice of sector angle \texorpdfstring{$\theta_n$}{theta n}}\label{secexpl}
\begin{figure}\label{figcircle}
\begin{center}
\includegraphics[width=0.4\textwidth]{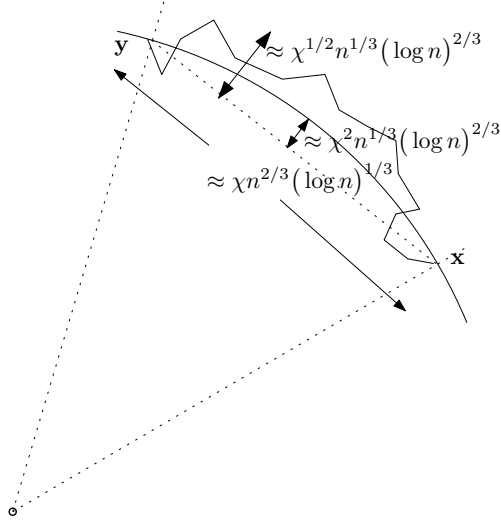} \\
\end{center}
\caption{The competition of fluctuation and curvature, heuristically illustrated: when $\bo{x}$ and $\bo{y}$ on a circle of radius $n$ have angular separation given by $\thetan = \thetan(\chi)$ as in (\ref{thtval}), the Gaussian orthogonal fluctuation of a subcritical path from $\bo{x}$ to $\bo{y}$ has a polynomially decaying probability of exceeding the orthogonal displacement of the circle between $\bo{x}$ and $\bo{y}$. If $\chi > 0$ is small, such fluctuations may be expected to appear often, among the order $n^{1/3}$ such sectors along the circle.}
\end{figure}
 To explain why the choice (\ref{thtval}) is natural, consider a circle $\mathcal{C}$ of radius $n$, such as the one depicted in Figure 4, as a simple example playing the role of the circuit $\gamsw$. 
Let $\bo{x},\bo{y} \in \mathcal{C}$ 
have angular displacement $\thetan$, where recall that $\thetan = \thetan(\ctheta)$. In this paragraph, 
for functions $f,g$ of $n$ and $\ctheta$,
we write $f = \Theta(g)$ to indicate that $0 < c \leq f/g \leq C < \infty$ for all $n \in \N$ and for $\ctheta \in (0,\delta)$, with $c$ and $C$ not depending on $\ctheta$, and with $\delta > 0$ some small constant. 

It is readily verified that 
the maximum distance $d_{\bo{x},\bo{y}}$ 
from a point in $A_{\bo{x},\bo{y}} \cap \mathcal{C}$ 
to $\big[ \bo{x},\bo{y} \big]$ is $\ctheta^2 \Theta \big( n^{1/3}( \log n )^{2/3} \big)$. 
Describing the formation of $\psi_{\bo{x},\bo{y}}$ as in the remark following Definition 
\ref{defrespsi}, the path $\gamma_{\bo{x},\bo{y}}$ in step $A$ has an $n^{-\epsilon}$ probability of fluctuating from $\big[ \bo{x},\bo{y} \big]$ away from the origin by 
$\phi(\epsilon) \vert\vert \bo{x} - \bo{y} \vert\vert^{1/2} \big( \log  \vert\vert \bo{x} - \bo{y} \vert\vert \big)^{1/2}$, where $\phi(\epsilon) > 0$ is a small constant (this is due to moderate deviations \cite{civ} of conditioned point-to-point connections in subcritical random cluster models having Gaussian behaviour). This last expression equals $\phi(\epsilon) \ctheta^{1/2} \Theta \big( n^{1/3} (\log n)^{2/3} \big)$. 
 For $\ctheta$ small, this fluctuation exceeds $d_{\bo{x},\bo{y}}$, 
so that the area condition in step $B$ 
of the formation of $\psi_{\bo{x},\bo{y}}$ is satisfied. 
So $\thetan$ has been tuned so that the resampling has a slow power-decaying probability that the point-to-point conditioned connection in step $A$ actually fluctuates enough  to satisfy the area requirement of step $B$. That is, the value of $\thetan$
resides on the scale at which the competition of local fluctuation and global curvature is finely balanced.
 Of course, we want a similar analysis without the heuristic assumption that the input circuit $\gamsw$ for $\res_j$ is circular.
\end{subsection}
\begin{subsection}{Finding a sufficient condition for the resampling to capture enough area}
We want to find an explicit criterion under which, given 
$\gamma_{\bo{x_j},\bo{y_j}}$ in step $A$, 
the additional conditioning  
\begin{equation}\label{addcon}
\Big\vert {\rm INT} \Big( \big( \gamsw \cap A_{\bo{x_j},\bo{y_j}}^c \big)  \cup \gamma_{\bo{x_j},\bo{y_j}} \Big) \Big\vert \geq n^2
\end{equation}
 appearing in step $B$ is satisfied. To this end, note that  
\begin{eqnarray}
 & & {\rm INT} \Big( ( \gamsw \cap A_{\bo{x_j},\bo{y_j}}^c ) \cup  \gamma_{\bo{x_j},\bo{y_j}}  \Big) 
\nonumber \\
 & = &
 {\rm INT} \Big(   \gamma_{\bo{x_j},\bo{y_j}} \cup \big[ \bo{0}, \bo{x_j} \big] 
                       \cup \big[ \bo{0}, \bo{y_j} \big] \Big)
 \cup
{\rm INT} \Big( ( \gamsw \cap A_{\bo{x_j},\bo{y_j}}^c ) \cup  
\big[ \bo{0}, \bo{x_j} \big] 
 \cup \big[ \bo{0}, \bo{y_j} \big] \Big), \label{gamswar}
\end{eqnarray}
since, if $\res_j$ is acting non-identically on $\omega$, then $\bo{0} \in {\rm INT}\big( \gamsw(\omega) \big)$
and
$\gamsw(\omega) \cap \ell_{\bo{0},\bo{z}}^+ = \{ \bo{z} \}$
 for  $\bo{z}=\bo{x},\bo{y}$. 

In the decomposition for ${\rm INT} \big( \gamsw \big)$ analogous to (\ref{gamswar}), the first term on the right-hand-side of (\ref{gamswar}) 
should be replaced by  
${\rm INT} \big(  \gamsw \big) \cap A_{\bo{x_j},\bo{y_j}}$. 
In defining $\res_j$, we consider only input $\omega$ 
satisfying 
$\big\vert  {\rm INT} \big( \gamsw(\omega) \big) \big\vert \geq n^2$.
The condition (\ref{addcon}) is thus implied by
$$ 
\Big\vert {\rm INT} \Big(   \gamma_{\bo{x_j},\bo{y_j}} \cup \big[ \bo{0}, \bo{x_j} \big] 
                       \cup \big[ \bo{0}, \bo{y_j} \big] \Big) \Big\vert 
\geq \Big\vert {\rm INT} \big(  \gamsw \big) \cap A_{\bo{x_j},\bo{y_j}} \Big\vert.
$$
We want to find an upper bound on the right-hand-side, to be used in a criterion sufficient for (\ref{addcon}) to be satisfied. We will use the following lemma. Although $\bo{x_j}$ and $\bo{y_j}$ are described as generic points in the statement, they will otherwise continue to denote the points selected in the definition \ref{defres} of the map $\res_j$.
\begin{lemma}\label{lemgamswjt}
Let $\omega \in \zoz$ be such that $\bo{0} \in {\rm INT} \big( \gamsw \big)$. 
Let $j \in \big\{ 0 ,\ldots, \mn \big\}$.
Write $\bo{z_j}$ for the element of  $\partial {\rm conv} \big( \gamsw \big)$ of argument $j \thetan$.
Let $\ell_j$ denote the tangent line of   $\partial {\rm conv} \big( \gamsw \big)$ at $\bo{z_j}$. 
(If $\bo{z_j}$ is an extremal point, we choose $\ell_j$ to be any line for which 
$\ell_j \cap \partial {\rm conv} \big( \gamsw \big)  = \{ \bo{z_j} \}$.)
Let $\bo{x_j},\bo{y_j} \in \gamsw \cap A_{\bo{z_j},\bo{z_{j+1}}}$.
Let $E_j$ be the pentagonal region bounded by the lines $\ell_{\bo{x_j},\bo{y_j}}$, $\ell_{\bo{0},\bo{x_j}}$, $\ell_j$, $\ell_{j+1}$ and $\ell_{\bo{0},\bo{y_j}}$.

Then
\begin{equation}\label{intatjt}
{\rm INT} \big(  \gamsw \big) \cap  A_{\bo{x_j},\bo{y_j}} \subseteq  T_{\bo{0},\bo{x_j},\bo{y_j}} \cup E_j,
\end{equation}
where recall that  $T_{\bo{0},\bo{x_j},\bo{y_j}}$ denotes the closed triangular region with vertices
$\bo{0}$, $\bo{x_j}$ and $\bo{y_j}$. 
\end{lemma}
\noindent{\bf Proof.} See Figure 5.  The region $T_{\bo{0},\bo{x_j},\bo{y_j}} \cup E_j$ is the bounded component $B$ of $A_{\bo{x_j},\bo{y_j}} \setminus \big( \ell_j \cup \ell_{j+1} \big)$. However, $\big( A_{\bo{x_j},\bo{y_j}} \setminus B \big) \cap {\rm conv}\big( \gamsw \big) = \emptyset$. \qed
\begin{figure}\label{figpent}
\begin{center}
\includegraphics[width=0.5\textwidth]{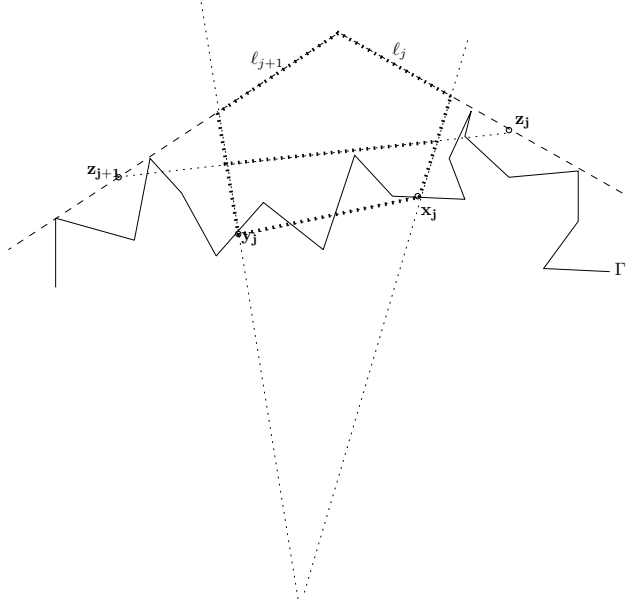} \\
\end{center}
\caption{Illustrating Lemma \ref{lemgamswjt}. The boundary of the pentagon $E_j$
is drawn as a bold dotted line, as is the line segment $E_j \cap \ell_{\bo{z_j},\bo{z_{j+1}}}$ which divides $E_j$ into $E_j^0$ (below) and $E_j^1$ (above).}
\end{figure}
The following criterion is therefore also sufficient for (\ref{addcon}): 
\begin{equation}\label{addintcon}
\big\vert {\rm INT} \big(  \gamma_{\bo{x_j},\bo{y_j}} \cup \big[ \bo{0}, \bo{x_j} \big] 
                       \cup \big[ \bo{0}, \bo{y_j} \big] \big) \big\vert 
 \geq  \big\vert T_{\bo{0},\bo{x_j},\bo{y_j}} \big\vert
  +   \big\vert E_j \big\vert. 
\end{equation}
It is natural then to seek an upper bound on  $\big\vert E_j \big\vert$. 
We write $E_j = E_j^0 \cup E_j^1$, where $E_j^0$ denotes those elements in $E_j$ lying on the same side of the line $\ell_{\bo{z_j},\bo{z_{j+1}}}$ as does $\bo{0}$, and where $E_j^1 = E_j \setminus E_j^0$.

For $\vert E_j^1 \vert$, we have:
\begin{lemma}\label{lemgbd}
Let $\omega \in \zoz$ be such that $\bo{0} \in {\rm INT} \big( \gamsw \big)$ and $\gamsw \subseteq B_{\cctwo n}$. 
Let $w_j$ denote the tangent vector of $\partial {\rm conv} \big( \gamsw \big)$ at $\bo{z_j}$ that points in the counterclockwise sense. 
The index set of sectors with {\rm moderate boundary turning} is defined by \hfff{mbt}
$$
\smtu = \Big\{ j \in \big\{1,\ldots,\mn \}: \vert\vert \bo{z_{j+1}} - \bo{z_j} \vert\vert \leq 40\pi \cctwo n/\mn, \ang \big( w_{j+1},w_j \big) \leq 40\pi/\mn \Big\}.
$$
%
Then $\big\vert \smtu \big\vert \geq 9\mn/10$. We have that
$j \in \smtu$ implies that 
$\big\vert E_j^1 \big\vert \leq 5^3 8^2 \cctwo^2 \ctheta^3 n \log n$.
\end{lemma}
\noindent{\bf Proof.} 
Let $T_j$ denote the triangle bounded by the lines $\ell_{\bo{z_j},\bo{z_{j+1}}}$, $\ell_j$ and $\ell_{j+1}$. Note that $E_j^1 \subseteq T_j$. We prove the statement of the lemma with $T_j$ in place of $E_j^1$.
Note that the sum of the angles of $T_j$ at $\bo{z_j}$ and at $\bo{z_{j+1}}$ is equal to $\ang \big( w_{j+1}, w_j \big)$. Provided that $\ang \big( w_{j+1}, w_j \big) \leq \pi/2$ 
(which holds if $j \in \smtu$), $T_j$ is contained in a right-angled triangle with hypotenuse $\big[ \bo{z_j}, \bo{z_{j+1}} \big]$ and another side contained in $\ell_{j+1}$. 
Hence, if  $j \in \smtu$,
\begin{equation}\label{tzsq}
\big\vert T_j \big\vert \leq \frac{1}{2} \vert\vert \bo{z_{j+1}} - \bo{z_j}  \vert\vert^2 \sin \ang \big( w_{j+1},w_j \big) \leq \frac{1}{2} \vert\vert \bo{z_{j+1}} - \bo{z_j}  \vert\vert^2 \ang \big( w_{j+1},w_j \big).
\end{equation}
Note that $\sum_{i=1}^{\mn} \vert\vert \bo{z_{i+1}} - \bo{z_i}  \vert\vert$ is the arclength of the convex set 
spanned by $\big\{ \bo{z_1},\ldots,\bo{z_m} \big\}$, which is contained in $B_{\cctwo n}$ (since this ball contains each $\bo{z_i}$). (In fact, we have omitted one side from consideration. This does not matter, since we are finding an upper bound on the sum.) It is readily verified that the arclength of the boundary of a convex set is at most that of any circle in whose interior the set is contained. Thus,  $\sum_{i=1}^{\mn} \vert\vert \bo{z_{i+1}} - \bo{z_i}  \vert\vert \leq 2 \pi \cctwo n$. Note also that  $\sum_{i=1}^{\mn-1} \ang \big( w_{i+1},w_i \big) \leq 2 \pi$. By definition, 
each element of $\big\{ 1,\ldots,\mn \big\} \setminus \smtu$ corresponds either to a long line segment, or a big angle (or both): the preceding two inequalities show that there are at most $\mn/20$ long line segments, and at most $\mn/20$ big angles.
Thus, $\big\vert \smtu \big\vert \geq 9\mn/10$. By (\ref{tzsq}) and (\ref{mtilde}), for all  $j \in \smtu$,
\begin{equation}\label{tisf}
\big\vert T_j \big\vert \leq \frac{1}{2} (40)^2 \pi^2 \cctwo^2 n^2 \frac{40 \pi}{\mn^3} \leq 
\frac{1}{2} \frac{(40)^3}{8} \cctwo^2 \ctheta^3 n \log n.
\end{equation}
Recalling that $E_j^1 \subseteq T_j$ completes the proof. \qed
We control $\vert E_j^0 \vert$ from above only if the circuit is close to its convex boundary in the sector $\secajn$. To quantify the hypothesis, we introduce some definitions:
\begin{definition}
For $\bo{u} \in \gamsw$, 
define the local roughness 
$\lr{u}$ of 
$\gamsw$ at 
$\bo{u}$ by means of 
$\lr{u} = d \big( \bo{u},\partial {\rm conv} \big( \gamsw \big) \big)$.
\end{definition}
\begin{definition}
Let $\functht > 0$ be chosen as a function of $\ctheta > 0$ in such a way that 
\hfff{functht} $\functht = o\big(\ctheta^{1/2}\big)$ as $\ctheta \to 0$.
Let $\omega \in \areatrap$. 
For $\ctheta > 0$ at a fixed value to be specified later, 
we say that the sector $\secajn$ is 
favourable under $\omega$ if there exists $\bo{u} \in \gamsw \cap \secajn$ such that $\lr{u} \geq \functht n^{1/3} \big( \log n \big)^{2/3}$.
We define \hfff{unfav} $\unfav = \unfav\big(\omega,\ctheta\big) \subseteq  \big\{1,\ldots,\mn \big\}$ to be the set of indices $j$ such that $\secajn$ is not favourable under $\omega$.
\end{definition}
\begin{lemma}\label{lemjiub}
Suppose that $\gamsw \cap B_{\ccone n} = \emptyset$.
For $j \in \smtu \cap \unfav$, $\big\vert E_j^0 \big\vert \leq 40 \cctwo \functht \ctheta n \log n$.
\end{lemma}
\noindent{\bf Proof.}
Let $H^j$ denote the component of $\R^2 \setminus \ell_{\bo{z_j},\bo{z_{j+1}}}$ containing $\bo{0}$. If $\bo{x_j},\bo{y_j} \in H^j$, then
the set $E_j^0$ is the quadrilateral bounded by the lines  
$\ell_{\bo{z_j},\bo{z_{j+1}}}$, $\ell_{\bo{0},\bo{x_j}}$,  $\ell_{\bo{0},\bo{y_j}}$ and  $\ell_{\bo{x_j},\bo{y_j}}$. If neither point lies in $H^j$, then $E_j^0 = \emptyset$ and there is nothing to prove. If exactly one point does, then $E_j^0$ is a triangle, and the ensuing proof applies almost verbatim to this case.  
Assume then that $\bo{x_j}, \bo{y_j} \in H^j$ (which is true in Figure 5). 
In a coordinate frame in which $\big[ \bo{z_j},\bo{z_{j+1}} \big]$ lies on the $x$-axis, with $\bo{0}$ (and $E_j^0$) in the upper-half-plane, $E_j^0$ has its longest intersection with any line $\big\{ y = c \big\}$ at $\big\{ y = 0 \big\}$. Thus, 
\begin{equation}\label{jjare}
\big\vert E_j^0 \big\vert \leq \vert\vert \bo{z_{j+1}} - \bo{z_j} \vert\vert
 \max \Big\{ d \big( \bo{y_j}, \big[ \bo{z_j},\bo{z_{j+1}} \big] \big) ,d \big( \bo{x_j}, \big[ \bo{z_j},\bo{z_{j+1}} \big] \big) \Big\}.
\end{equation}
We now show that
\begin{equation}\label{xjyj}
 d \big( \bo{z}, \big[ \bo{z_j},\bo{z_{j+1}} \big] \big) 
 \leq \lr{z} \qquad \textrm{for $\bo{z} \in \big\{ \bo{x_j},\bo{y_j} \big\}$}.
\end{equation}
Indeed, write $\lr{x_j} = d \big(\bo{x_j} ,\bo{v} \big)$ with $\bo{v} \in \partial {\rm conv} \big( \gamsw \big)$. If  $\bo{v} \in \secajn$, then there exists $\bo{v'} \in \big[ \bo{x_j},\bo{v} \big] \cap \big[ \bo{z_j},\bo{z_{j+1}} \big]$, since  $\partial {\rm conv} \big( \gamsw \big) \cap \secajn$ is separated from $\bo{0}$ by $\ell_{\bo{z_j},\bo{z_{j+1}}}$. Hence, $\lr{x_j} \geq d \big(\bo{x_j}, \bo{v'} \big) \geq d \big( \bo{x_j}, \big[ \bo{z_j},\bo{z_{j+1}} \big] \big)$. If, on the other hand, $\bo{v} \not\in \secajn$, then $\ang \big( \bo{v}, \bo{x_j} \big) \geq \thetan/8$
because $j \thetan + \thetan/8 \leq \argu(\bo{x_j}) < \argu(\bo{y_j}) \leq (j+1)\thetan - \thetan/8$. (It was in order to arrange this angular separation that we chose to search for $\bo{x_j}$ and $\bo{y_j}$ away from the boundaries of the sector $\secajn$.) We have then that 
$$
\vert\vert \bo{x_j} - \bo{v} \vert\vert \geq 
\vert\vert \bo{x_j} \vert\vert \sin \big( \thetan/8 \big) \geq 4^{-1} \pi^{-1} \ccone \ctheta n^{2/3} \big( \log n \big)^{1/3},   
$$
by the use of $\vert\vert \bo{x_j} \vert\vert \geq \ccone n$, which follows from $\bo{x_j} \in \gamsw$. Hence, $\lr{x_j}$ exceeds the right-hand-side of the last displayed equation. This contradicts $j \in \unfav$ (the lower bound on $\lr{x_j}$ is much too large), and completes the derivation of (\ref{xjyj}) for $\bo{z} = \bo{x_j}$. (The case that $\bo{z} = \bo{y_j}$ is identical.)

For $i \in \smtu \cap \unfav$,  
$\vert\vert \bo{z_{i+1}} - \bo{z_i} \vert\vert \leq 40 \pi \cctwo n \mn^{-1}$. 
We apply this and (\ref{xjyj}) in  (\ref{jjare}), and then bound
$\lr{x_j}$ and $\lr{y_j}$ by means of $j \in \unfav$, and use $(\ref{mtilde})$, to obtain the statement in the case that  $\bo{x_j},\bo{y_j} \in H^j$. \qed
Recall that (\ref{addintcon}) is sufficient for (\ref{addcon}).  Bringing together Lemma \ref{lemgamswjt}, $E_j = E_j^0 \cup E_j^1$, Lemma \ref{lemgbd} and Lemma \ref{lemjiub}, we find that
if 
$\omega \in \areatrap$ satisfies $\bo{0} \in {\rm INT} \big( \gamsw \big)$ and 
$\gamsw \subseteq B_{\cctwo n} \setminus B_{\ccone n}$, and
$j \in \big\{ 1,\ldots,\mn \big\}$ satisfies $j \in \smtu \cap \unfav$, then a sufficient condition on the path $\gamma_{\bo{x_j},\bo{y_j}}$ for (\ref{addcon}) to be satisfied is
\begin{equation}\label{accrit}
  \Big\vert {\rm INT} \Big(  \gamma_{\bo{x_j},\bo{y_j}} \cup \big[0, \bo{x_j} \big] \cup \big[0, \bo{y_j} \big] \Big) \Big\vert - 
\big\vert T_{\bo{0},\bo{x_j},\bo{y_j}}\big\vert
 \geq 
\Big( 5^3 8^2  \cctwo^2 \ctheta^3 +  40 \cctwo  \functht \ctheta \Big) n \log n.
\end{equation}
\end{subsection}
\begin{subsection}{Securing enough area and local deviation under sector resampling}
We will successively apply the maps $\res_j$, seeking to obtain an output that is favourable for the associated sector. The next two definitions will be used to define the event in terms of the resampling performed by $\res_j$ that we will show secures a favourable output.
\begin{definition}\label{defgac}
Let $\bo{x},\bo{y} \in \Z^2$ and let $\epsilon > 0$.
\hfff{gac}
Let the set ${\rm GAC}\big( \bo{x},\bo{y},\epsilon \big)$ of $\epsilon$-good area capture configurations in $A_{\bo{x},\bo{y}}$ denote the subset of 
$\omega \in \{ 0,1 \}^{\axye}$
such that the following conditions apply:
\begin{itemize}
 \item $\bo{x}  \build\leftrightarrow_{}^{\axy} \bo{y}$ under $\omega$,
 \item 
$\overline\gamma_{\bo{x},\bo{y}} \subseteq \cluh{\bo{x}} \cap \clumh{\bo{y}}$,
 \item  ${\rm diam} \big( \gamma_{\bo{x},\bo{y}} \big) \leq \cgac \vert\vert
\bo{x} - \bo{y} \vert\vert$,
 \item writing  \hfff{ixy} $I_{\bo{x},\bo{y}}\big( \gamma_{\bo{x},\bo{y}} \big) \subseteq \R^2$ for the bounded component of $\axy  \setminus \gamma_{\bo{x},\bo{y}}$, 
$$
 \Big\vert I_{\bo{x},\bo{y}} \big( \gamma_{\bo{x},\bo{y}} \big) \Big\vert
  \geq \big\vert  T_{\bo{0},\bo{x},\bo{y}} \big\vert + \epsilon
  \vert\vert \bo{x} - \bo{y} \vert\vert^{3/2} \big( \log \vert\vert \bo{x} -
  \bo{y} \vert\vert \big)^{1/2},
$$
where $T_{\bo{0},\bo{x},\bo{y}}$ is specified in Definition \ref{defto}.
\end{itemize}
\end{definition}
\begin{definition}
Let $\bo{x},\bo{y} \in \Z^2$, and let $\gamma \subseteq \axy$ denote a 
connected set for which $\{ \bo{x} \} \cup \{ \bo{y} \} \subseteq \gamma$. Let $\delta > 0$.
We say that $\gamma$ has \hfff{sid} $\delta$-{\rm significant inward deviation} if there exists $\bo{z} \in \gamma$ for which  $d\big( \bo{z}, \axy^c \big) \geq \vert\vert \bo{x} - \bo{y} \vert\vert \sin(\qzero)/3$ and
$$
d \Big(  \bo{z} , \partial {\rm conv} \big( \big[ \bo{0}, \bo{x} \big] \cup\big[ \bo{0}, \bo{y} \big] \cup \gamma    \big) \Big) \geq \delta \vert\vert \bo{x} - \bo{y} \vert\vert^{1/2} \Big( \log \vert\vert \bo{x} - \bo{y} \vert\vert \Big)^{1/2}.
$$
We write ${\rm SID} \big( \bo{x}, \bo{y},\delta \big)$ for the subset of configurations $\omega \in \{ 0,1 \}^{\axye}$ in $\axye$ such that $\bo{x} \build\leftrightarrow_{}^{\axy} \bo{y}$ under $\omega$, and for which the outermost open path $\gamma_{\bo{x},\bo{y}}$ from $\bo{x}$ to $\bo{y}$ in $\axy$  has $\delta$-significant inward deviation. 
\end{definition}
\begin{definition}\label{defactsuc}
Let $\omega \in \zoz$. We fix $\epso > 0$ at a small value to be specified later. We say that $\res_j$ acts successfully on $\omega$ if $\res_j$ selects 
the pair $(\bo{x_j},\bo{y_j})$ successfully, and if the resampling $\psi_{\bo{x_j},\bo{y_j}}$ applied by $\res_j$ realizes the event 
 ${\rm GAC} \big( \bo{x_j}, \bo{y_j},\epso \big) \cap \sid \big( \bo{x_j}, \bo{y_j},\epso \big)$. 
\end{definition}
Note that, if $\res_j$ acts successfully on an input $\omega$ for which $\gamsw \cap B_{\ccone n} = \emptyset$, then
\begin{equation}\label{angtht}
 \ang \big( \bo{x_j} , \bo{y_j} \big) \geq \thetan/4.
\end{equation}
Indeed, any successfully selected pair  $\big( \bo{x_j} , \bo{y_j} \big)$ satisfies this bound. 
To see this, recall 
from Definition \ref{defres} that this pair are selected by choosing two directions $U_j^-$ and $U_j^+$, with an angle between the two directions of at least $\thetan/2$. We then pivot about the origin counterclockwise from direction $U_j^-$ until the endpoint of the edge $v_j^-$ is reached, to find the point $\bo{x_j}$.
An analogous clockwise turning is made from direction $U_j^+$ in order to find $\bo{y_j}$.
The angle through which the respective turnings are made is bounded above by the angular width as viewed from $\bo{0}$ of $v_j^-$, or of $v_j^+$. Each of these edges is disjoint from $B_{\ccone n}$, so that each has an angular width of at most $2(\ccone n)^{-1}$. Hence, $\ang \big( \bo{x_j} , \bo{y_j} \big) \geq \thetan/2 - 4 (\ccone n)^{-1}$, so that (\ref{angtht}) follows from the definition of $\thetan$ in (\ref{thtval}).

The resamplings $\res_j$ have a slow power decaying probability of acting successfully, as long as $j \in \smtu \cap \unfav$:
\begin{lemma}\label{lemactsuc}
Let $j \in \N$ be such that $\secajn$ lies in the right-hand half-plane. Suppose that $\epsilon \in (0,2/3)$.
Let $\omega \in \zoz$ satisfy 
\begin{eqnarray}
& & \areatrap  \cap \big\{ \marcl \leq n^{\epsilon - 1}/2 \big\} \cap 
 \big\{ {\rm MLR} \big( \gamsw \big) \leq n^{2/3} \big\} \nonumber \\ 
& & \quad \cap 
\big\{  \gamsw \subseteq B_{\cctwo n} \setminus B_{\ccone n} \big\} \cap \big\{ \bo{0} \in {\rm INT} \big( \gamsw \big) \big\}, \nonumber
\end{eqnarray}
and let $j \in \smtu \cap \unfav$. 
Then there exists a constant $\clemgac > 0$ 
such that $\res_j$ acts successfully with probability at least $(4\cctwo)^{-2} \sin^2 \big( \qzero/2 \big)  n^{-2\epsilon - \clemgac \epso^2}$.
\end{lemma}
\noindent{\bf Proof.} Successful selection $(\bo{x_j},\bo{y_j})$ is 
dealt with by Lemma \ref{lemsucsel}. It remains to show that,  given $\bo{x_j}$ and $\bo{y_j}$,  $\res_j$ acts successfully with probability at least
 $n^{- \clemgac \epso^2}$.
We adopt the terminology of the remark following Definition \ref{defrespsi}, with the choice $\bo{x} = \bo{x_j}$ and $\bo{y} = \bo{y_j}$ being made. The path $\gamma_{\bo{x_j},\bo{y_j}}$ produced by the conditioning in step $A$ realizes the event ${\rm GAC} \big( \bo{x_j}, \bo{y_j},\epso \big) \cap \sid \big( \bo{x_j}, \bo{y_j},\epso \big)$ with probability at least $n^{- \clemgac \epso^2}$. 
This is a consequence of the bound presented in the following lemma, with the choice of $\bo{x} = \bo{x_j}$ and $\bo{y} = \bo{y_j}$ being made. The proof of Lemma \ref{lemgacmo} will be given after the present one is completed. Moreover, there is a little work required to show that the hypotheses (\ref{xyincone}) and (\ref{xyinctwo}) are satisfied for the choice  $\bo{x} = \bo{x_j}$ and $\bo{y} = \bo{y_j}$.
We will finish the rest of the proof of  Lemma \ref{lemactsuc} before 
justifying that these hypotheses are satisfied. 
\begin{lemma}\label{lemgacmo}
Let $P = \P_{\beta,q}$, with $\beta < \beta_c$  and $q \geq 1$.  There exists $\clemgac > 0$ and $n_0:(0,\infty) \to (0,\infty)$ such that the following holds. Let $\delta > 0$ and let $n \in \N$ 
satisfy $n \geq n_0(\delta)$. 
Let $\bo{x},\bo{y} \in \Z^2$ 
satisfy $\arg(\bo{x}) < \arg(\bo{y})$, $\vert\vert \bo{x} \vert\vert, \vert\vert \bo{y} \vert\vert \leq \clemgac n$, $\vert\vert \bo{y} - \bo{x} \vert\vert \geq  \clemgac \log n$, 
\begin{equation}\label{xyincone}
\bo{y} \in  C_{\pi/2 - \qzero}^F \big( \bo{x} \big) 
\end{equation}
and
\begin{equation}\label{xyinctwo}
\bo{x} \in  C_{\pi/2 - \qzero}^B \big( \bo{y} \big).
\end{equation}
Let $\omega \in \{0,1\}^{E(\Z^2)
  \setminus \axye}$ be arbitrary. 
Then
$$
P_{\omega} \Big( {\rm GAC} \big( \bo{x}, \bo{y},\delta \big) \cap \sid \big( \bo{x}, \bo{y},\delta \big)   \Big) \geq n^{-\clemgac \delta^2}  P \Big( \bo{x}  \build\leftrightarrow_{}^{\axy} \bo{y} \Big).
$$
\end{lemma}
We now argue that the 
event on which we condition in step $B$ of the formation of the law $\psi_{\bo{x_j},\bo{y_j}} \big( \omega \big)$
is automatically satisfied, if the $\gamma_{\bo{x_j},\bo{y_j}}$ that arises in step $A$ realizes the event ${\rm GAC} \big( \bo{x_j}, \bo{y_j},\epso \big) \cap \sid \big( \bo{x_j}, \bo{y_j},\epso \big)$. This suffices for the lemma, because it shows that the probability that $\gamma_{\bo{x_j},\bo{y_j}}$ realizes this event is at least as high after step $B$ as at the end of step $A$. We must verify the two requirements of step $B$. For the area bound, recall the sufficient criterion (\ref{accrit}). 

By (\ref{thtval}), (\ref{angtht}) and $\bo{x_j},\bo{y_j} \not\in B_{\ccone n}$, it follows that
\begin{equation}\label{distxjyj}
\vert\vert \bo{x_j} - \bo{y_j} \vert\vert \geq \frac{\ccone \ctheta}{\pi} n^{2/3} 
\big( \log n \big)^{1/3}. 
\end{equation}
Thus, by the occurrence of ${\rm GAC} \big( \bo{x_j}, \bo{y_j},\epso \big)$, the left-hand-side of (\ref{accrit}) exceeds $$
\frac{1}{2} (2/3)^{1/2}  
\big( \ccone/\pi \big)^{3/2} \epso \ctheta^{3/2}  n \log n.$$ 
By our choice 
$\functht = o\big(\ctheta^{1/2}\big)$ as $\ctheta \to 0$, we see that (\ref{accrit}) may be ensured by fixing $\ctheta > 0$ small enough. 

Note that the other condition in step $B$, that
$\overline\gamma_{\bo{x_j},\bo{y_j}} \subseteq \cluh{\bo{x_j}} \cap \clumh{\bo{y_j}}$, is satisfied, is implied by the occurrence of ${\rm GAC} \big( \bo{x_j}, \bo{y_j},\epso \big)$.

It remains to demonstrate that the pair $\big( \bo{x_j},\bo{y_j} \big)$ does satisfy the hypotheses of Lemma \ref{lemgacmo}. 
This is straightforward, except for the hypotheses (\ref{xyincone}) and (\ref{xyinctwo}). Regarding these, note firstly that we know that these two statements hold, with $\qzero/2$ in place of $\qzero$, from 
$\bo{x_j},\bo{y_j} \in \rggamswcl$ and $\ang \big( \bo{x_j},\bo{y_j} \big) \leq \theta_n  \leq c_0/2$.
The statements  (\ref{xyincone}) and (\ref{xyinctwo}) are stronger than this, however. It is in order to obtain them that we have introduced the hypothesis that ${\rm MLR} \big( \gamsw \big) \leq n^{2/3}$.

Indeed, suppose that $(\ref{xyincone})$ fails for the choice $\bo{x} = \bo{x_j}$ and $\bo{y} = \bo{y_j}$. (The case involving $(\ref{xyinctwo})$ is similar.)
We will show that this is inconsistent with our hypotheses on $\gamswcl$. Before discussing the peculiarities to be found in the circuit $\gamsw$ under our assumption, we make an observation regarding the circuit's convex boundary. Let $\bo{z} \in \partial {\rm conv} \big( \gamsw \big)$, and let $\bo{w_z}$ denote the tangent vector of $\partial {\rm conv} \big( \gamsw \big)$ at $\bo{z}$ (that points in the counterclockwise direction along the curve). Then we claim that 
\begin{equation}\label{wyminy}
\ang \big( \bo{w_z} , - \bo{z} \big) \geq \frac{\ccone}{\cctwo}. 
\end{equation}
Indeed, the tangent line $t_{\bo{z}}$ to $\partial {\rm conv} \big( \gamsw \big)$ at $\bo{z}$ 
cannot meet $B_{\ccone n}$, because all of $\gamsw$ lies on one side of $t_{\bo{z}}$, and so 
$t_{\bo{z}} \cap B_{\ccone n} \not= \emptyset$ would force the circuit $\gamsw$ into $B_{\ccone n}$, since
this circuit encircles the origin. We are assuming, however, that $\gamsw \cap B_{\ccone n} = \emptyset$.
The point $\bo{v}$ of closest approach of $t_{\bo{z}}$ to $\bo{0}$ satisfies
$$
 \dist \bo{v} \dist = \dist \bo{z} \dist \sin \ang \big( \bo{w_z}, - \bo{z} \big)
 \leq \dist \bo{z} \dist \ang \big( \bo{w_z}, - \bo{z} \big).
$$
Hence, $\dist \bo{v} \dist \geq \ccone n$ and
$\bo{z} \in {\rm conv} \big( \gamsw \big) \subseteq B_{\cctwo n}$ imply that (\ref{wyminy}) indeed holds.

Turning to consider $\gamsw$ under the assumption that (\ref{xyincone}) fails, note that
$\dist \bo{x_j} - \bo{y_j} \dist \geq 2 \pi^{-1} \ccone n \thetan/4$, by means of (\ref{angtht}) and
$\bo{x_j},\bo{y_j} \in \gamswcl \subseteq B_{\ccone n}^c$.  By (\ref{thtval}), then,
\begin{equation}\label{xjyjlog}
 \dist \bo{x_j} - \bo{y_j} \dist \geq c n^{2/3} \big( \log n \big)^{1/3}.
\end{equation}
By ${\rm MLR}\big( \gamsw \big) \leq n^{2/3}$ and $\bo{x_j},\bo{y_j} \in \gamsw$
(which is implied by $\bo{x_j}$ and $\bo{y_j}$ being cluster regeneration points of $\gamswcl$),
we may locate $\bo{x'},\bo{y'} \in \partial {\rm conv} \big( \gamsw \big)$ for which 
$\max \big\{ \dist \bo{x'} - \bo{x_j} \dist, \dist  \bo{y'} - \bo{y_j} \dist  \big\} \leq n^{2/3}$.
In view of (\ref{xjyjlog}), we have that $\ang \big(  \bo{y_j} - \bo{x_j}  ,  \bo{y'} - \bo{x'}  \big) \leq C \big( \log n \big)^{-1/3}$.

We are assuming that (\ref{xyincone}) fails for the choice $\bo{x} = \bo{x_j}$ and $\bo{y} = \bo{y_j}$. From this, the inequality that we have just derived, $\dist \bo{x'} - \bo{x_j} \dist = o(n)$ and $\dist \bo{x_j} \dist \geq \ccone n$, we learn that 
\begin{equation}\label{yprnotinc}
 \bo{y'} \not\in C_{\pi/2 - 4\qzero/3}^F \big( \bo{x'} \big),
\end{equation}
for $n$ sufficiently high. Without loss of generality, $\bo{y'}$ lies on the side of 
$C_{\pi/2 - 4\qzero/3}^F \big( \bo{x'} \big)$ for which 
\begin{equation}\label{zeroypr}
 \big[ \bo{0}, \bo{y'} \big] \cap C_{\pi/2 - 4\qzero/3}^F \big( \bo{x'} \big) = \emptyset.
\end{equation}
The counterclockwise-pointing tangent vector $\bo{w_{y'}}$ satisfies
$\ang \big( \bo{w_{y'}} , - \bo{y'} \big) \leq \ang \big( \bo{y'} - \bo{x'} , - \bo{y'} \big)$,
due to the counterclockwise turning of the circuit $\gamsw$ as it is traversed counterclockwise from $\bo{x'}$ to
$\bo{y'}$. However,
$$
\ang \big( \bo{y'} - \bo{x'} , - \bo{y'} \big) \leq 
\ang \big( \bo{y'} - \bo{x'} , - \bo{x'} \big) +
\ang \big( \bo{x'} , \bo{y'} \big) \leq 4 \qzero/3 + O \big( n^{-1/3} (\log n)^{1/3} \big).
$$
In the latter inequality, the bound on 
$\ang \big( \bo{y'} - \bo{x'} , - \bo{x'} \big)$ is due to (\ref{yprnotinc}) and (\ref{zeroypr})
and that on 
$\ang \big( \bo{x'} , \bo{y'} \big)$ is an easy consequence of (\ref{thtval}), the definition of $\bo{x'}$ and $\bo{y'}$, and $\bo{x_j},\bo{y_j} \not\in B_{\ccone n}$.
From the condition $\qzero < \frac{3 \ccone}{4 \cctwo}$ that we imposed on $\qzero$ in Definition \ref{defmar}, we see that this contradicts (\ref{wyminy}). This establishes (\ref{xyincone}) with the choice 
$\bo{x} = \bo{x_j}$ and $\bo{y} = \bo{y_j}$, and completes the proof of Lemma \ref{lemactsuc}. \qed
\noindent{\bf Proof of Lemma \ref{lemgacmo}.} 
A proof of the statement for the quantity 
$P_{\omega} \big( {\rm GAC} ( \bo{x}, \bo{y},\delta ) \big)$ is given in Lemma 3.1 of \cite{hammondone}.
The event $\sid \big( \bo{x}, \bo{y},\delta  \big)$ requires a deviation toward the origin on the part of $\gamma_{\bo{x},\bo{y}}$, and, as such, it appears to be negatively correlated with $\gacc \big( \bo{x}, \bo{y},\delta  \big)$. Nonetheless, the lemma follows by a minor variation of the proof of Lemma 3.1 of \cite{hammondone} that we now describe.

We will make use of a coordinate frame for $\R^2$ in which $\ell_{\bo{x},\bo{y}}$ is horizontal, with origin equal to $\bo{y}$ and with $\bo{x}$ having positive $x$-coordinate,  and with $\bo{0}$ in the usual coordinates lying in the lower half-plane. 
Set $h = \vert\vert \bo{y} - \bo{x} \vert\vert$. (We omit integer-rounding from our notation, and assume that $h$ and related quantities are integers.) 
Using the new coordinate system, we write $\bo{x_1} = \big( h/4, 10\delta \sqrt{h}  ( \log h )^{1/2} \big)$, 
$\bo{x_2} = \big( h/2, 5\delta \sqrt{h}  ( \log h )^{1/2} \big)$ and $\bo{x_3} = \big( 3h/4, 10\delta \sqrt{h}  ( \log h )^{1/2} \big)$. Further set $\bo{x_0} = \bo{y}$ and $\bo{x_4} = \bo{x}$. 
\begin{figure}\label{figguide}
\begin{center}
\includegraphics[width=0.45\textwidth]{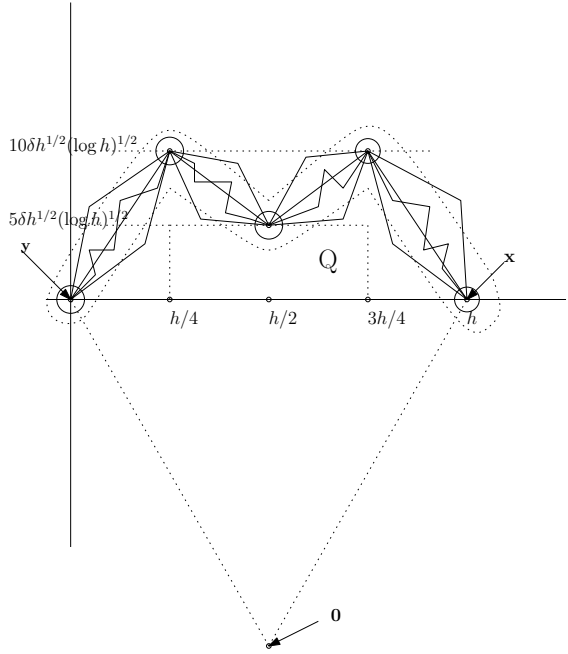} \\
\end{center}
\caption{Illustrating the proof of Lemma \ref{lemgacmo}. The five circles are the boundaries of the  radius-$K$ balls about the successive points $\bo{y} = \bo{x_0}$ on the left up to $\bo{x} = \bo{x_4}$ on the right. The four diamond-shaped regions each enclosing a line segment $[\bo{x_i},\bo{x_{i+1}}]$  are the $R_i$ for $i=0,1,2,3$ from left to right. 
The dotted curve that surrounds the union of the segments $[\bo{x_i},\bo{x_{i+1}}]$ is the boundary of the region $N$. An instance of $\overline{\gamma}_{\bo{x},\bo{y}}$ that realizes the event $\outfluc$ is depicted.}
\end{figure} 

For $i \in \big\{ 0, 1,2,3 \big\}$, 
let $R_i =  W_{\bo{x_{i + 1}} - \bo{x_i},\qzero/4} \big( \bo{x_i}  \big) \cap  W_{\bo{x_i} - \bo{x_{i+1}},\qzero/4} \big( \bo{x_{i+1}}\big)$. Note that we may find $K \in \N$ such that there exists an infinite simple (lattice) path from $\bo{0}$ in $B_K \cup W$, for every aperture-$\qzero/2$ cone $W$ with apex at $\bo{0}$. We fix $K \in \N$ at such a value, independently of the value of $h$.
Set $B_0 = B_K(\bo{x_0}) \cap A_{\bo{x},\bo{y}}$, $B_1 =  B_K(\bo{x_1})$,  $B_2 =  B_K(\bo{x_2})$, $B_3 =  B_K(\bo{x_3})$ and  $B_4 = B_K(\bo{x_4}) \cap A_{\bo{x},\bo{y}}$. For such $i$, let
$H_i$ denote the event that $\bo{x_i} \leftrightarrow \bo{x_{i+1}}$ in $B_i \cup R_i \cup B_{i+1}$, with the common connected component $\overline{\gamma}_i$ of $\bo{x_i}$ and $\bo{x_{i+1}}$ in  $B_i \cup R_i \cup B_{i+1}$  intersecting  
  $\partial \big( B_i \cup R_i \cup B_{i+1} \big)$ only in $\partial \big( B_i \cup B_{i+1}\big)$ and satisfying
$$
\sup \Big\{ d \big( \bo{v}, \big[  \bo{x_i}, \bo{x_{i+1}} \big] \big): \bo{v} \in \overline{\gamma}_i \Big\}
\leq 10 \vert\vert  \bo{x_{i+1}} - \bo{x_i} \vert\vert^{1/2}.
$$
Let $J$ denote the event that $\overline\gamma_{\bo{x},\bo{y}} \cap B_K(\bo{x}) \subseteq 
C^F_{\pi/2 - \qzero/2}(\bo{x})$ and
 $\overline\gamma_{\bo{x},\bo{y}} \cap B_K(\bo{y}) \subseteq 
C^B_{\pi/2 - \qzero/2}(\bo{x})$. 
 
Set $\outfluc = H_0 \cap H_1 \cap H_2 \cap H_3 \cap J$. 

Set $L$ equal to the union of the line segments  
$\big[\bo{x_i},\bo{x_{i+1}} \big]$ for $0 \leq i \leq 3$.
Let $N$ denote the $10 h^{1/2}$-neigbourhood of $L$. Note that $\outfluc$ implies that
 $\bo{x} \build\leftrightarrow_{}^{\axy} \bo{y}$, with $\bo{x_i} \in \overline{\gamma}_{\bo{x},\bo{y}}$ for $1 \leq  i \leq 3$. Moreover, $\overline{\gamma}_{\bo{x},\bo{y}} \subseteq N$. Hence,
$$
 I_{\bo{x},\bo{y}} \big( L \big) 
\subseteq I_{\bo{x},\bo{y}}\big(\gamma_{\bo{x},\bo{y}}\big) \cup N.
$$
Let $Q$ denote the rectangle that, in the chosen coordinates, has the form 
$Q = \big[ h/4,3h/4 \big] \times \big[ 0, 5 \delta h^{1/2} \big( \log h \big)^{1/2} \big]$.
We have that the sets $T_{\bo{0},\bo{x},\bo{y}}$ and $Q$ are disjoint, each of them being a subset of 
$I_{\bo{x},\bo{y}}\big( L \big)$. Noting that $\vert N \vert \leq 20 h^{3/2}$, and that
$\vert Q \vert = (5/2) \delta h^{3/2} \big( \log h  \big)^{1/2}$, it follows that, for $n$ high enough,
$$
  \Big\vert I_{\bo{x},\bo{y}}\big( \gamma_{\bo{x},\bo{y}} \big) \Big\vert
 \geq  \big\vert T_{\bo{0},\bo{x},\bo{y}} \big\vert + \delta  h^{3/2} \big( \log h  \big)^{1/2},
$$
since $h = \dist \bo{x} - \bo{y} \dist \to \infty$
as $n \to \infty$. To confirm
$\outfluc \subseteq \gacc \big( \bo{x},\bo{y},\delta \big)$, it remains to verify the second condition listed in Definition \ref{defgac}. 
To this end, set $\overline{B} = \cup_{i=0}^3 \big( R_i \cup B_i \big) \cup R_4$.
Note that $\outfluc$ entails that $\overline{\gamma}_{\bo{x},\bo{y}} \subseteq \overline{B}$, and hence that
 $\gamma_{\bo{x},\bo{y}} \subseteq \overline{B}$. From (\ref{xyincone}) and (\ref{xyinctwo}) and the definitions of the constituent sets of $\overline{B}$,  any point $\bo{u} \in \overline{B} \setminus \big( B_0 \cup B_4 \big)$ 
 is such that each of the angles $\ang \big( \bo{u} - \bo{x}, \bo{x}^{\perp} \big)$
 and $\ang \big( \bo{u} - \bo{x}, - \bo{y}^{\perp} \big)$ is at most $\pi/2 - \qzero/2$.
Hence, $\overline\gamma_{\bo{x},\bo{y}}$ satisfies the containment in this second condition except possibly in $B_0 \cup B_4$; but, in this region, the containment is assured by  the occurrence of $J$. Hence, indeed, we have  the second condition listed in Definition \ref{defgac}.

We will now verify that
\begin{equation}\label{outflinc}
 \outfluc \subseteq \sid \big( \bo{x},\bo{y},\delta \big).
\end{equation}
By $\outfluc \subseteq  \big\{ \gamma_{\bo{x},\bo{y}} \subseteq \overline{B} \big\}$, we may, if $\outfluc$ occurs, 
find $\bo{z_i} \in B_i \cap \gamma_{\bo{x},\bo{y}}$ for $i \in \{1,2,3\}$, since $\bo{x}$
and $\bo{y}$ are not connected in any of the sets $\overline{B} \setminus R_i$, for $i=1,2,3$.

We will check that $\bo{z_2}$ realizes the event  $\sid \big( \bo{x},\bo{y},\delta \big)$.

Note that 
$$
 {\rm conv} \Big( \big[ \bo{0},\bo{x} \big] \cup \big[ \bo{x},\bo{z_1} \big] \cup
   \big[ \bo{z_1},\bo{z_3} \big] \cup \big[ \bo{z_3},\bo{y} \big] \cup  \big[ \bo{0},\bo{y} \big]  \Big)
\subseteq 
{\rm conv} \Big( \big[ \bo{0},\bo{x} \big] \cup \gamma_{\bo{x},\bo{y}} \cup \big[ \bo{0},\bo{y} \big] \Big),
$$
since $\bo{z_1},\bo{z_3} \in \gamma_{\bo{x},\bo{y}}$.
Note further that $\bo{z_2}$ belongs to the left-hand-side of this inclusion, by construction, and that its distance to the boundary of this set is attained on the interval $\big[ \bo{z_1},\bo{z_3} \big]$. Hence,
$$
d \Big( \bo{z_2}, \partial 
{\rm conv} \big( \big[ \bo{0},\bo{x} \big] \cup \gamma_{\bo{x},\bo{y}} \cup \big[ \bo{0},\bo{y} \big] \big) 
\Big) \geq d \big(  \bo{z_2}, \big[ \bo{z_1},\bo{z_3} \big] \big).
$$
Now, 
$d \big(  \bo{z_2}, \big[ \bo{z_1},\bo{z_3} \big] \big) \geq 5 \delta h^{1/2} \big( \log h \big)^{1/2} - 2K$.
Hence, 
$$
d \Big( \bo{z_2}, \partial {\rm conv} \big( \big[ \bo{0},\bo{x} \big] \cup \gamma_{\bo{x},\bo{y}} \cup \big[ \bo{0},\bo{y} \big]  
\big) \Big) \geq 4 \delta  h^{1/2} \big( \log h \big)^{1/2}. 
$$
To confirm (\ref{outflinc}),
it remains to check that
\begin{equation}\label{eqdzerot}
  d \big( \bo{z_2}, \axy^c \big) \geq \dist \bo{x} - \bo{y} \dist \sin (\qzero)/3.
\end{equation}
Let $\bo{z'}$ denote the midpoint of $\big[ \bo{x},\bo{y} \big]$. Note that
\begin{equation}\label{eqzprz}
 d \big( \bo{z_2}, \bo{z'} \big) \leq 5 \delta  h^{1/2} \big( \log h \big)^{1/2} + K.
\end{equation}
Note that, since $\bo{x} \in C_{\pi/2 - \qzero}^B (\bo{y})$, 
the line segment $\big[ \bo{x},\bo{y} \big]$ makes an angle of at least $\qzero$
with $\ell_{\bo{0},\bo{y}}$. From $\bo{z'} \in \big[ \bo{x},\bo{y} \big]$, we find that
$d \big( \bo{z'}, \ell_{\bo{0},\bo{y}} \big) \geq d \big( \bo{z'}, \bo{y} \big) \sin (\qzero)
 =  d \big( \bo{x},\bo{y} \big) \sin (\qzero)/2$. By (\ref{eqzprz}),
we obtain, for $n$ sufficiently high, (since $d(\bo{x},\bo{y}) = h \to \infty$ as $n \to \infty$), 
$$
 d \big( \bo{z_2}, \ell_{\bo{0},\bo{y}} \big) \geq 
d \big( \bo{x},\bo{y} \big) \sin (\qzero)/2.
$$
By this, and the same inequality on  $d \big( \bo{z_2}, \ell_{\bo{0},\bo{x}} \big)$, we obtain 
(\ref{eqdzerot}), and hence (\ref{outflinc}).

The following assertion is analogous to (3.21) in \cite{hammondone}, and may be obtained by the argument leading to that result.  There exist constants $C > c > 0$ such that, for $n$ sufficiently high, 
\begin{equation}\label{eqoutflbd}
  P_\omega \Big( \outfluc \Big) \geq c h^{-C\delta^2} h^{-3/2} P \Big( \bo{x} \leftrightarrow \bo{y} \Big).
\end{equation}

Similarly to the paragraph that follows (3.26) in \cite{hammondone}, we could conclude now, if the factor of $n^{-3/2}$ were absent on the right-hand-side of (\ref{eqoutflbd}), since $h \leq 2 \clemgac n$. This factor arises due to three ``local limit theorem'' factors of order $n^{-1/2}$ arising from the insistence that, under $\outfluc$, $\overline\gamma_{\bo{x},\bo{y}}$ pass through the three bounded neighbourhoods $B_1,B_2$ and $B_3$ at the x-coordinates $h/4$,$h/2$ and $3h/4$. We may consider variants of $\outfluc$, in which the vertical coordinates of each of $\bo{x_1}$,$\bo{x_2}$ and $\bo{x_3}$, assumes a value differing from its original one by at most $h^{1/2}$. By considering only choices in which these vertical differences are multiples of $2K+1$, we obtain disjoint events, numbering on the order of $h^{3/2}$, each of which satisfies the bound (\ref{eqoutflbd}), and is contained in $\gacc \big( \bo{x},\bo{y},\delta  \big) \cap \sid \big( \bo{x},\bo{y},\delta  \big)$. This completes the proof. \qed
We now show that successful action of $\res_j$ gives rise to a favourable outcome in the associated sector. 
\begin{lemma}\label{lemsucfav}
Let $j \in \{1,\ldots,\mn \}$, and let $\omega \in \zoz$ satisfy $\areatrap \cap \big\{ \gamsw \cap B_{\ccone n} = \emptyset \big\}$. Fix $\chi > 0$ in (\ref{thtval}) at a sufficiently small value.  Then, if $\res_j$ acts successfully on $\zoz$,  its output is favourable in the sector $\secajn$. 
\end{lemma}
\noindent{\bf Proof.} Write $\omega' = \res_j(\omega)$.  Under 
$\sid \big( \bo{x_j}, \bo{y_j},\epso \big)$, there exists $\bo{z} \in \gamma_{\bo{x_j},\bo{y_j}}(\omega')$ such that   
$$
d \Big( \bo{z}, \partial {\rm conv} \big(  \big[ \bo{0}, \bo{x_j} \big]  \cup \big[ \bo{0}, \bo{y_j} \big] \cup \gamma_{\bo{x_j},\bo{y_j}}(\omega') \big) \Big) \geq \csid \vert\vert \bo{x_j} - \bo{y_j}  \vert\vert^{1/2} \big( \log  \vert\vert \bo{x_j} - \bo{y_j}  \vert\vert \big)^{1/2}.
$$
By Lemma \ref{leminv}, $\gamsw(\omega') = \big( \gamsw(\omega) \cap A_{\bo{x_j},\bo{y_j}}^c \big) \cup \gamma_{\bo{x_j},\bo{y_j}}(\omega')$.
Note that 
\begin{eqnarray}
d \Big(  \bo{z}  , \partial {\rm conv} \big( \gamsw(\omega') \big) \Big) 
   & = & d \Big(  \bo{z}  , \partial {\rm conv} \big( ( \gamsw(\omega) \cap A_{\bo{x_j},\bo{y_j}}^c ) \cup \gamma_{\bo{x_j},\bo{y_j}}(\omega') \big) \Big) \nonumber \\
 & = & 
\min \bigg\{ d \Big(  \bo{z}  , \partial {\rm conv} \big( ( \gamsw(\omega) \cap A_{\bo{x_j},\bo{y_j}}^c ) \cup \gamma_{\bo{x_j},\bo{y_j}}(\omega')  \big) \cap A_{\bo{x_j},\bo{y_j}} \Big), \nonumber \\
 & & \qquad \qquad d \Big(  \bo{z}  , \partial {\rm conv} \big( ( \gamsw(\omega) \cap A_{\bo{x_j},\bo{y_j}}^c ) \cup \gamma_{\bo{x_j},\bo{y_j}}(\omega')  \big) \cap A_{\bo{x_j},\bo{y_j}}^c \Big)
   \bigg\}. \nonumber
\end{eqnarray}
The first quantity in the minimum is at least $d \Big( \bo{z}, \partial {\rm conv} \big(  \big[ \bo{0}, \bo{x_j} \big]  \cup \big[ \bo{0}, \bo{y_j} \big] \cup \gamma_{\bo{x_j},\bo{y_j}}(\omega')  \big) \Big)$, because the replacement of  
$\big[ \bo{0}, \bo{x_j} \big]  \cup \big[ \bo{0}, \bo{y_j} \big]$ by 
 $\gamsw(\omega) \cap A_{\bo{x_j},\bo{y_j}}^c$ in ${\rm conv}\big( \cdot \cup \gamma_{\bo{x_j},\bo{y_j}}(\omega')  \big)$ cannot cause this set to become smaller. The second quantity is at least $d \big( \bo{z},  A_{\bo{x_j},\bo{y_j}}^c  \big)$ which exceeds  
 $\vert\vert \bo{x_j} - \bo{y_j}  \vert\vert \sin(\qzero)/3$, by assumption. 
Hence,
\begin{eqnarray}
   & &  d \Big(  \bo{z}  , \partial {\rm conv} \big( ( \gamsw(\omega) \cap A_{\bo{x_j},\bo{y_j}}^c ) \cup \gamma_{\bo{x_j},\bo{y_j}}(\omega')  \big) \Big) \nonumber \\
& \geq & \min \Big\{   \csid \vert\vert \bo{x_j} - \bo{y_j}  \vert\vert^{1/2} \big( \log  \vert\vert \bo{x_j} - \bo{y_j}  \vert\vert \big)^{1/2}  , \vert\vert \bo{x_j} - \bo{y_j}  \vert\vert  \sin(\qzero)/3 \Big\} \nonumber \\
 & \geq &  
    \frac{1}{2} \big( \frac{2\ccone}{3\pi} \big)^{1/2}  \epso  \ctheta^{1/2} n^{1/3} (\log n)^{2/3}, \nonumber
\end{eqnarray}
the second inequality by (\ref{distxjyj}). From $\functht = o\big(\ctheta^{1/2}\big)$, we know that 
$\functht <  (1/2) \big( 2\ccone/(3\pi) \big)^{1/2}  \epso  \ctheta^{1/2}$ holds if we choose  $\ctheta > 0$  to be small enough,
so that $\bo{z} \in \gamsw(\omega') \cap \secajn$
ensures that $\secajn$ is favourable under $\omega' = \res_j(\omega)$. \qed
\end{subsection}
\begin{subsection}{Local deviation cannot be undone by resampling distant sectors}
As we progressively apply the maps $\res_j$ in  counterclockwise order, we need to check that a favourable outcome returned by one map cannot be undone by a later one. For this purpose, it is convenient to impose that 
${\rm MFL}\big(\gamsw\big)$ is unusually high in neither the input nor the output configuration of $\res_j$.
\begin{lemma}\label{lemconvakpl}
Let $\csix  \in \{ 1,\ldots,\mn\}$, and let $\cfive  \in (0,\infty)$
satisfy $\cfive < \frac{\ctheta}{2} \csix$.
Set $\cfour = \csix + 1$.
Let $k,j \in \{ 1,\ldots, \mn \}$
satisfy $\vert j - k \vert   \geq \cfour$. Let $\omega \in \zoz$ satisfy $\areatrap$ and $\gamsw(\omega) 
\subseteq  B_{\ccone n}^c$. 
Set $\omega' = \res_j(\omega)$. Suppose that
$$
\max \Big\{ {\rm MFL} \big( \gamsw(\omega) \big),  {\rm MFL} \big( \gamsw(\omega') \big)  \Big\} \leq \cfive n^{2/3} \big( \log n \big)^{1/3}.
$$
Then
the sector $A_k^n$ is favourable under $\omega$ if and only if it is favourable under $\omega'$.
\end{lemma}
The following two lemmas are convenient tools for the proof of Lemma \ref{lemconvakpl}. The quantities $\csix$,$\cfive$ and $\cfour$ are fixed as in the statement of Lemma \ref{lemconvakpl}.
\begin{lemma}\label{lemverangsep}
Suppose that $\gamsw$ satisfies $\gamsw \cap B_{\ccone n} = \emptyset$. 
If a line segment $[\bo{u},\bo{v}]$ of $\partial {\rm conv} \big( \gamsw \big)$ intersects two sectors $A_k^n$ and $A_l^n$ for which $\vert k - l \vert \geq 2$, then $\vert\vert \bo{v} - \bo{u} \vert\vert \geq (\ctheta/2) \vert k - l \vert n^{2/3} \big( \log n \big)^{1/3}$.
\end{lemma}
\noindent{\bf Proof.} Since $\bo{u},\bo{v} \in \gamsw$ and $\gamsw \cap B_{\ccone n} = \emptyset$, we find that  $\vert\vert \bo{v} - \bo{u} \vert\vert \geq 2 \pi^{-1} \ang \big( \bo{u}, \bo{v} \big) \ccone n$. From $\bo{u} \in A_k^n$ and 
$\bo{v} \in A_l^n$, we find that $\ang \big( \bo{u}, \bo{v} \big) \geq \vert k - l  - 1 \vert 
\ctheta n^{-1/3} \big( \log n \big)^{1/3}$, as required. \qed
\begin{lemma}\label{lemconvak}
Let $k,j \in \{ 1,\ldots, \mn \}$
satisfy $\vert j - k   \vert \geq \csix$. Let $\omega \in \zoz$ satisfy $\areatrap$. Set $\omega' = \res_j(\omega)$ and further suppose that
$$
\max \Big\{ {\rm MFL} \big( \gamsw(\omega) \big),  {\rm MFL} \big( \gamsw(\omega') \big)  \Big\} \leq \cfive n^{2/3} \big( \log n \big)^{1/3}.
$$
Then  
$$
\partial {\rm conv}\big( \gamsw(\omega') \big) \cap A_k^n
= \partial {\rm conv}\big( \gamsw(\omega) \big) \cap A_k^n.
$$ 
\end{lemma}
\noindent{\bf Proof.}
We assume that $j \leq k - \csix$ without loss of generality. 
Let $\big[ \bo{w},\bo{w'} \big]$ (with $\argu(\bo{w}) < \argu(\bo{w'})$)
denote that edge in the polygon  
$\partial {\rm conv} \big( \gamsw(\omega) \big)$ that intersects the clockwise boundary of $A_k^n$.  
By Lemma \ref{lemverangsep}, the hypothesis on ${\rm MFL}\big( \gamsw(\omega) \big)$,  and   
$(\ctheta/2) \csix > \cfive$,
we have that $\bo{w} \in \cup_{l=j+1}^{k-1}A_l^n$.
If  
\begin{equation}\label{gamswjemp} 
\gamsw(\omega') \cap \secajn \cap \ell_{\bo{w},\bo{w'}} = \emptyset,
\end{equation}
then ${\rm conv} \big( \gamsw(\omega') \big)$ is unchanged from ${\rm conv} \big( \gamsw(\omega) \big)$ counterclockwise to $\bo{w}$,
so that the statement of the lemma holds.
If (\ref{gamswjemp}) does not hold, then  $\partial {\rm conv} \big( \gamsw(\omega') \big)$ contains a line segment 
$\big[ \bo{v}, \bo{v'} \big]$ with $\bo{v} \in \secajn$ 
and $\argu(\bo{v'}) \geq \argu(\bo{w'})$, so that $\bo{v'} \in A_{k'}^n$ for some $k' \geq k$. This circumstance is depicted in Figure 7.
\begin{figure}\label{figlemconvak}
\begin{center}
\includegraphics[width=0.6\textwidth]{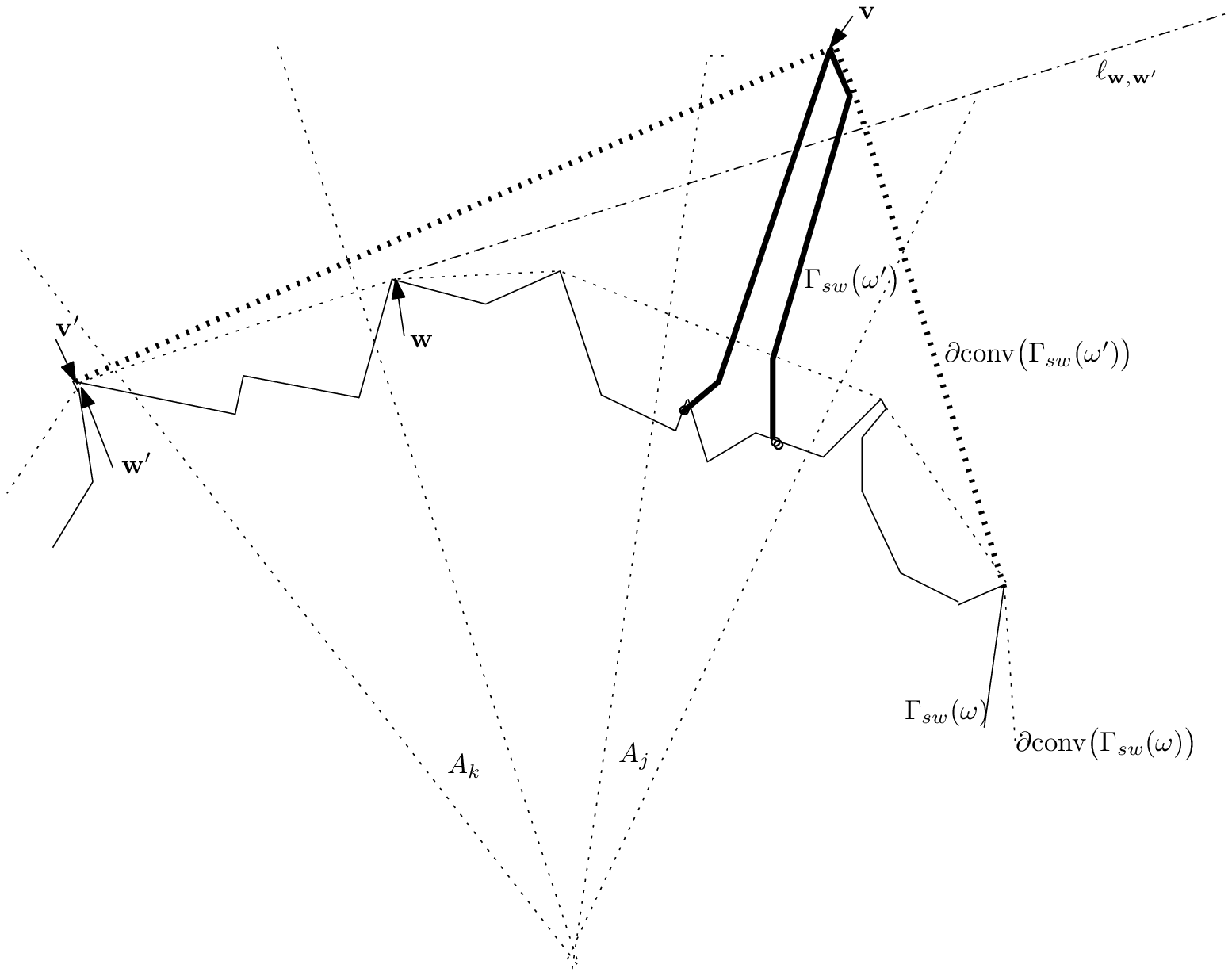} \\
\end{center}
\caption{An illustration of the case that 
$\gamsw(\omega') \cap \secajn \cap \ell_{\bo{w},\bo{w'}} \not= \emptyset$ 
in the proof of Lemma \ref{lemconvak}. The bold resampled circuit segment in $\secajn$ cuts across $\ell_{\bo{w},\bo{w'}}$ to create a modified convex boundary (shown as bold dotted) that contains a line segment $\big[ \bo{v},\bo{v'} \big]$ intersecting both $\secajn$ and $A_k^n$.}
\end{figure}
By Lemma \ref{lemverangsep} and  $(\ctheta/2) \csix > \cfive$, we obtain ${\rm MFL} \big(  \gamsw(\omega') \big)
\geq \dist \bo{v} - \bo{v'} \dist 
 > \cfive n^{2/3} \big( \log n \big)^{1/3}$, which contradicts the hypothesis. \qed
\noindent{\bf Proof of Lemma \ref{lemconvakpl}.} 
By Lemma \ref{lemconvak},
\begin{equation}\label{convaeq}
  \partial {\rm conv} \Big( \gamsw(\omega') \Big) \cap \Big( A_{k-1}^n \cup A_k^n \cup A_{k+1}^n \Big) = 
    \partial {\rm conv} \Big( \gamsw(\omega) \Big) \cap \Big( A_{k-1}^n \cup A_k^n \cup A_{k+1}^n \Big).
\end{equation}
Suppose that $A_k^n$ is not favourable under $\omega$. We wish to show that any $\bo{v} \in A_k^n \cap  \gamsw(\omega')$ satisfies   ${\rm LR} \big( \bo{v} , \gamsw(\omega') \big) \leq \functht  n^{1/3} \big( \log n \big)^{2/3}$.
Let $\bo{v'} \in  \partial {\rm conv} \big( \gamsw(\omega) \big)$ satisfy $d\big( \bo{v'} , \bo{v} \big) = {\rm LR} \big( \bo{v} , \gamsw(\omega) \big)$. We have that $\bo{v'} \in  A_{k-1}^n \cup A_k^n \cup A_{k+1}^n$, for otherwise
\begin{equation}\label{vprv}
d\big( \bo{v'} , \bo{v} \big) \geq \vert\vert \bo{v} \vert\vert
 \sin \ang \big( \bo{v'},\bo{v} \big) \geq \vert\vert \bo{v} \vert\vert
 \sin \Big(  \ctheta n^{-1/3} \big( \log n \big)^{1/3} \Big)
 \geq 2 \pi^{-1} \ccone \ctheta   n^{2/3} \big( \log n \big)^{1/3}, 
\end{equation}
(the third inequality by $\gamsw(\omega) \subseteq B_{\ccone n}^c$). 
This bound is impossible, since 
$$
d\big( \bo{v'} , \bo{v} \big) =  {\rm LR} \big( \bo{v} , \gamsw(\omega) \big)
 \leq \functht  n^{1/3} \big( \log n \big)^{2/3},$$ 
due to $A_k^n$ being favourable under $\omega$. 
By (\ref{convaeq}), $\bo{v'} \in   \partial {\rm conv} \big( \gamsw(\omega') \big)$, so that, since $\bo{v} \in \gamsw(\omega')$, we indeed have ${\rm LR} \big( \bo{v} , \gamsw(\omega') \big) \leq 
d\big( \bo{v'} , \bo{v} \big) \leq \functht n^{1/3}  \big( \log n \big)^{2/3}$.
We find then that  $A_k^n$ is not favourable under $\omega'$.

Suppose that  $A_k^n$ is favourable under $\omega$.
Let $\bo{v} \in A_k^n \cap \gamsw(\omega)$ satisfy ${\rm LR} \big( \bo{v} , \gamsw(\omega) \big) \geq  \functht  n^{1/3} \big( \log n \big)^{2/3}$.
Note that  $\bo{v} \in \gamsw(\omega')$, because $\gamsw(\omega')$ coincides with $\gamsw(\omega)$ in the complement of $\secajn$. We wish to show that  ${\rm LR} \big( \bo{v} , \gamsw(\omega') \big) \geq  \functht  n^{1/3} \big( \log n \big)^{2/3}$, because this will demonstrate that  $A_k^n$ is favourable under $\omega'$. Let 
$\bo{v'} \in  \partial {\rm conv} \big( \gamsw(\omega') \big)$ satisfy $d\big( \bo{v'} , \bo{v} \big) = {\rm LR} \big( \bo{v} , \gamsw(\omega') \big)$. 
If $\bo{v'} \not\in  A_{k-1}^n \cup A_k^n \cup A_{k+1}^n$, then (\ref{vprv}) certainly shows that 
$d\big( \bo{v'} , \bo{v} \big)$ is greater than $\functht  n^{1/3} \big( \log n \big)^{2/3}$. In the other case,  
$\bo{v'} \in   \partial {\rm conv} \big( \gamsw(\omega) \big)$ by (\ref{convaeq}), so that $d\big( \bo{v'} , \bo{v} \big) \geq {\rm LR} \big( \bo{v} , \gamsw(\omega) \big)$, which we know to be at least $\functht  n^{1/3} \big( \log n \big)^{2/3}$. Thus,   $A_k^n$ is indeed favourable under $\omega'$. \qed
\end{subsection}
\begin{subsection}{Considerations for defining the complete resampling procedure}
The naive approach would be to define the overall resampling procedure \hfff{res} $\res$ by iteratively applying each of the maps $\res_j$ in a counterclockwise order. However, we have already noted that we must insist that the sectors $\secajn$ be in the right-hand half-plane. There are two further problems. Firstly, favourable outcomes in a sector are known by Lemma \ref{lemconvakpl} to be locked in for the long term only if $\cfour$ further sectors are treated without destroying the favourable outcome. 
Secondly, Lemma \ref{lemactsuc} provides a lower bound on a favourable outcome in a sector only if the sector has a low local roughness in the input. We do not control the outcome if the maximal local roughness in the sector is higher, and, even though this is the favourable property that we seek, it could, 
in principle, be consistently destroyed by the resampling. One might object that, if maximum local roughness is high enough in the input, we already have the desired conclusion, making a random resampling, or, indeed, any other argument, unnecessary. However, a problematic case is one in which, given 
$\areatrap$,
with fifty percent probability, the maximal local roughness in sector $\secajn$ exceeds $\functht n^{1/3} (\log n)^{2/3}$ for all $j$, 
and, in the rest of the conditional probability space, $ {\rm MLR}(\gamsw) \leq \functht n^{1/3} (\log n)^{2/3}$. Under a resampling of every sector successively,
the first part of the space might be mapped into the second (and {\it vice versa}), however far-fetched this may seem. Our solution to this difficulty is to define $\res$
by attempting to resample sectors only with a positive probability, so that sectors already realizing the lower bound on 
${\rm MLR}\big(\gamsw\big)$ 
may be preserved in the final output.  
\end{subsection}
\begin{subsection}{The proof of Theorem \ref{thmmlrubd}}\label{secthmpf}
We are now ready to formally define our complete resampling procedure $\res$ and to analyse its action.

Let $\mn'$  denote the number of $j$, $1 \leq j \leq \mn$, such that $\secajn$ lies in the first quadrant. Note that, by (\ref{mtilde}), 
\begin{equation}\label{mnpr}
\mn' \sim \frac{\pi}{2\ctheta} n^{1/3} \big(\log n  \big)^{-1/3}. 
\end{equation}
It is this initial segment of the sequence of sectors that we will use to define $\res$.  
The operation $\res$ will be defined in terms of the parameters $\qone,\qtwo$ and $\qthr$ from Lemma \ref{lemconvakpl}, with their form being specified later. 

Let $\big( \Omega, \mathcal{F}, \mathbb{P} \big)$ denote a probability space on which are defined 
a random configuration $\phi$ having the law $P\big(\cdot \big\vert \areatrap \big)$, and the operations $\res_j$, $1 \leq j \leq \mn'$, each of these acting independently. (The use of the notation $\mathbb{P}$ extends its original use in Definition \ref{defres}.) A further independent sequence of length $\mn'$ of biased coins with probability $1/\cfour$ of landing heads is provided under  $\big( \Omega, \mathcal{F},\mathbb{P} \big)$.
 The operation $\res$ is defined under $\big( \Omega, \mathcal{F},\mathbb{P} \big)$ as follows. 
The input $\omega_0$ is taken to be $\phi$.
The operation has  $\mn'$ further stages. At the $j$-th stage, $1 \leq j \leq \mn'$, with probability $1/\cfour$ (using the $j$-th coin), procedure $\res_j$ acts on the present configuration, so that the configuration $\omega_j$ at the end of stage $j$ is given by $\omega_j = \res_j(\omega_{j-1})$. In the other case, we take $\omega_j = \omega_{j-1}$, and say that no action is taken. 
Fixing notation, let $\gamsw(j) = \gamsw(\omega_j)$ and $\gamswcl(j) = \gamswcl(\omega_j)$ . Set $\smtu_j = \smtu(\omega_j)$ and $\unfav_j =  \unfav(\omega_j)$. 

By Lemma \ref{lemresjinv}, the law  $P \big( \cdot \big\vert \areatrap \big)$ is invariant under $\res$. We will analyse the measure  $P \big( \cdot \big\vert \areatrap \big)$ by considering it as the law of the output $\omega_{\mn'}$ of $\res$ under the measure $\mathbb{P}$. 

The analysis of the action will be made under the condition that the realization of $\res$ lies in a 
space $\mathcal{G}$ of ``good'' outcomes. To define this,
for $0 \leq i \leq \mn'$, set 
$$
\mathcal{G}_{1,i} = 
\Big\{ {\rm MFL} \big( \gamsw(i) \big) \leq \cfive     n^{2/3} \big(\log n  \big)^{1/3} \Big\},
$$
$$
  \mathcal{G}_{2,i} = \Big\{ {\rm MLR} \big(  \gamsw(i) \big) \leq n^{2/3}   \Big\},
$$
$$
\mathcal{G}_{3,i} =  \Big\{ \theta_{{\rm RG}_{\qzero/2,c_0/2}}^{\rm MAX} \big( \gamswcl(i) \big)  \leq  n^{\epsilon - 1}/2  \Big\},
$$
and
$$
\mathcal{G}_{4,i} = \Big\{  \gamsw(i) \subseteq B_{\cctwo n} \setminus  B_{\ccone n}   \Big\}
  \cap \Big\{ \bo{0} \in {\rm INT} \big( \gamsw(i) \big) \Big\}.
$$
(Recall that $\epsilon \in (0,2/3)$ was introduced in Definition \ref{defres}.)
We also set $\mathcal{G}_i = \cap_{j=1}^4 \mathcal{G}_{j,i}$,
and $\mathcal{G}_{(i)} = \cap_{j=0}^i \mathcal{G}_j$. We then define $\mathcal{G} = \mathcal{G}_{(\mn')}$.
We make a definition of the set of good outcomes so that this event is probable: 
with $\epsone \in (0,2/3)$, we write (with integer-rounding omitted) 
\begin{equation}\label{qchoose}
\qthr = n^{\epsone} \qquad \textrm{and} \qquad  \qtwo = \frac{\chi}{4} n^{\epsone},
\end{equation}
so that these quantities  satisfy the hypotheses stated in Lemma \ref{lemconvakpl}. As we derive after completing the rest of the proof, if $\epsilon \in (0,1/6)$, then 
\begin{equation}\label{gbound}
 \mathbb{P} \big( \mathcal{G}^c \big) \leq \exp \Big\{ - c n^{\frac{3\epsone}{2}} \log n  \Big\} + \exp \Big\{ - c n^{\epsilon} \Big\}.
\end{equation}
\begin{lemma}\label{lemkh}
For $1 \leq j \leq \mn' - \cfour$, the occurrence of $\mathcal{G}$ implies that
\begin{equation}\label{statk}
  \unfav_j \cap \Big\{ j+ \cfour, \ldots, \mn' \Big\} =   \unfav_0 \cap \Big\{ j+ \cfour, \ldots, \mn' \Big\}
\end{equation}
and that
\begin{equation}\label{stath}
  \smtu_j \cap \Big\{ j+ \cfour, \ldots, \mn' \Big\} =   \smtu_0 \cap \Big\{ j+ \cfour, \ldots, \mn' \Big\}.
\end{equation}
\end{lemma}
\noindent{\bf Proof.} 
Let the pair $(j,k)$, $j,k \in \{ 1,\ldots,\mn'\}$ satisfy $j+ \cfour \leq k \leq \mn'$. 
The occurrence of $\mathcal{G}$ means that Lemma \ref{lemconvakpl} may be applied to each of the first $j$ stages of the formation of $\res$. It tells us that $A_k^n$ is favourable under $\omega_j$ if and only if it is favourable under $\omega_0$. This is (\ref{statk}).
Note that, for such  $(j,k)$, the condition 
$k \in \smtu_j$ 
is determined by the data 
$\partial {\rm conv} \big( \gamsw(j) \big) \cap A_k^n$. 
However, 
$\partial {\rm conv} \big( \gamsw(j) \big) \cap A_k^n 
=\partial {\rm conv} \big( \gamsw(0) \big) \cap A_k^n$,
because Lemma \ref{lemconvak} is also applicable to each of  the first $j$ stages of the formation of $\res$. Thus, we obtain (\ref{stath}). \qed

Recall from Lemma \ref{lemgbd} that,  if $\gamsw(0) \subseteq B_{\cctwo n}$ and 
$\bo{0} \in  {\rm INT} \big( \gamsw(0) \big)$, 
then  $\big\vert \smtu_0 \big\vert \geq 9\mn/10$, for large $n$. For such an input,  
we also have that $\big\vert \smtu_0 \cap \big\{ 1,\ldots,\mn'\big\} \big\vert \geq \mn'/2$,
since $\mn' \geq (\mn/4) \big( 1 + o(1) \big)$. By the occurrence of $\mathcal{G}_{4,0}$, we may thus find 
define a subset
$\smtup \subseteq \smtu_0 \cap \big\{ 1,\ldots,\mn' \big\}$ 
with the following properties: $\big\vert \smtup \big\vert \geq \frac{\mn'}{4\cfour}$, and each pair of consecutive elements of $\smtup$ differ by at least $2\cfour+1$. 
We enumerate the two sets $\smtup \cap \unfav_0 = \big\{ p_1,\ldots, p_{r_1} \big\}$ and  $\smtup \cap \unfav_0^c = \big\{ q_1,\ldots, q_{r_2} \big\}$, where $r_1 + r_2 = \vert \smtup \vert$. 
For $1 \leq r \leq r_1$, let $P_r$ denote the event that, in the action of $\res$, at stage $p_r$, $\res_{p_r}$ is chosen to act, and that it acts successfully, while at the $\cfour$ stages preceding, and at the $\cfour$ stages following, the $p_r$-th stage, no action is taken. 
For $1 \leq r \leq r_2$, let $Q_r$ denote the event that, in the action of $\res$, at stage $q_r$, $\res_{q_r}$ does not act, and neither does it do so at the $\cfour$ stages preceding, and the $\cfour$ stages following, the $q_r$-th stage.

We claim that:\\
\noindent{\bf Claim A.} For each $r \in \{ 1,\ldots,r_1 \}$, $\mathcal{G} \cap P_r$ entails $\res$ returning an output $\omega_{\mn'}$ under which  $A_{p_r}^n$ is favourable.\\
\noindent{\bf Claim B.}  For $r \in \{ 1,\ldots,r_2 \}$, $\mathcal{G} \cap Q_r$ entails $\res$ returning an output $\omega_{\mn'}$  under which  $A_{q_r}^n$ is favourable.\\
\noindent{\bf Proof of Claim A.}
Note that, by Lemma \ref{lemkh}, $p_r \in \unfav_{p_r - \cfour} \cap \smtu_{p_r - \cfour}$, since $p_r \in \smtu_0 \cap \unfav_0$. Given $P_r$, then, we have that $p_r \in \unfav_{p_r - 1} \cap \smtu_{p_r -1}$, because $\omega_{p_r - 1} = \omega_{p_r - \cfour}$.  
Applying Lemma \ref{lemsucfav} to the successful action of $\res_{p_r}$ on $\omega_{p_r - 1}$, we find that $A_{p_r}^n$ is favourable under $\omega_{p_r}$. 
On the event $P_r$, this remains the case under $\omega_{p_r + \cfour}$, because $\omega_{p_r + \cfour} = \omega_{p_r}$. That $A_{p_r}^n$ is now favourable in the configuration is now permanent: indeed, 
we may show that $A_{p_r}^n$ is favourable under $\omega_j$ for $p_r + \cfour \leq j \leq \mn'$
inductively. The $j$-indexed inductive step is trivial if $\res_j$ does not act. If it does act, then Lemma \ref{lemconvakpl} applies to  this action to give the inductive step.
Lemma \ref{lemconvakpl} is applicable because $\big\vert {\rm INT}\big( \gamsw(j) \big) \big\vert \geq n^2$ holds for $j = 0$ by assumption on $\omega_0$, and, for $j > 0$, by the construction of $\res_j$, while the other hypotheses are satisfied by the occurrence of $\mathcal{G}$. \qed
\noindent{\bf Proof of Claim B.}
As in the preceding proof, we have that $q_r \not\in \unfav_{q_r - \cfour}$. The inaction of $\res_j$ in stages between $q_r - \cfour$ and $q_r + \cfour$ means that  $q_r \not\in \unfav_{q_r + \cfour}$. Lemma \ref{lemconvakpl} is then applicable as in the earlier case. \qed
Claims A and B imply that 
\begin{equation}\label{eqronertwo}
   \Big(  \bigcup_{i=1}^{r_1} P_i \, \cup \, \bigcup_{i=1}^{r_2} Q_i  \Big) \cap \mathcal{G}
  \subseteq 
\Big\{   {\rm MLR} \big( \gamsw(\mn') \big) \geq \functht n^{1/3} \big( \log n \big)^{2/3} \Big\}
\end{equation}
Note then that
\begin{eqnarray}
 & &  \Big(  \bigcup_{i=1}^{r_1} P_i \, \cup \, \bigcup_{i=1}^{r_2} Q_i  \Big)^c \cap \mathcal{G} \label{incpqg} \\
 & \subseteq &
   \bigcap_{i=1}^{\frac{\mn'}{8 \cfour}} \Big( 
\big\{  r_1 \geq i \big\} \cap
P_i^c \cap \mathcal{G}_{(p_i - \cfour)} \Big)   
  \, \cup \,
 \bigcap_{i=1}^{\frac{\mn'}{8 \cfour}} \Big( 
\big\{  r_2 \geq i \big\} \cap Q_i^c \cap \mathcal{G}_{(q_i - \cfour)} \Big)  \nonumber
\end{eqnarray}
since, as we have noted, $\mathcal{G}$ implies that $r_1 + r_2 \geq \mn'/(4\cfour)$.

We claim that, for any $K \in \N$,
given $\big\{ r_1 \geq K \big\} \cap \mathcal{G}_{(p_K - \cfour)}$ and the values of $\mathbf{1}_{P_1},\ldots,\mathbf{1}_{P_{K-1}}$, the conditional probability that $P_K$ occurs is at least 
$(4\cctwo)^{-2}  \sin^2 \big( \qzero/2 \big)  n^{-2\epsilon - \clemgac \epso^2} \cfour^{-1} \big( 1 - 1/\cfour \big)^{2\cfour}$.
Indeed, the event on which we condition here is measurable with respect to $\big\{ \omega_0,\ldots,\omega_{p_K - \qthr} \big\}$, and, if it occurs, then $\omega_{p_K - \qthr}$
satisfies the event in the hypothesis of Lemma \ref{lemactsuc}. The claim then follows by this lemma, since, similarly to the proof of Claim A, we have that 
$p_K \in \unfav_{p_K - \qthr} \cap \smtu_{p_K - \qthr}$ under the event on which we condition, with 
$p_K \in \unfav_{p_K - 1} \cap \smtu_{p_K - 1}$ if no action is taken at stages $p_K - \qthr,\ldots, p_K - 1$.
From this claim, we find that
\begin{eqnarray}
  & & \mathbb{P} \Big( 
 \bigcap_{i=1}^{\frac{\mn'}{8 \cfour}}  
\big\{  r_1 \geq i \big\} \cap
P_i^c \cap \mathcal{G}_{(p_i - \cfour)} \Big) \label{pgineq} \\ 
& \leq & 
\Big( 1 - (4\cctwo)^{-2}  \sin^2 \big( \qzero/2 \big)  n^{-2\epsilon - \clemgac \epso^2} \cfour^{-1} \big( 1 - 1/\cfour \big)^{2\cfour}
 \Big)^{\frac{\mn'}{8\cfour}}. \nonumber
\end{eqnarray}
A similar assertion for the sequence of events $\big\{ Q_i \big\}$ yields
\begin{equation}\label{qgineq}
  \mathbb{P} \Big( 
\bigcap_{i=1}^{\frac{\mn'}{8 \cfour}} 
\big\{  r_2 \geq i \big\} \cap Q_i^c \cap \mathcal{G}_{(q_i - \cfour)} \Big) 
\leq
\Big( 1 -  \big( 1 - 1/\cfour \big)^{2\cfour + 1} \Big)^{\frac{\mn'}{8\cfour}}.
\end{equation}
Recall that the entire resampling procedure, whose law we indicate 
by $\mathbb{P}$, begins with a copy of $P \big( \cdot \big\vert \areatrap \big)$ as the input $\omega_0$, and ends up, after its numerous attempted modifications, spitting out another copy of $P \big( \cdot \big\vert \areatrap \big)$. Bearing this in mind,
we now reach the equation summarizing our approach, of viewing the conditional measure 
$P \big( \cdot \big\vert \areatrap \big)$ as the output of the complete resampling $\res$:
\begin{eqnarray}
 & & P \Big( {\rm MLR} \big( \gamsw \big) \geq  
\functht  n^{1/3} \big( \log n \big)^{2/3} \Big\vert \areatrap \Big) \label{mlrlowbd} \\
 & = & \mathbb{P} \Big( {\rm MLR} \big( \gamsw(\mn') \big) \geq  
\functht n^{1/3} \big( \log n \big)^{2/3} \Big) \nonumber \\ 
 & \geq & 1 - \mathbb{P} \big( \mathcal{G}^c \big)
   - \mathbb{P} \bigg( \Big(  \bigcup_{i=1}^{r_1} P_i \, \cup \, \bigcup_{i=1}^{r_2} Q_i  \Big)^c \cap \mathcal{G} \bigg) \nonumber 
\end{eqnarray} 
The inequality is due to (\ref{eqronertwo}).


By using
(\ref{gbound}), (\ref{incpqg}), (\ref{pgineq}), (\ref{qgineq}), (\ref{mnpr}) and (\ref{qchoose}),
\begin{eqnarray}
 & & P \Big( {\rm MLR} \big( \gamsw \big) <  
\functht  n^{1/3} \big( \log n \big)^{2/3} \Big\vert \areatrap \Big) \label{mlrareatrap} \\
  & \leq &  \exp \Big\{ - c n^{\frac{3\epsone}{2}} \log n \Big\}
   + \exp \big\{ - c n^{\epsilon} \big\} \nonumber \\
 & & +
   \exp \Big\{ - c n^{1/3 - 2\epsilon - 2 \epsone - \clemgac \epso^2} \big( \log n \big)^{-1/3} \Big\}
  +  \exp \Big\{ - c n^{1/3 -  \epsone} \big( \log n \big)^{-1/3} \Big\} \nonumber 
\end{eqnarray}
The minimum 
$$
\min \Big\{  \frac{3\epsone}{2}, \epsilon, 1/3 - 2 \epsilon - 2 \epsone - \clemgac \epso^2 \Big\}
$$
is attained by choosing $\epsilon = \frac{1}{13} - \frac{3 \clemgac}{13} \epso^2$
and $\epsone = \frac{2\epsilon}{3}$. 
The quantity $\epso > 0$, introduced in Definition \ref{defactsuc}, may be fixed at an arbitrarily small value.
Hence, for any $\epsilon > 0$, and for $n \in \N$ sufficiently high,
\begin{equation}\label{mlrareatrapnew}
 P \Big( {\rm MLR} \big( \gamsw \big) <  
\functht  n^{1/3} \big( \log n \big)^{2/3} \Big\vert \areatrap \Big) 
   \leq  \exp \Big\{ -n^{\frac{1}{13} - \epsilon}  \Big\}.
\end{equation}
Recall that $\functht > 0$ is a fixed positive number, specified by making a sufficiently small choice of $\ctheta > 0$ (in terms of $\epso$). In this way, the statement of Theorem \ref{thmmlrubd}
is the analogue of (\ref{mlrareatrapnew}) for the circuit $\cir$ under the conditional measure
$P \big( \cdot \big\vert \areacon \big)$. We obtain this statement from (\ref{mlrareatrap}) by applying Lemma \ref{lemlas} with the choice $\mathcal{M} = \big\{ \Gamma: {\rm MLR}\big( \Gamma\big) < \functht 
  n^{1/3}  \big( \log n \big)^{2/3} \big\}$.

It remains to establish (\ref{gbound}). Note that 
\begin{equation}\label{gmnbd}
 \mathbb{P} \big( \mathcal{G}^c \big)  \leq  \sum_{i=0}^{\mn'}
  \sum_{j=1}^4 
\mathbb{P} \big( \mathcal{G}_{j,i}^c \big)  
  =  \big( \mn' + 1 \big)   \sum_{j=1}^4 
\mathbb{P} \big( \mathcal{G}_{j,0}^c \big) 
\end{equation}
the equality due to Lemma \ref{lemresjinv}. Note that
$$
\mathbb{P} \big( \mathcal{G}_{1,0}^c \big) 
 = P \Big(    {\rm MFL} \big( \gamsw \big) > \cfive     n^{2/3} \big(\log n  \big)^{1/3}   \Big\vert \areatrap  \Big).
$$
The probability of the analogous event for $\cir$ under $P \big( \cdot \big\vert \areacon \big)$
appears in \cite{hammondone} and has been quoted as Theorem \ref{thmmflbd}.
By setting $\mathcal{M} = \big\{ \Gamma: {\rm MFL}(\Gamma) > t n^{2/3} \big(\log n  \big)^{1/3}  \big\}$,
we may apply Lemma \ref{lemoswtrans} to Theorem \ref{thmmflbd}
to learn that the statement of this theorem holds for the circuit $\gamsw$
under the measure $P \big( \cdot \big\vert \areatrap \big)$. 

We also have that
$$
\mathbb{P} \big( \mathcal{G}_{2,0}^c \big) 
 = P \Big(    {\rm MLR} \big( \gamsw \big) \leq n^{2/3}  \Big\vert \areatrap  \Big)  \leq \exp \big\{ - n^{1/6} \big( \log n \big)^{-C} \big\}.
$$
In this case, we apply Lemma \ref{lemoswtrans} to Theorem \ref{thmmlrbd} with the maximum possible choice $t = O\big( n^{5/36} (\log n)^{-C} \big)$.

We have that 
$$
\mathbb{P} \big( \mathcal{G}_{3,0}^c \big) 
 = P \Big(   \marcl >  n^{\epsilon - 1}     \Big\vert \areatrap  \Big).
$$
For any choice of $\epsilon \in (0,1)$, Theorem \ref{thmmaxrgclgw} yields
$$
\mathbb{P} \big( \mathcal{G}_{3,0}^c \big)  \leq \exp \big\{ - c n^{\epsilon} \big\}.
$$
Finally, 
$$
\mathbb{P} \big( \mathcal{G}_{4,0}^c \big) 
 = P \Big(  \Big\{ \gamsw \not\subseteq B_{\cctwo n} \setminus B_{\ccone n} \Big\} 
    \cup
 \Big\{
\bo{0} \not\in {\rm INT} \big( \gamsw \big)  \Big\}
  \Big\vert \areatrap  \Big)
$$
is bounded in Lemma \ref{lemmacat}. Substituting these bounds into (\ref{gmnbd}), and using (\ref{mnpr}), 
we obtain that 
$\mathbb{P} \big( \mathcal{G}^c \big)$ is at most
$$
C  n^{1/3} \big(\log n  \big)^{-1/3}
  \Big(  \exp \big\{ - c \cfive^{3/2} \log n \big\} + 
\exp \big\{ - n^{1/6} \big( \log n \big)^{-C} \big\} +
\exp \big\{ - c n^{\epsilon} \big\} + 
\exp \big\{ - c' n \big\} \Big).
$$
From $\qtwo = n^{\epsone}$ with $\epsone > 0$ and $\epsilon \in (0,1/6)$,
we obtain (\ref{gbound}). \qed
\end{subsection}
\begin{subsection}{The lower bound on maximum facet length}
This is, in essence, a straightforward corollary of Theorem \ref{thmmlrubd}. We discuss the argument in outline before reaching the formal proof. A more detailed heuristic discussion of the same ideas may be found in Section 1.2.2 of \cite{hammondone}. Given the lower bound $\mlr \geq c n^{1/3} \big( \log n \big)^{2/3}$, the occurrence of $\mfl = o \big( n^{2/3} \big( \log n \big)^{1/3} \big)$ would force the existence of a facet in the convex boundary $\delconv$ of length 
$o \big( n^{2/3} \big( \log n \big)^{1/3} \big)$ but for which the subpath $\sentier$ of $\cir$ running from one endpoint of this facet to the other
deviates from the facet by at least $c n^{1/3} \big( \log n \big)^{2/3}$. (To find such a facet, we consider that facet for which the corresponding section of the circuit contains the point attaining maximum local roughness.) The subpath $\sentier$ is a subcritical open path that deviates from the line segment interpolating its endpoints to such a degree that this fluctuation lies in a superpolynomially small part of the probability space of a point-to-point connection between these endpoints. (Such connections have Gaussian fluctuation \cite{civ}.) As such, it is probabilistically cheaper to form an alternative subcritical path from one endpoint of the facet to another, because an open path between the two points with a typical fluctuation will trap enough area.

Formally, the apparatus required to give the proof of Theorem \ref{thmmflubd} has been set up in \cite{hammondone}.
The proof of Theorem 1.1 of \cite{hammondone} uses the information that $\mfl$ cannot be too large to deduce that $\mlr$ is also not too large. Equally, it may be used to show that, if $\mfl$ has even a slight probability of being small, then so does $\mlr$.  \\
\noindent{\bf Proof of Theorem \ref{thmmflubd}}
This is a matter of modifying a few details in the proof of Theorem 1.1  of \cite{hammondone}, that appears in Section 5. 
Consider the inclusion \cite{hammondone}:(5.67), with the appearance of $t^{2 - \delta}$ on the right-hand-side being replaced by $t^{2 + \delta}$, and  where now we take $t > 0$ to be small.  The analogue of \cite{hammondone}:(5.69) is the assertion that there exists $t_0 > 0$ such that, for $0 < t < t_0$, for which \cite{hammondone}:(5.69) 
holds, with $\delta$ replaced by $-\delta$.  Applying the analogue of \cite{hammondone}:(5.69), the inequality $\mlrs \leq \mfl$, the deduction made after \cite{hammondone}:(5.80) that $A = \emptyset$,  \cite{hammondone}:(5.82)
and 
Lemma \ref{lemmac}
to the analogue of \cite{hammondone}:(5.67),
we obtain
\begin{eqnarray}
 & &  P \Big(  \mlr \geq n^{1/3} \big( \log n \big)^{2/3} t , \acon \Big) \label{mlrinte} \\
 & \leq & 
 P \Big(  \mfl \geq n^{1/3} \big( \log n \big)^{2/3} t^{2 + \delta} , \acon \Big) \nonumber \\
  & & + \bigg( n^{2/3} \big( \log n \big)^{1/3} t^{2 + \delta} \exp \Big\{ -
  c  \min \big\{  t^{-\delta} \log n, n^{1/3} \big\} \Big\}
+ \exp \Big\{  - c n^{1/6} \Big\}  \bigg)   P \Big( \acon \Big).  \nonumber
\end{eqnarray}
At this stage, we require the version of Theorem \ref{thmmlrubd} for the measure $P \big( \cdot \big\vert \acon \big)$. This version follows from Theorem \ref{thmmlrubd} and the bound
\begin{equation}\label{eqcomp}
P \Big( \acon \Big\vert \areacon \Big) \geq (1/2) \pi^{-1} \cctwo^{-2} n^{-2}.
\end{equation}
 This last bound is valid, because, as Proposition \ref{propglobdis} states, $\centre\big( \cir \big) \in {\rm INT} \big( \cir \big)$ is highly probable under $P \big( \cdot \big\vert \areacon \big)$; moreover, $\cir \subseteq B_{\cctwo n}$ is also probable, by Lemma \ref{lemmac}. These have the consequence that a translate of a typical realization of $P \big( \cdot \big\vert \areacon \big)$ by some vector $\bo{v} \in \Z^2 \cap B_{\cctwo n}$ will realize $\acon$, whence (\ref{eqcomp}).

Setting $t$ in (\ref{mlrinte}) to be sufficiently small, dividing this inequality by $P \big( \acon \big)$ and applying our new version of 
Theorem \ref{thmmlrubd}, we find that the left-hand-side of the resulting inequality is close to one. 
In this way, we obtain the statement of the theorem 
under the conditional measure $P \big( \cdot \big\vert \acon \big)$. Lemma 2.8 (of \cite{hammondone}), that translates upper bounds for events under $P \big( \cdot \big\vert \areacon \big)$ in terms of the corresponding quantity under $P \big( \cdot \big\vert \acon \big)$,  provides the means to obtain the theorem. \qed
\end{subsection}
\end{section}
\begin{section}{The technical tools}\label{sectools}
The proofs of several elements of our approach have been deferred and appear in this section.
\begin{subsection}{The circuits after resampling}\label{sectoolsone}
\noindent{\bf Proof of Lemma \ref{leminv}.}
Set $B_{\bo{x},\bo{y}} \subseteq \zoz$, $B_{\bo{x},\bo{y}} = \areatrap \cap \big\{ \bo{x}, \bo{y} \in {\rm RG} \big( \gamswcl \big) \big\} \cap 
\big\{ \bo{0} \in {\rm INT} \big( \gamsw \big) \big\}$. We must show that 
the law $P\big( \cdot \big\vert B_{\bo{x},\bo{y}} \big)$
is invariant under $\psi_{\bo{x},\bo{y}}$. 
To this end, let $\omega \in  B_{\bo{x},\bo{y}}$, and write $\omega' = \psi_{\bo{x},\bo{y}}(\omega)$. 
Recall that $\overline\gamma_{\bo{x},\bo{y}}(\omega)$ denote
the common $\omega-$open cluster of $\bo{x}$ and $\bo{y}$ in $\axy$,
and  $\gamma_{\bo{x},\bo{y}}(\omega)$
the outermost $\omega-$open path from $\bo{x}$ to $\bo{y}$ in $\axy$.
We also use these notations for the configuration $\omega'$. 
To derive
\begin{equation}\label{gsomtwo}
 \gamsw \big( \omega' \big)
 = \Big(  \gamsw \big( \omega \big)\cap \axy^c  \Big) \cup \gamma_{\bo{x},\bo{y}}(\omega'),
\end{equation}
we will begin by showing that
\begin{equation}\label{gsomone}
 \gamswcl \big( \omega' \big)
 = \Big(  \gamswcl \big( \omega \big)\cap \axy^c  \Big) \cup 
   \overline\gamma_{\bo{x},\bo{y}}(\omega').
\end{equation}
To this end, note that $\bo{x},\bo{y} \in {\rm RG} \big( \gamswcl(\omega)\big)$ implies that
\begin{equation}\label{lplusxy}
\ell_{\bo{0},\bo{z}}^+ \cap  \gamswcl(\omega) = \big\{ \bo{z} \big\}
\qquad \textrm{for $\bo{z} = \bo{x},\bo{y}$.}
\end{equation}
From $\gamsw(\omega) \subseteq \gamswcl(\omega)$, (\ref{lplusxy}) and 
$\bo{0} \in {\rm INT}\big( \gamsw(\omega) \big)$, we find that 
 $\gamsw(\omega) \cap A_{\bo{x},\bo{y}}^c$ is a path from $\bo{x}$ to $\bo{y}$, and hence is connected. By this, (\ref{lplusxy}), 
$\gamsw(\omega) \subseteq \gamswcl(\omega)$ and
 the connectedness of $\gamswcl(\omega)$,
we find that
 $\gamswcl(\omega) \cap A_{\bo{x},\bo{y}}^c$ is connected. Using this alongside $\omega'\big\vert_{E(\Z^2) \setminus \axye} = \omega\big\vert_{E(\Z^2) \setminus \axye}$, we find that the $\omega'$-open cluster of $\csw$ in $\axy^c$ equals 
 $\gamswcl(\omega) \cap A_{\bo{x},\bo{y}}^c$. The $\omega'$-open cluster 
 $\gamswcl(\omega')$ of $\csw$ therefore contains 
 $\gamswcl(\omega) \cap A_{\bo{x},\bo{y}}^c$; it also contains 
 $\overline\gamma_{\bo{x},\bo{y}}(\omega')$, by (\ref{lplusxy}). However, since 
 $\overline\gamma_{\bo{x},\bo{y}}(\omega') \subseteq \cluh{\bo{x}} \cap \clumh{\bo{y}}$ implies that 
$\overline\gamma_{\bo{x},\bo{y}}(\omega') \cap \partial A_{\bo{x},\bo{y}} = \big\{ \bo{x},\bo{y}\big\}$, the set $\gamswcl(\omega')$ has no further elements. Thus, 
(\ref{gsomone}).

To establish (\ref{gsomtwo}), 
we will firstly confirm that
\begin{equation}\label{eqtwoparts}
  \gamsw \big( \omega' \big) \cap \partial A_{\bo{x},\bo{y}} = 
 \big\{ \bo{x},\bo{y} \big\}, \qquad \textrm{and that $\gamsw \big( \omega' \big)$ intersects the interior of $A_{\bo{x},\bo{y}}$.}
\end{equation}
We know that
 $\gamsw(\omega) \cap A_{\bo{x},\bo{y}}^c$ is an open path from $\bo{x}$ to $\bo{y}$ under $\omega'\big\vert_{E(\Z^2) \setminus \axye}$, since $\omega'\big\vert_{E(\Z^2) \setminus \axye} = \omega\big\vert_{E(\Z^2) \setminus \axye}$.
This implies that the right-hand-side of (\ref{gsomtwo}) is an $\omega'$-open circuit.
It is southwest-centred at $\csw$: indeed,
$$
 \csw = {\rm SW} \big( \gamsw(\omega) \big) = 
{\rm SW}  \Big(  \big( \gamsw(\omega)   \cap \axy^c \big) \cup \gamma_{\bo{x},\bo{y}}(\omega') \Big),
$$
by the definition of $\gamsw(\omega)$ in the first equality, and each point in $\axy \subseteq \big\{ \bo{z} \in \R^2: z_1 \geq 0 \big\}$ being lexicographically greater than $\csw$, as well as 
$\gamma_{\bo{x},\bo{y}} \big( \omega' \big)  \subseteq \axy$.

The right-hand-side of (\ref{gsomtwo}) traps $\bo{0}$ in its interior. Therefore, so must $\gamsw ( \omega' )$, since ${\rm INT}\big( \gamsw ( \omega' ) \big)$ contains the interior of any $\omega'$-open circuit that is southwest-centred at $\csw$. 

That $\bo{0} \in {\rm INT}\big( \gamsw ( \omega' ) \big)$ forces
$\gamsw ( \omega' )$ to intersect the interior of $A_{\bo{x},\bo{y}}$. From (\ref{gsomone}), 
we know that $\ell_{\bo{0},\bo{z}}^+ \cap \gamswcl (\omega') = \{ \bo{z} \}$ for $\bo{z} = \bo{x},\bo{y}$.
By $\gamsw(\omega') \subseteq \gamswcl(\omega')$ and $\bo{0} \in {\rm INT} \big(  \gamsw(\omega') \big)$,
we infer that $\gamsw(\omega') \cap \partial \axy = \big\{ \bo{x}, \bo{y} \big\}$. We have derived (\ref{eqtwoparts}). 

By (\ref{eqtwoparts}), and each point in $\axy$ being lexicographically at least $\csw$,
we see that $\gamsw(\omega') \cap \axy$ is the outermost $\omega'$-open path from $\bo{x}$ to $\bo{y}$
in $\axy$, i.e.
\begin{equation}\label{gamtxy} 
   \gamsw(\omega') \cap \axy = \gamma_{\bo{x},\bo{y}}(\omega').
\end{equation}
Let $\zeta_\omega$ denote the outermost $\omega$-open path  $\zeta$ from $\bo{x}$ to $\bo{y}$
in $\axy^c$ for which ${\rm SW} \big( \zeta \big) = \csw$, and write $\zeta_{\omega'}$ analogously. Then note that, by  (\ref{eqtwoparts}), $\gamsw(\omega') \cap \axy^c = \zeta_{\omega'}$.  
By
 $\omega'\big\vert_{E(\Z^2) \setminus \axye} = \omega\big\vert_{E(\Z^2) \setminus \axye}$,
we see that $\zeta_{\omega'} = \zeta_{\omega}$. However, $\zeta_{\omega} = \gamsw(\omega) \cap \axy^c$, since, as we have seen, 
$\gamsw(\omega) \cap \axy^c$ is a path from $\bo{x}$ to $\bo{y}$. That is, 
$\gamsw(\omega') \cap \axy^c = \gamsw(\omega) \cap \axy^c$, which, alongside 
(\ref{gamtxy}), gives (\ref{gsomtwo}). 

We record the formulae
\begin{equation}\label{gstwo}
 \gamsw \big( \omega \big)
 = \Big(  \gamsw \big( \omega \big)\cap \axy^c  \Big) \cup \gamma_{\bo{x},\bo{y}}(\omega),
\end{equation}
and
\begin{equation}\label{gsone}
 \gamswcl \big( \omega \big)
 = \Big(  \gamswcl \big( \omega \big) \cap \axy^c \Big) \cup 
   \overline\gamma_{\bo{x},\bo{y}}(\omega).
\end{equation}
that are analogous to (\ref{gsomtwo}) and (\ref{gsomone}).

To verify that $P \big( \cdot \big\vert B_{\bo{x},\bo{y}} \big)$
is invariant under $\psi_{\bo{x},\bo{y}}$, let $\omega_0 \in \{0,1\}^{E(\Z^2) \setminus \axye}$ be of the form $\omega_0 = \omega \big\vert_{E(\Z^2) \setminus \axye}$ for some $\omega \in  B_{\bo{x},\bo{y}}$. It suffices to check that
\begin{equation}\label{eqpphiinv}
  P \Big( \cdot \Big\vert  B_{\bo{x},\bo{y}},  
\omega\big\vert_{E(\Z^2) \setminus \axye} = \omega_0  \Big) = 
  \Big( P \circ \psi_{\bo{x},\bo{y}}^{-1} \Big) \Big( \cdot \Big\vert  B_{\bo{x},\bo{y}},  
\omega\big\vert_{E(\Z^2) \setminus \axye} = \omega_0  \Big).
\end{equation}
Under the measure on the left-hand-side of (\ref{eqpphiinv}),
$\omega\big\vert_{\axye}$ has the conditional distribution of the marginal of $P$ in $\axye$ given 
$\omega\big\vert_{E(\Z^2) \setminus \axye} = \omega_0$
and
$\big\{ \bo{x}  \build\leftrightarrow_{}^{\axy} \bo{y} \big\}
 \cap \big\{  \overline\gamma_{\bo{x},\bo{y}} \subseteq \cluh{\bo{x}} \cap \clumh{\bo{y}} \big\} \cap \big\{ \big\vert {\rm INT} \big( (\gamsw(\omega_0) \cap \axy^c) \cup  \gamma_{\bo{x},\bo{y}}(\omega) \big) \big\vert \geq n^2  \big\}$. This is true for the following two reasons. Firstly,  given 
$\omega\big\vert_{E(\Z^2) \setminus \axye} = \omega_0$, the event 
$\big\{ \bo{x}  \build\leftrightarrow_{}^{\axy} \bo{y} \big\}
 \cap \big\{  \overline\gamma_{\bo{x},\bo{y}} \subseteq \cluh{\bo{x}} \cap \clumh{\bo{y}} \big\}$
is characterized by $\big\{ \bo{x},\bo{y} \in {\rm RG} \big( \gamswcl \big) \big\} \cap \big\{ \bo{0} \in  {\rm INT} ( \gamsw ) \big\}$, (this by virtue of the definition of 
 ${\rm RG} \big( \gamswcl \big)$, (\ref{gsone}), $\argu(\bo{y}) < \argu(\bo{x}) + c_0/2$ and
 $\cluh{\bo{x}} \cap \clumh{\bo{y}} \subseteq \axy$). Secondly, by the use of (\ref{gstwo}). Under the measure on the right-hand-side, 
$\omega\big\vert_{\axye}$ has exactly the same conditional distribution: this is proved by use of  (\ref{gsomtwo}) and (\ref{gsomone}) alongside the definition of $\psi_{\bo{x},\bo{y}}$ (which was designed to ensure this property). \qed
\end{subsection}
\begin{subsection}{Proofs related to recentering the circuit}
We present here the proofs of Lemmas \ref{lemlas}, \ref{lemoswtrans} and \ref{lemmacat}.
We start with a simple assertion.  
\begin{lemma}\label{lematcn}
There exist $0 < c , C < \infty$ such that
$$
P \Big( \gamsw \not\subseteq B_{C n} \Big\vert \areatrap \Big)
 \leq \exp \big\{ - c n \big\}. 
$$
\end{lemma}
\noindent{\bf Proof.}
We begin by arguing that there exists $c > 0$ for which $P(\areatrap) \geq \exp \big\{ - cn \big\}$.
Let $\Gamma$ denote any circuit for which $\big\vert {\rm INT}(\Gamma) \big\vert \geq n^2$ and 
${\rm SW}(\Gamma) = \csw$. By definition, we have that $P(\areatrap)$ is at least the probability that $\Gamma$ is open. Such a $\Gamma$ may be chosen so that $\big\vert E(\Gamma) \big\vert$ has order $n$. Hence, the bounded energy property of $P$ implies that indeed  $P(\areatrap) \geq \exp \big\{ - cn \big\}$ for some $c > 0$.

It suffices then to argue that
\begin{equation}\label{gamswupbd}
 P \Big( \gamsw \not\subseteq B_{Cn} \Big) \leq \exp \big\{ - C' n \big\}
\end{equation}
for all $n$ sufficiently high, where $C' > 0$ may be chosen so that $C' \to \infty$ as $C \to \infty$.
It readily follows from Lemma \ref{lemmac} that $\dist \csw \dist \leq \cctwo n$. As such, for $C > \cctwo$, the event $\big\{ \gamsw \not\subseteq B_{Cn} \big\}$
entails the existence of an open path of length at least $(C - \cctwo)n$ that contains an element of $B_{\cctwo n}$. Setting $C = 2 \cctwo$, we see that $P \big( \gamsw \not\subseteq B_{2 \cctwo n}\big) \leq 2\pi \big( \cctwo n \big)^2 \exp \big\{ - c \cctwo n \big\}$, by use of the exponential decay of connectivity property satisfied by $P$ (with the constant $c > 0$ being from the statement of this property in Section \ref{sechyp}). In this way, we obtain (\ref{gamswupbd}). \qed
\noindent{\bf Proof of Lemma \ref{lemlas}.} Let $\bo{v} \in \Z^2$ be such that 
$$
p_{\bo{v}} : = P \Big( \sw{\cir} = \csw + \bo{v} \Big\vert \areacon, \cir \in \mathcal{M} \Big)
$$
is maximal. It is a consequence of Lemma \ref{lemmac} that either $p_{\bo{v}} \geq \frac{1}{2 \pi \cctwo^2 n^2}$ or 
\begin{equation}\label{eqwwr}
P \big( \cir \in \mathcal{M} \big\vert \areacon  \big) \leq \exp \big\{ - c n \big\}, 
\end{equation}
for some small constant $c > 0$. Indeed, suppose that the latter alternative fails for  $c = \epsilon/2$, 
where $\epsilon > 0$ is specified in Lemma \ref{lemmac}. Then we learn from Lemma \ref{lemmac} that
$$
P \Big( \cir \not\subseteq B_{\cctwo n} \Big\vert \areacon, \cir \in \mathcal{M} \Big) \leq \exp \big\{ - (\epsilon/2) n \big\}.
$$  
Under the law $P\big( \cdot \big\vert  \areacon, \cir \in \mathcal{M} \big)$, 
the random vertex $\sw{\cir}$ is thus likely to lie in $B_{\cctwo n}$, whence $p_{\bo{v}}  \geq \frac{1}{2 \pi \cctwo^2 n^2}$.

The form of the statement of the lemma permits us to exclude the case that (\ref{eqwwr}) holds. For any $\omega \in \zoz$
and $\bo{x} \in \Z^2$, we write $\omega_{\bo{x}} : = \omega \big( \cdot + \bo{x}  \big)$
for the translation of $\omega$ by $- \bo{x}$. Note that if 
$\omega \in  \big\{ \sw{\cir} = \csw + \bo{v} \big\} \cap \big\{ \areacon \big\} \cap \big\{ \cir \in \mathcal{M} \big\}$, then
$\omega_{\bo{v}} \in \areatrap \cap \big\{ \gamsw \in \mathcal{M} \big\}$, because, as is readily verified, 
$\gamsw(\omega_{\bo{v}}) = \cir(\omega) - \bo{v}$. From this, and the translation-invariance of $P$,
we learn that
$$
P \Big( \areatrap \cap \big\{ \gamsw \in \mathcal{M} \big\} \Big)
 \geq 
P \Big( \big\{ \sw{\cir} = \csw + \bo{v} \big\}  \cap \big\{ \areacon \big\} \cap \big\{ \cir \in \mathcal{M} \big\} \Big).
$$
By the assumed bound on $p_{\bo{v}}$,
\begin{equation}\label{areatrapineq}
P \Big( \areatrap \cap \big\{ \gamsw \in \mathcal{M} \big\} \Big)
 \geq \frac{1}{2\pi \cctwo^2} n^{-2} 
P \Big( \Big\{ \areacon \Big\} \cap \Big\{ \cir \in \mathcal{M} \Big\} \Big).
\end{equation}
We divide (\ref{areatrapineq}) by $P \big( \areatrap \big)$ and use 
\begin{equation}\label{atubd}
 P \big( \areatrap \big) \leq 10 \pi \cctwo^2 n^2 P \big( \areacon \big)
\end{equation}
to obtain 
$$
10 \pi \cctwo^2 n^2 P \Big( \gamsw \in \mathcal{M} \Big\vert \areatrap \Big)
 \geq \frac{1}{2\pi \cctwo^2} n^{-2}  
  P \Big( \cir \in \mathcal{M} \Big\vert \areacon \Big),
$$
which yields the statement of the lemma in the case that remains (when (\ref{eqwwr}) fails).

It remains to verify (\ref{atubd}). By Lemma \ref{lematcn}, we have that
$$
P \big( \areatrap \big) \leq 2 P \big( \areatrap, \gamsw \subseteq B_{\cctwo n} \big).
$$
The event on the right-hand-side entails that the ball $B_{\cctwo n}$ contain a circuit that traps an area of at least $n^2$. For any such configuration, some translate, by an integer vector of norm at most $\cctwo n$, will realize   $\areacon$. 
Thus, the bound (\ref{atubd}) follows from the translation invariance of $P$.
\qed
\noindent{\bf Proof of Lemma \ref{lemoswtrans}.} The proof is very similar to that of Lemma \ref{lemlas}. In place of (\ref{atubd}),
we use
\begin{equation}\label{ineqatac}
P \big( \areatrap \big) \geq \frac{1}{2\pi \cctwo^2 n^2} P \Big( \areacon \Big),
\end{equation}
which follows from 
$$
P \Big( \areatrap \Big\vert \areacon  \Big) \geq 
P \Big( \sw{\cir} = \csw \Big\vert \areacon  \Big) \geq 
\frac{1}{2\pi \cctwo^2 n^2},
$$
where the second inequality is due to the Definition \ref{defcsw} of $\csw$. \qed
The next two lemmas give Lemma \ref{lemmacat}.
\begin{lemma}\label{lemgtsat}
There exists $c > 0$ such that, for any $\epsilon \in (0,c)$, and for all $n \in \N$ sufficiently high,
$$
P \Big(  {\rm GD}(\gamsw) \geq \epsilon n \Big\vert \areatrap \Big) \leq \exp \big\{ - c \epsilon n \big\}.
$$
Under this conditional measure, $\centre(\gamsw) \in {\rm INT}\big( \gamsw \big)$ except with exponentially decaying probability in $n$.
\end{lemma}
\noindent{\bf Proof.}
The result follows from Proposition \ref{propglobdis} and Lemma \ref{lemoswtrans}. \qed
%
\begin{lemma}\label{lemab}
For all $\cona > 0$, there exists $\conb > 0$ such that, for $n$ sufficiently high,
$$
P \Big( \centre(\gamsw) \not\in B_{\cona n}  \Big\vert \areatrap \Big) \leq \exp \big\{ - \conb n \big\}.
$$
\end{lemma}
\noindent{\bf Proof.} 
We parametrize the Wulff curve by the polar angle, writing 
$\partial \wulff (\theta)$ for the unique element of $\partial \wulff$ of argument
$\theta$ (for $\theta \in [0,2\pi)$). We claim that, for all $\delta > 0$,
there exists $\epsilon > 0$ such that, for $n$ sufficiently high, 
and for any circuit $\Gamma$ satisfying  $\centre(\Gamma) \in {\rm INT}(\Gamma)$ and  
${\rm GD}(\Gamma) \leq \epsilon n$, we have that
\begin{equation}\label{eqswcn}
\vert\vert   {\rm SW} \big( \Gamma \big) - \centre\big(\Gamma\big)  -n \partial \wulff \big( \pi \big) \vert\vert \leq  \delta n.
\end{equation} 
Indeed, 
there exists a constant $\cstar > 0$ such that
\begin{equation}\label{gamcen}
\Gamma - \centre\big( \Gamma \big) \subseteq \big( n + \cstar \epsilon n \big) \wulff \setminus 
 \big( n - \cstar \epsilon n \big) \wulff.
\end{equation}
(In essence, this follows from the definition of ${\rm GD}(\Gamma)$ and the convexity of $\wulff$. A detailed proof of (\ref{gamcen}) is given in Lemma 3.3 of \cite{hammondthr}.)
From (\ref{gamcen}), we see that 
$\Gamma - \centre\big( \Gamma \big)$, being a circuit whose interior contains $\bo{0}$, contains a point of $x$-coordinate at most $- \big( n - \cstar \epsilon n \big) \dist \partial \wulff (\pi)\dist$
but no point of $x$-coordinate at most $- \big( n + \cstar \epsilon n \big) \dist \partial \wulff (\pi)\dist$. Hence,
$$
- \big( n + \cstar \epsilon n \big) \dist \partial \wulff (\pi)\dist
  \leq {\rm SW}\big( \Gamma  \big)_1 - \centre\big(\Gamma \big)_1 \leq
 - \big( n - \cstar \epsilon n \big) \dist \partial \wulff (\pi)\dist
$$
By the strict convexity of $\partial \wulff$ (Lemma \ref{lemozstr}), for all $\delta > 0$, 
there exists $\epsilon > 0$, such that
\begin{eqnarray}
 & & \Big( \big( n + \cstar \epsilon n \big) \wulff \setminus 
 \big( n - \cstar \epsilon n \big) \wulff \Big) \cap \Big\{ (x,y) \in \R^2: x \leq - \big(  n - \cstar \epsilon n  \big)  \dist \partial \wulff (\pi)\dist  \Big\} \nonumber \\
 & \subseteq & \big\{ (x,y) \in \R^2: \big\vert y \big\vert < \delta n/2 \big\}. \nonumber
\end{eqnarray}
Hence,  
$\vert\vert   {\rm SW} \big( \Gamma \big) - \centre\big(\Gamma\big)  -n \partial \wulff \big( \pi \big) \vert\vert \leq  \sqrt{(2 \cstar \eps)^2 + (\delta/2)^2} n$, which implies (\ref{eqswcn}), if we impose that $\epsilon < \frac{\delta}{4\cstar}$.

Note that Lemma \ref{lemgtsat}, (\ref{eqswcn}) and the Definition \ref{defcsw} of $\csw$ imply that 
$$
\vert\vert   \csw   - n \partial \wulff \big( \pi \big) \vert\vert \leq \delta n.
$$
Note that, by (\ref{eqswcn}), any circuit $\Gamma$ for which 
$\centre(\Gamma) \not\in B_{\cona n}$, $\centre(\Gamma) \in {\rm INT}(\Gamma)$ and ${\rm GD}(\Gamma) \leq \epsilon n$ satisfies  
$\vert\vert {\rm SW} \big( \Gamma \big) -  n \partial W \big( \pi \big) \vert\vert \geq \cona n  - \delta n$. 
Hence, if $\cona > 2 \delta$, the events 
$\big\{ {\rm GD}(\gamsw) \leq \epsilon n \big\} \cap \big\{ \centre(\gamsw) \in {\rm INT}(\gamsw) \big\}$ and  $\big\{ \centre(\gamsw) \not\in B_{\cona n} \big\}$ 
are disjoint. For this reason, Lemma \ref{lemgtsat} implies the result. \qed
Note that Lemma \ref{lemmacat} follows from Lemmas \ref{lematcn} and \ref{lemab}, with a relabelling of constants.
\end{subsection}
\begin{subsection}{Regeneration structure and recentering}
We now explain how to derive Theorem \ref{thmmaxrgclgw} from Theorem \ref{thmmaxrgcl}.
The centre $\centre(\gamsw)$ of the circuit $\gamsw$ arising from the measure $P \big( \cdot \big\vert \areatrap \big)$ is not necessarily $\bo{0}$, but it is nearby to $\bo{0}$. We establish this as a preliminary to the proof of Theorem \ref{thmmaxrgclgw}:
\begin{lemma}\label{lemrgtran}
Let $c_2 > 0$ satisfy $\frac{3 \pi c_2}{\ccone} < \min \big\{ \qzero, c_0 \big\}$.  
For $\bo{x} \in \Z^2$ and $q,c > 0$, we use the notation ${\rm RG}_{q,c,\bo{x}}\big(\Gamma  \big)$ to denote $\bo{x} + {\rm RG}_{q,c}\big(\Gamma - \bo{x} \big)$. Let $\Gamma$ be a circuit such that $\Gamma \cap B_{\ccone n} = \emptyset$. Let $\bo{x} \in \Z^2$ satisfy $\vert\vert \bo{x} \vert\vert \leq c_2 n$.
Then
$$
 {\rm RG}_{\qzero,c_0}\big(\Gamma  \big)
\subseteq 
{\rm RG}_{\qzero/2,c_0/2,\bo{x}}\big(\Gamma  \big).
$$
\end{lemma}  
\noindent{\bf Proof.} We begin with a useful claim.

Let $\bo{v},\bo{x} \in \Z^2$ satisfy  $\vert\vert \bo{x} \vert \vert \leq c_2 n$ and
 $\vert\vert \bo{v} \vert \vert \geq \ccone n$. 
We claim that 
\begin{equation}\label{xwbx}
  \Big( \bo{x} + W_{\bo{v} - \bo{x},c_0/2} \Big) \cap B_{\ccone n}^c \subseteq W_{\bo{v},c_0}.
\end{equation}
Let $\ell_1$ and $\ell_2$ denote the clockwise and counterclockwise boundaries of $W_{\bo{x},c_0}$. 
Let $\ell_1^*$ and $\ell_2^*$ denote the corresponding boundaries of $\bo{x} + W_{\bo{v} - \bo{x},c_0/2}$. Write $\bo{w_1},\bo{w_2},\bo{w_1^*},\bo{w_2^*}$, for the intersections of  $\ell_1, \ell_2,\ell_1^*,\ell_2^*$ with $\partial B_{\ccone n}$. For (\ref{xwbx}), it suffices to show that
\begin{equation}\label{xwbxo}
 \arg \big( \bo{w_1} \big) < \arg \big( \bo{w_1^*} \big) <\arg \big( \bo{w_2^*} \big) <\arg \big( \bo{w_2} \big) 
\end{equation}
and that
\begin{equation}\label{xwbxt}
 \big( \ell_1^* \cup \ell_2^* \big) \cap \big( \ell_1 \cup \ell_2 \big) \cap 
 B_{\ccone n}^c = \emptyset.
\end{equation}  
\begin{figure}\label{figregreloc}
\begin{center}
\includegraphics[width=0.5\textwidth]{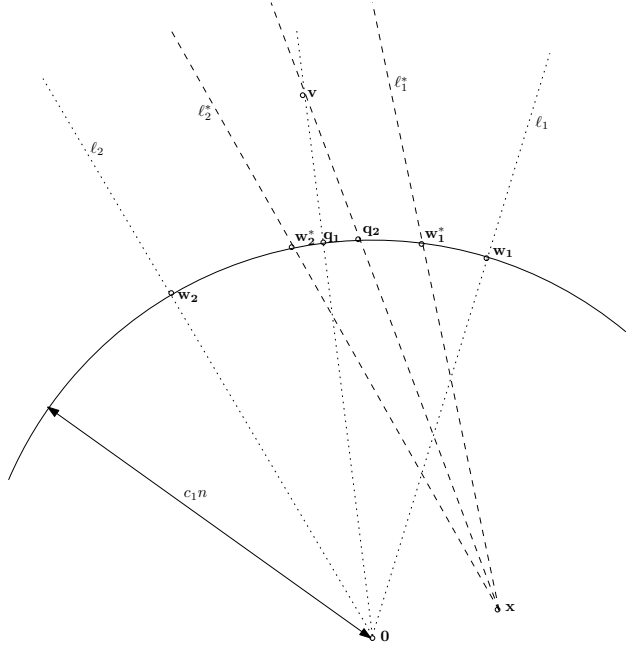} \\
\end{center}
\caption{Illustrating the derivation of (\ref{xwbx}).}
\end{figure}
Set $\bo{q_1}$ equal to the point of intersection 
$\ell_{\bo{0},\bo{v}}^+ \cap \partial B_{\ccone n}$
and $\bo{q_2}$ equal to the point 
$\ell_{\bo{x},\bo{v}}^+ \cap \partial B_{\ccone n}$.
To show (\ref{xwbxo}), we begin by bounding $\ang\big( \bo{q_1} , \bo{q_2} \big)$.
By considering the triangle $T_{\bo{0},\bo{v},\bo{q_2}}$, we have that
$\ang\big( \bo{q_1} , \bo{q_2} \big) + \ang\big( \bo{v} , \bo{v} - \bo{x} \big)    =  \ang\big( \bo{q_2} , \bo{q_2} - \bo{x} \big)$. Hence,
$\ang\big( \bo{q_1} , \bo{q_2} \big) \leq \ang\big( \bo{q_2} , \bo{q_2} - \bo{x} \big)$. However, 
$\sin \ang\big( \bo{q_2} , \bo{q_2} - \bo{x} \big) \leq \frac{\dist \bo{x} \dist}{\dist \bo{q_2}\dist} \leq \frac{c_2}{\ccone}$, so that  
$\ang\big( \bo{q_2} , \bo{q_2} - \bo{x} \big) \leq \frac{\pi c_2}{2 \ccone}$. Thus, 
$\ang\big( \bo{q_1} , \bo{q_2} \big) \leq   \frac{\pi c_2}{2 \ccone}$.

Next, we bound
$\ang\big( \bo{q_2}, \bo{w_1^*} \big)$. We have that
\begin{equation}\label{angqineq}
\ang\big( \bo{q_2}, \bo{w_1^*} \big) \leq 
\ang\big( \bo{q_2}, \bo{q_2} - \bo{x} \big) +  \ang\big( \bo{q_2} - \bo{x}, \bo{w_1^*} - \bo{x} \big)
 + \ang\big(  \bo{w_1^*} - \bo{x},\bo{w_1^*} \big).
\end{equation}
Note that, for any $\bo{u}$ satisfying $\dist \bo{u} \dist \geq \ccone n$,
\begin{equation}\label{anguxineq}
\ang \big( \bo{u} , \bo{u} - \bo{x} \big) \leq \frac{\pi c_2}{2 \ccone}.
\end{equation}
Since $\dist \bo{q_2} \dist, \dist \bo{w_1^*} \dist = \ccone n$, we may use this, alongside 
$\ang\big( \bo{q_2} - \bo{x}, \bo{w_1^*} - \bo{x} \big) = c_0/2$, to obtain from (\ref{angqineq})
that
$\ang\big( \bo{q_2}, \bo{w_1^*} \big) \leq c_0/2 + \pi c_2/c_1$.
We also have that 
$\argu(\bo{w_1}) = \argu(\bo{q_1}) - c_0$. Hence, 
\begin{eqnarray}
 & & \argu(\bo{w_1^*}) \geq  \argu(\bo{q_2}) - c_0/2 - \frac{\pi c_2}{2 \ccone}  \geq \argu(\bo{q_1}) - 
\ang \big( \bo{q_1},\bo{q_2}  \big) - c_0/2 - \frac{\pi c_2}{2 \ccone}  \nonumber \\ 
& = & 
 \argu(\bo{w_1}) + c_0 - \frac{3\pi c_2}{2 \ccone} - c_0/2 > \argu(\bo{w_1}), \nonumber
\end{eqnarray}
since  $3\pi c_2 < c_0 c_1$. We similarly have that
$\argu(\bo{w_2}) > \argu(\bo{w_2^*})$. Note that $\dist \bo{x} \dist < c_1 n$
implies that $\argu(\bo{w_1^*}) < \argu(\bo{w_2^*})$. Hence, we have (\ref{xwbxo}).

We may now derive the statement of the lemma.
Let $\bo{v_1}$,  $\bo{v_2}$, $\bo{v_1^*}$, $\bo{v_2^*}$ denote the direction vectors of $\ell_1$ and  $\ell_2$ (away from $\bo{0}$) and $\ell_1^*$ and $\ell_2^*$
(away from $\bo{x}$). Given (\ref{xwbxo}), it suffices for (\ref{xwbxt}) that
\begin{equation}\label{xwbxth}
 \arg \big( \bo{v_1} \big) \leq \arg \big( \bo{v_1}^* \big) \leq \arg \big( \bo{v_2}^* \big) \leq \arg \big( \bo{v_2} \big). 
\end{equation}
Note that $\arg \big( \bo{v_1} \big) = \arg \big( \bo{v} \big) - c_0$. Note further that
$\big\vert \arg \big( \bo{v_1^*} \big) - \arg \big( \bo{v} \big) \big\vert \leq 
\big\vert \arg \big( \bo{v_1^*} \big) - \arg \big( \bo{v} - \bo{x} \big) \big\vert
 + \ang \big( \bo{v}, \bo{v} - \bo{x} \big)$.
However, 
$\big\vert \arg \big( \bo{v_1^*} \big) - \arg \big( \bo{v} - \bo{x} \big) \big\vert = c_0/2$,
while 
$\ang \big( \bo{v} , \bo{v} - \bo{x} \big) \leq \frac{\pi c_2}{2 \ccone}$
by
(\ref{anguxineq}) and $\dist \bo{v} \dist \geq c_1 n$.
Hence, 
$$
\argu \big( \bo{v_1^*} \big) \geq \argu \big( \bo{v} \big) - c_0/2 - \frac{\pi c_2}{2 \ccone}
 \geq \argu \big( \bo{v_1} \big) + c_0/2  - \frac{\pi c_2}{2 \ccone}.
$$
Hence, $\argu \big( \bo{v_1^*} \big) \geq \argu \big( \bo{v_1} \big)$ follows from 
$\pi c_2 \leq c_0 \ccone$. The other inequalities being derived either identically or trivially, we obtain 
(\ref{xwbxth}).


Let $\bo{v} \in  {\rm RG}_{\qzero,c_0}\big(\Gamma  \big)$. We must show that
\begin{equation}\label{gaminc}
\Gamma \cap \Big( \bo{x} +  W_{\bo{v} - \bo{x},c_0/2} \Big)
\subseteq 
 \Big( \bo{x} + C_{\pi/2 - \qzero/2}^F \big( \bo{v} - \bo{x} \big)
 \Big) \cup
 \Big( \bo{x} + C_{\pi/2 - \qzero/2}^B \big( \bo{v} - \bo{x} \big)
 \Big).
\end{equation} 
Let $\bo{y} \in 
\Gamma \cap \big( \bo{x} +  W_{\bo{v} - \bo{x},c_0/2} \big)$.
By (\ref{xwbx}) and $\Gamma \cap B_{\ccone n}^c = \emptyset$, we find that 
$\bo{y} \in \Gamma \cap W_{\bo{v},c_0}$. From $\bo{v} \in {\rm RG}_{q_0,c_0}\big( \Gamma \big)$,
we obtain $\bo{y} \in \clu{\bo{v}} \cup \clum{\bo{v}}$. 

For (\ref{gaminc}), it suffices then to show that
$\clu{\bo{v}} \subseteq \bo{x} + C_{\pi/2 - \qzero/2}^F \big( \bo{v} - \bo{x} \big)$
and  
$\clum{\bo{v}} \subseteq \bo{x} + C_{\pi/2 - \qzero/2}^B \big( \bo{v} - \bo{x} \big)$.
These two statements are implied by $\ang\big( \bo{v} , \bo{v} - \bo{x} \big) \leq \qzero/2$. This follows from (\ref{anguxineq}), $\dist \bo{v} \dist \geq c_1 n$ and  $\pi c_2 < \qzero \ccone$. \qed
\noindent{\bf Proof of Theorem \ref{thmmaxrgclgw}.}
By Lemma \ref{lemab} and the form of (\ref{eqmarcl}), we may suppose that there exists 
$\bo{v} \in B_{\cona n}$ such that
\begin{eqnarray}
 & & P \Big(  \marcl > u/n, \centre\big(\gamsw\big) = \bo{v} 
  \Big\vert \areatrap \Big) \nonumber \\
 & \geq &  \frac{1}{2\pi \cona^2 n^2}
  P \Big(  \marcl > u/n
  \Big\vert \areatrap \Big),\label{eqthtmax}
\end{eqnarray}
where $\cona > 0$ is an arbitrary constant.
For $\bo{u} \in \Z^2$, recall that $\Gamma_{\bo{u}}$ denotes the outermost open circuit whose interior contains $\bo{u}$. We write $\overline{\Gamma}_{\bo{u}}$ for the open cluster in which $\Gamma_{\bo{u}}$ is contained, analogously to the notation $\gamswcl$.
We claim that, given 
\begin{equation}\label{eqposp}
 \Big\{  \marcl > u/n, \centre\big(\gamsw\big) = \bo{v}, \areatrap \Big\},
\end{equation}
the conditional probability that $\overline{\Gamma}_{\bo{v}}$ is equal to $\gamswcl$ is bounded away from zero, uniformly in $n$. Indeed, conditionally on (\ref{eqposp}), $\overline{\Gamma}_{\bo{v}} = \gamswcl$ provided that there is no open cluster disjoint from $\gamswcl$ that encircles $\gamswcl$. 
To see this, for $\bo{x} \in \Z^d$, let $C(\bo{x})$ denote the event that there exists an infinite closed path emanating from $\bo{x}$ whose lexicographically maximal element is $\bo{x}$. Let $H^+_{\bo{x}} \subseteq \Z^2$ denote the set of vertices whose lexicographical order is at least that of $\bo{x}$.  We assert that there exists $c > 0$ such that, for each $\bo{x} \in \Z^2$ and for all $\omega' \in \{ 0,1\}^{E(H^+_{\bo{x}})}$,
\begin{equation}\label{eqcxlbd}
  P \Big( C\big(\bo{x}\big) \Big\vert \omega \big\vert_{E(H^+_{\bo{x}})} = \omega' \Big) \geq c.
\end{equation}
\begin{figure}\label{figoutpath}
\begin{center}
\includegraphics[width=0.3\textwidth]{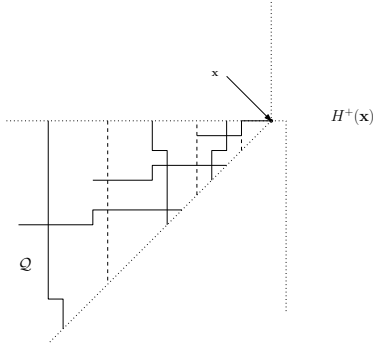} \\
\end{center}
\caption{A sketch that sketches the proof of (\ref{eqcxlbd}). To see that the unconditioned probability  $P \big( C(\bo{x}) \big)$ is positive, we divide the indicated wedge $\mathcal{Q}$ into a series of dyadic scales. We insist that, for each one of these, a vertical closed path crosses the scale, and a horizontal closed path crosses the scale and its neighbours. By exponential decay of connectivity, the absence of such paths has a probability that decays exponentially in the index of the scale. However, the presence of all such paths ensures a closed infinite path emanating from $\bo{x}$ that avoids $H^+_{\bo{x}}$. The probability of the intersection of the events for the different scales is at least the product of the probabilities by the FKG inequality, so that the intersection has positive probability, as required. Moreover, the same argument may be undertaken in the presence of  conditioning on data in $E(H^+_{\bo{x}})$, due to the ratio-weak-mixing property enjoyed by $P$ (that was recalled in Section \ref{sechyp}). This is because the edge-set regions $E(\mathcal{Q})$ and $E\big(H^+_{\bo{x}}\big)$ are ``well-separated'' in the spirit  of  Definition \ref{defrwm}, making what goes on inside one close to independent of what happens in the other.}
\end{figure}
The proof of (\ref{eqcxlbd}) is outlined in Figure 9 and the accompanying text.
Note that, conditionally on (\ref{eqposp}), the event $C\big( {\rm SW}(\gamswcl) - e_1 \big)$ ensures that there is no open cluster disjoint from $\gamswcl$ that encircles $\gamswcl$.  (Recall that $e_1$ denotes the unit vector in the positive $x$-direction.) As such, to prove the claim made involving (\ref{eqposp}), it suffices to establish that $C\big( {\rm SW}(\gamswcl) - e_1 \big)$ given (\ref{eqposp}) has a probability that is bounded below uniformly in the conditioning. This, however, follows from (\ref{eqcxlbd}) with the choice of $\bo{x} =  {\rm SW}(\gamswcl) - e_1$ being made, since the event (\ref{eqposp}) is measurable with respect to the configuration in
$E\big(H^+_{{\rm SW}(\gamswcl) - e_1}\big)$. This completes the proof of the claim involving (\ref{eqposp}).

We make use of the following result, whose proof appears after the end of the argument.
\begin{lemma}\label{lemrgtool}
For $\bo{v} \in \Z^2$ and $\omega \in \zoz$,
let $\omega_{\bo{v}} : = \omega \big( \bo{v} + \cdot \big)$
denote the configuration $\omega$ translated by $- \bo{v}$.

Let $\bo{v} \in \Z^2$ satisfy $\dist \bo{v} \dist \leq \cnew n$, with $\cnew = \min \big\{ \ccone  \sin(\qzero/2), \ccone \qzero/(2\pi), \ccone  c_0/(5\pi), \cctwo \}$. 
\begin{enumerate}
 \item 
Suppose that $\omega$ realizes the event
\begin{equation}\label{omevl}
 \Big\{ c_0/2 > \marcl > u/n, \centre\big(\gamswcl\big) = \bo{v}, \areatrap, \overline{\Gamma}_{\bo{v}} = \gamswcl, \gamswcl \subseteq B_{\cctwo n} \setminus B_{\ccone n} \Big\}.
\end{equation}
Then
$$
\omega_{\bo{v}} \in 
\Big\{ \marclo >  
\frac{\ccone \sin(\qzero/8) u}{\pi \cctwo n}, \areacon  \Big\}.
$$
\item
Suppose that $\omega$ realizes the event
$$
 \Big\{  \marcl \geq c_0/2, \centre\big(\gamswcl\big) = \bo{v}, \areatrap, \overline{\Gamma}_{\bo{v}} = \gamswcl, \gamswcl \cap  B_{\ccone n} = \emptyset \Big\}.
$$
Then
$$
\omega_{\bo{v}} \in 
\Big\{ \marclo >  
\frac{c_0}{4}, \areacon  \Big\}.
$$
\end{enumerate}
\end{lemma}
By the two parts of Lemma \ref{lemrgtool} and the translation invariance of $P$,
\begin{eqnarray} 
& & 
P \Big( \marcl > u/n, \centre\big(\gamswcl\big) = \bo{v}, \areatrap, \overline{\Gamma}_{\bo{v}} = \gamswcl, \gamswcl \subseteq B_{\cctwo n} \setminus B_{\ccone n} \Big) \nonumber \\
 & \leq & P \Big( \marclo > 
\frac{\ccone \sin(\qzero/8) u}{\pi \cctwo n} , \areacon  \Big), \nonumber
\end{eqnarray}
where we choose $\cona = \cnew$ so that $\bo{v} \in \Z^2$ is a permissible choice for the application of Lemma \ref{lemrgtool}.
Hence, 
\begin{eqnarray}
 & &  P \Big(   \marcl > u/n, \centre\big(\gamswcl\big) = \bo{v}
  \Big\vert \areatrap \Big) \nonumber \\
 & \leq &  2 \pi \cctwo^2 n^2 
C P \Big(   \marclo > \frac{\ccone \sin(\qzero/8) u}{\pi \cctwo n} 
  \Big\vert \areacon \Big) + \exp \big\{ - c n \big\}  \nonumber \\ 
 & \leq & C  n^2 \exp \big\{ - c u \big\},  \nonumber  
\end{eqnarray}
where the first inequality used (\ref{ineqatac}), the statement involving (\ref{eqposp}),
Lemma \ref{lemmacat}, and Theorem \ref{thmmaxrgcl}.
The second inequality used Theorem \ref{thmmaxrgcl} and $u \leq cn$. By (\ref{eqthtmax}),
we obtain
$$
P \Big(   \marcl > u/n \Big\vert \areatrap \Big) \leq 4 \pi^2 \cctwo^2 \cona^2 n^4 \exp \big\{ - cu \big\},
$$
as required. \qed
\noindent{\bf Proof of Lemma \ref{lemrgtool}.}
Note firstly that, by Lemma \ref{lemrgtran}
and $\bo{v} \in B_{\cnew n}$, $\pi \cnew/\ccone < \min \{ \qzero,c_0 \}$,
\begin{equation}\label{rginc}
 {\rm RG}_{\qzero,c_0,\bo{v}} \big( \gamswcl \big)
  \subseteq {\rm RG}_{\qzero/2,c_0/2} \big( \gamswcl \big).
\end{equation}
\noindent{\bf Proof of part (i).}
Let $\bo{r_1},\bo{r_2} \in \rggamswclex$ satisfy 
$\argu(\bo{r_1}) < \argu(\bo{r_2})$. We claim that
\begin{equation}\label{vrinc}
  \bo{v} \in B_{\ccone \sin ( \qzero/2) n} \, \, \, \textrm{implies that $\argu\big( \bo{r_1} - \bo{v} \big) < \argu\big( \bo{r_2} - \bo{v} \big)$.}
\end{equation}
It suffices to prove (\ref{vrinc}) for $\bo{r_1},\bo{r_2} \in \rggamswclex$
that are adjacent as viewed from $\bo{0}$, i.e., for which 
$A_{\bo{r_1},\bo{r_2}} \cap \rggamswclex = \big\{ \bo{r_1},\bo{r_2} \big\}$.
By $\marcl < c_0/2$, we have that 
$\argu(\bo{r_2}) - \argu(\bo{r_1}) < c_0/2$,
so that
$\bo{r_2} \in \rggamswclex$ implies that $\bo{r_1} \in C^B_{\pi/2 - \qzero/2}(\bo{r_2})$. This implies that 
$\argu\big( \bo{r_1} - \bo{w} \big) < \argu\big( \bo{r_2} - \bo{w} \big)$
for any $\bo{w}$ lying in the component of $\R^2 \setminus \big(  C^F_{\pi/2 - \qzero/2}(\bo{r_2}) \cup C^B_{\pi/2 - \qzero/2}(\bo{r_2}) \big)$ containing $\bo{0}$. Note that the point of closest approach to $\bo{0}$ in 
$\partial \big(  C^F_{\pi/2 - \qzero/2}(\bo{r_2}) \cup C^B_{\pi/2 - \qzero/2}(\bo{r_2}) \big)$
has distance from $\bo{0}$ at least  $\dist \bo{r_2} \dist \sin \big( \qzero/2 \big) \geq \ccone \sin \big( \qzero/2 \big) n$, the inequality due to $\bo{r_2} \in \gamswcl \subseteq B_{\ccone n}^c$.
We have verified (\ref{vrinc}).

By $\marcl > u/n$, we may find 
$\bo{r},\bo{r^*} \in \rggamswclex$, with
 $\argu(\bo{r^*}) > \argu(\bo{r}) + u/n$, 
such that 
$A_{\bo{r},\bo{r^*}} \cap \rggamswclex = \big\{ \bo{r},\bo{r^*} \big\}$.
The cyclic ordering of the elements of $\rggamswclex$ being the same whether these elements are viewed from $\bo{0}$
or from $\bo{v} \in B_{\ccone \sin(\qzero/2) n}$ (by (\ref{vrinc})), we find that
$$
 \Big(  \bo{v} + A_{\bo{r}-\bo{v},\bo{r^*}-\bo{v}}\Big) \cap {\rm RG}_{\qzero/2,c_0/2}\big( \gamswcl \big) = \big\{ \bo{r},\bo{r^*} \big\}.
$$
By (\ref{rginc}), then, 
$$
 \Big(  \bo{v} + A_{\bo{r}-\bo{v},\bo{r^*}-\bo{v}}\Big) \cap {\rm RG}_{\qzero,c_0,\bo{v}} \big( \gamswcl \big) \subseteq \big\{ \bo{r},\bo{r^*} \big\},
$$
implying that 
\begin{equation}\label{thtrineq}
  \theta^{\rm MAX}_{{\rm RG}_{\qzero,c_0}} \big( \gamswcl - \bo{v} \big)
 \geq \ang \big( \bo{r} - \bo{v}, \bo{r^*} - \bo{v} \big).
\end{equation}
We wish to bound from below the right-hand-side. 
To this end, note that
\begin{equation}\label{eqrstr}
  \dist \bo{r^*} - \bo{r} \dist \geq 2 \pi^{-1} \ccone n \ang \big( \bo{r},\bo{r^*} \big)
  \geq 2 \pi^{-1} \ccone u,
\end{equation}
the first inequality by $\bo{r},\bo{r^*} \in \gamswcl \subseteq B_{\ccone n}^c$
and the second by $\ang (\bo{r},\bo{r^*}) \geq u/n$.
We claim that 
\begin{equation}\label{rminvinc}
 \bo{r} - \bo{v} \in C^B_{\pi/2 - \qzero/4}\big( \bo{r^*} - \bo{v} \big).
\end{equation}
Indeed, by $\marcl < c_0/2$, we have that 
$\argu(\bo{r}) > \argu(\bo{r^*}) - c_0/2$, so that
$\bo{r^*} \in \rggamswclex$ yields 
\begin{equation}\label{eqrcb}
  \bo{r} \in C^B_{\pi/2 - \qzero/2} \big( \bo{r^*} \big).
\end{equation}
From (\ref{eqrcb}), it suffices for (\ref{rminvinc}) to show that
\begin{equation}\label{eqrrstq}
 \ang \big( \bo{r^*}, \bo{r^*} - \bo{v} \big) \leq \qzero/4.
\end{equation}
From $\bo{v} \in B_{\cnew n}$ and $\dist \bo{r^*} \dist \geq \ccone n$,
we have that $\sin \ang \big( \bo{r^*},\bo{r^*} - \bo{v} \big) \leq \cnew/\ccone$,
and thus
$\ang \big( \bo{r^*},\bo{r^*} - \bo{v} \big) \leq \pi\cnew/(2\ccone)$,
so that (\ref{eqrrstq}) follows from $\cnew \leq 2 \pi^{-1} \ccone \qzero/4$.
By (\ref{eqrstr}), (\ref{rminvinc}) and the upcoming Lemma \ref{lemdistang}, we obtain
$$
\ang \big( \bo{r^*} - \bo{v}, \bo{r} - \bo{v} \big)
  \geq \frac{2 \ccone  \sin(\qzero/8) u}{\pi \dist \bo{r^*} - \bo{v} \dist}.
$$
By $\bo{r^*} \in \gamswcl \subseteq B_{\cctwo n}$ and $\bo{v} \in B_{\cnew n} \subseteq B_{\cctwo n}$, we learn that
$$
\ang \big( \bo{r^*} - \bo{v}, \bo{r} - \bo{v} \big)
  \geq \frac{\ccone \sin(\qzero/8) u}{\pi \cctwo n}.
$$
From (\ref{thtrineq}), then, we obtain
$$
  \theta^{\rm MAX}_{{\rm RG}_{\qzero,c_0}} \big( \gamswcl - \bo{v} \big)
 \geq \frac{\ccone \sin(\qzero/8) u}{\pi \cctwo n}.
$$
This implies the first statement of the lemma. \\
\noindent{\bf Proof of part (ii).}
Firstly, note that 
\begin{equation}\label{eqrginfo}
  \bo{r} \in \gamswcl  \qquad \implies \qquad \ang \big( \bo{r},\bo{r} - \bo{v} \big) < \frac{c_0}{8}.
\end{equation}
Indeed,   $\bo{r} \in \gamswcl$ implies that $\dist \bo{r} \dist \geq \ccone n$, so that
$\sin \ang \big( \bo{r},\bo{r} - \bo{v} \big) \leq \dist \bo{v} \dist/\dist \bo{r} \dist \leq c_3/\ccone$, whence 
  $\ang \big( \bo{r},\bo{r} - \bo{v} \big) \leq  \frac{\pi c_3}{2\cctwo} < c_0/8$.

Let $\bo{r},\bo{r'} \in \rggamswclex$ satisfy
 $\argu(\bo{r'}) - \argu(\bo{r}) =  \theta^{\rm MAX}_{{\rm RG}_{\qzero/2,c_0/2}} \big( \gamswcl  \big)$,
with $\rggamswclex$ containing no element whose argument value lies between 
$\argu(\bo{r})$ and $\argu(\bo{r'})$.

Let $A^*$ denote the cone with apex at $\bo{v}$
whose clockwise boundary $\ell$ has argument 
$\argu(\bo{r} - \bo{v}) + c_0/8$, and whose counterclockwise boundary $\ell'$
 has argument 
$\argu(\bo{r'} - \bo{v}) - c_0/8$.

Note firstly that 
$$
\argu(\ell') - \argu(\ell) = \argu(\bo{r'} - \bo{v}) - \argu(\bo{r} - \bo{v}) - c_0/4
 \geq \argu(\bo{r'}) - \argu(\bo{r}) - c_0/2,
$$
since $\ang (\bo{s},\bo{s} - \bo{v}) \leq c_0/4$ for $\bo{s} = \bo{r},\bo{r'}$, 
by (\ref{eqrginfo}). Thus,
\begin{equation}\label{eqthtczero}
   \theta^{\rm MAX}_{{\rm RG}_{\qzero,c_0}} \big( \gamswcl - \bo{v} \big) > \frac{c_0}{4}
\end{equation}
will follow from the claim that
\begin{equation}\label{eqaoclaim}
  \textrm{the interior of $A^*$ is disjoint from ${\rm RG}_{\qzero,c_0,\bo{v}}\big( \gamswcl \big)$.}
\end{equation}
To see this, consider any $\bo{s} \in {\rm RG}_{\qzero,c_0,\bo{v}}\big( \gamswcl \big)$.
By (\ref{rginc}), we know that  
 $\bo{s} \in {\rm RG}_{\qzero/2,c_0/2}\big( \gamswcl \big)$,
so that
$\argu(\bo{s}) \not\in \big[ \argu(\bo{r}) , \argu(\bo{r'}) \big]$.
Suppose that $\argu(\bo{s}) \leq \argu(\bo{r})$.
Note then that
$\argu(\bo{s} - \bo{v}) \leq \argu(\bo{s}) + c_0/8 \leq \argu(\bo{r}) + c_0/8$,
the first inequality by (\ref{eqrginfo}).
That is, viewed from the vantage point of $\bo{v}$, $\bo{s}$ lies clockwise of $\ell$, the clockwise boundary of $A^*$.
Likewise, any such $\bo{s}$ for which $\argu(\bo{s}) \geq \argu(\bo{r'})$
lies counterclockwise to $\ell'$. This establishes (\ref{eqaoclaim}) and thus (\ref{eqthtczero}), 
and so concludes the proof of the second part of the lemma. \qed
We used the following straightforward result, Lemma 2.3 of \cite{hammondone}.
\begin{lemma}\label{lemdistang}
Let $q > 2c > 0$. 
If $\bo{x},\bo{y} \in \R^2$ satisfy $\ang \big( \bo{x}, \bo{y} \big) \leq c$   and $\bo{y} \in 
 C^F_{\pi/2 - q}\big( \bo{x}  \big) \cup C^B_{\pi/2 - q}\big( \bo{x}  \big)$, 
then $\vert\vert \bo{y} - \bo{x}  \vert\vert \leq \csc \big( q/2 \big) \vert\vert \bo{x} \vert\vert \ang \big( \bo{x} , \bo{y} \big)$.
\end{lemma}
\end{subsection}
\end{section}
\begin{section}{Concluding remarks}\label{secconcrem}
We conclude by briefly discussing some questions raised by our approach. 
\begin{subsection}{Rates of convergence to equilibrium}
We introduced a time-inhomogeneous Markov chain $\res$ that successively resamples sections of the conditioned circuit. The individual steps of the resampling, as described by the two-step formation appearing in the remark after Definition \ref{defrespsi}, 
act, in effect, by forming a subcritical  point-to-point connection between the endpoints of the circuit section where a modification is proposed, and then testing the result to see if the new circuit traps enough area. As we discussed in Section \ref{secexpl}, the value of the sector angle $\thetan$ given in (\ref{thtval}) is tuned so that these proposed point-to-point connections have a slow power-law decay probability of being accepted. 
This tuning of an acceptance probability to within a desired range is reminiscent of mixing time analysis of the Metropolis algorithm, where a rapid approach to the invariant measure is often achieved by tuning the acceptance rate to be $\Theta(1)$. It would be interesting to try to determine the mixing time of such procedures as $\res$ for sampling the measure $P \big( \cdot \big\vert \areacon \big)$, and to show that choices of $\thetan$ of the order of $n^{-1/3 + o(1)}$ achieve an optimal rate of convergence. 
\end{subsection}
\begin{subsection}{Rates of decay for the probability of small maximum local deviation}
The decay rate appearing in Theorem \ref{thmmlrubd} is almost certainly not optimal. The appearance of the term $1/13$ is a consequence of our controlling such quantities as ${\rm MFL}$ from above for the purpose of understanding the action of $\res$. It is natural to try to improve this decay rate. In this regard, it would be sensible to weaken the definition of the space $\mathcal{G}$  of good configurations for the action of $\res$ that appears in Section \ref{secthmpf}. After all, one long facet in the convex boundary only really prevents us from analysing the action of $\res$ in a neighbourhood of that facet. An alternative is to 
control from above not ${\rm MFL}$, but only the sum of the lengths of all long facets. Similar comments apply to controlling angular gaps in the circuit regeneration structure. Theorems \ref{thmmaxrgcl} and \ref{thmmflbd} would again provide the needed control, with some additional effort.

The polynomial decay rate in Theorem \ref{thmmflubd} is also an artefact of the method of proof. It would not be difficult to obtain a stretched exponential bound: each instance of local roughness of at least $c n^{1/3} \big( \log n \big)^{2/3}$ being realized by a point whose argument-value lies in the interval of argument-values given by a convex boundary facet of length $o \big( n^{2/3} ( \log n )^{1/3} \big)$ entails a polynomial cost, as we saw in the proof of Theorem \ref{thmmflubd}. These costs may be multiplied if there are many instances of local roughness taking at least this value, and the proof of Theorem \ref{thmmlrubd} may be adapted to show that indeed there are.
\end{subsection}
\end{section}
\bibliographystyle{plain}
\bibliography{mlrbib}
\end{document}